\documentclass[10pt]{article}
\usepackage{latexsym}
\usepackage[all]{xy}

\usepackage{amsmath}
\usepackage{amssymb}
\usepackage{amsfonts}

\newtheorem{thm}{Theorem}[section]
\newtheorem{prop}[thm]{Proposition}
\newtheorem{lem}[thm]{Lemma}
\newtheorem{rmk}[thm]{Remark}
\newtheorem{df}[thm]{Definition}
\newtheorem{cor}[thm]{Corollary}

\newcommand{\s}{\infty}

\newcommand{\rh}{\mathbb{R}\underline{\mathrm{Hom}}}
\newcommand{\z}{\mathbb{Z}}

\begin{document}

\title{\textbf{Derived Azumaya algebras and generators for twisted derived categories}}
\bigskip
\bigskip

\author{\bigskip\\
Bertrand To\"en \\
\small{Institut de Math\'ematiques et de Mod\'elisation de Montpellier}\\
\small{Universit\'e de Montpellier 2}\\
\small{Case Courrier 051}\\
\small{Place Eug\`ene Bataillon}\\
\small{34095 Montpellier Cedex}\\
\small{France}\\
\small{e-mail: btoen@math.univ-montp2.fr}}

\bigskip

\date{December 2011}

\maketitle

\begin{abstract}
We introduce a notion of derived Azumaya algebras over ring and schemes generalizing
the notion of Azumaya algebras of \cite{gr}. We prove that 
any such algebra $B$ on a scheme $X$ provides a class $\phi(B)$ in $H^{1}_{et}(X,\mathbb{Z})\times
H^{2}_{et}(X,\mathbb{G}_{m})$. We prove that for $X$ a quasi-compact and quasi-separated
scheme $\phi$ defines a bijective correspondence, and in particular that any 
class in $H^{2}_{et}(X,\mathbb{G}_{m})$, torsion or not, can be represented by 
a derived Azumaya algebra on $X$. Our result is a consequence of a more general 
theorem about the existence of compact generators in \emph{twisted derived categories, 
with coefficients in any local system of reasonable dg-categories}, generalizing 
the well known existence of compact generators in derived categories of quasi-coherent sheaves
of \cite{bv}. A huge part of this paper concerns the  treatment of twisted
derived categories, as well as the proof that the existence of compact generator locally for the
fppf topology implies the existence of a global compact generator. We present 
explicit examples of derived Azumaya algebras that are not represented by 
classical Azumaya algebras, as well as applications of our main result
to the localization for twisted algebraic K-theory and to the
stability of saturated dg-categories by direct push-forwards along smooth and proper maps.
\end{abstract}

\tableofcontents

\section*{Introduction}

In his seminal paper \cite{gr}, A. Grothendieck studied Azumaya algebras over schemes, which are
locally free sheaves of $\mathcal{O}$-algebras $A$, such that 
the natual morphism
$$A\otimes_{\mathcal{O}}A^{op} \longrightarrow \underline{End}_{\mathcal{O}}(A,A)$$
is an isomorphism. He constructed, for any such  algebra $A$ on a scheme $X$, a cohomology class
$\gamma(A) \in H^{2}_{et}(X,\mathbb{G}_{m})$. The class $\gamma(A)$ is always a torsion element, 
and Grothendieck asked the question wether all torsion classes in $H^{2}_{et}(X,\mathbb{G}_{m})$
arise by this construction. The importance of this question lies in the fact that Azumaya algebras 
have description purely in terms of vector bundles on $X$, never mentioning the \'etale
topology. A positive answer to his question would thus provide a description of the torsion part of the 
group $H^{2}_{et}(X,\mathbb{G}_{m})$ uniquely in terms envolving vector bundles on the scheme $X$. 

Since then, this question has been studied by several authors and many results, positive or
negative, have been proven. The most important of these results is certainly the theorem 
of O. Gabber (in \cite{gab}, and \cite{dj}), stating that any torsion class in $H^{2}_{et}(X,\mathbb{G}_{m})$ 
can be realized by an Azumaya algebra when $X$ is a (quasi-compact and separated) scheme 
admitting a ample  line bundle (e.g. a quasi-projective scheme over some ring $k$).
On the side of negative results, there are known 
examples of non-separated schemes 
for which not all torsion classes in $H^{2}_{et}(X,\mathbb{G}_{m})$ come from Azumaya algebras
(see \cite[Cor. 3.11]{kr}). There are also schemes $X$ with non-torsion $H^{2}_{et}(X,\mathbb{G}_{m})$, 
such as the famous example of Mumford with $X$ a normal complex algebraic surface, for which 
non-torsion elements persist locally for the Zariski topology on $X$ (see \cite{gr}). Non-reduced
scheme, such as $X[\epsilon]$ for $X$ a smooth variety, also provide plenty of examples
of schemes with non-torsion $H^{2}_{et}(X,\mathbb{G}_{m})$ (e.g. first-order deformations
of the trivial Azumaya algebra, which are classified by $H^{2}(X,\mathcal{O}_{X})$). 

The main purpose of this work is to continue further the study of the cohomology group
$H^{2}_{et}(X,\mathbb{G}_{m})$, and more particularly of the part which is not 
represented by Azumaya algebras, including the non-torsion part. For this, we introduce 
a notion of \emph{derived Azumaya algebras}, which are differential graded analogs 
of Azumaya algebras, which, as a first rought approximation, can be thought as 
\emph{Azumaya algebra objects in the perfect derived category}. 
By definition, a derived Azumaya algebra over a scheme $X$, 
is a sheaf of associative $\mathcal{O}_{X}$-dg-algebras $B$, satisfying the following two conditions:
\begin{enumerate}
\item The underlying complex of $\mathcal{O}_{X}$-modules is a perfect complex
on $X$, and we have $B\otimes_{\mathcal{O}_{X}}^{\mathbb{L}}k(x)\neq 0$ for all point $x\in X$
(which is equivalent to say that the restriction of $B$ over any affine open sub-scheme
$Spec\, A \subset X$, provides a compact generator of the derived category $D(A)$, see remark \ref{r0}).
\item The natural morphism of complexes of $\mathcal{O}_{X}$-modules 
$$B\otimes_{\mathcal{O}_{X}}^{\mathbb{L}}B^{op} \longrightarrow \mathbb{R}\underline{\mathcal{H}om}_{\mathcal{O}_{X}}(B,B)$$
is a quasi-isomorphism. 
\end{enumerate}
Our first main result state that derived Azumaya algebras are
locally trivial for the \'etale topology (see proposition \ref{p4}), 
and that the Morita equivalence classes of derived Azumaya algebras
on a scheme $X$ embedds into $H^{1}_{et}(X,\mathbb{Z}) \times H^{2}_{et}(X,\mathbb{G}_{m})$
(see corollary \ref{c5}). Any derived Azumaya algebra $B$ on a scheme (or more generally stack) $X$ 
then possesses a characteristic class 
$\phi(B)\in H^{2}_{et}(X,\mathbb{G}_{m})$. Our second main result is the following theorem.

\begin{thm}\label{ti1}{\emph{(See corollary \ref{c6})}}
If $X$ is a quasi-compact and quasi-separated scheme, then any class in $H^{2}_{et}(X,\mathbb{G}_{m})$
is of the form $\phi(B)$ for some derived Azumaya algebra $B$ over $X$. 

More precisely, there exists a bijection
$$\phi : dBr(X) \simeq H^{1}_{et}(X,\mathbb{Z})\times H^{2}_{et}(X,\mathbb{G}_{m}),$$
where $dBr(X)$ is the group of Morita equivalence classes of derived Azumaya algebras over $X$.
\end{thm}
 
As a first consequence we get interesting derived Azumaya algebras in some concrete examples. 
To start with we can take $X$ the normal surface considered by Mumford, and $A$ its local 
ring at the singular point. It is a normal local $\mathbb{C}$-algebra of dimension $2$, 
with a very big $H^{2}_{et}(Spec\, A, \mathbb{G}_{m})$. For any class $\alpha \in H^{2}_{et}(Spec\, A, \mathbb{G}_{m})$, torsion or not, we do get by the previous theorem a derived Azumaya algebra
$B_{\alpha}$ over $A$ representing this class. The dg-algebra $B$ is Morita equivalent to 
an actual (non derived) Azumaya algebra if and only if $\alpha$ is torsion. Another example
is to start with $X$ a smooth and projective complex manifold, with $H^{2}(X,\mathcal{O}_{X})\neq 0$. 
Any non-zero element $\alpha \in H^{2}(X,\mathcal{O}_{X})$, can be interpreted as a class
in $H^{2}_{et}(X[\epsilon],\mathbb{G}_{m})$ restricting to zero on $X$, and this provides
by our theorem a derived Azumaya algebra $B_{\alpha}$ over $X[\epsilon]$, which is
a first order deformation of $\mathcal{O}_{X}$, up to Morita equivalence. In both cases, 
we obtain examples of derived Azumaya algebras which are not equivalent to 
Azumaya algebras, but represent interesting, or at least natural, cohomology classes. \\

The proof of theorem \ref{ti1} is done in two independant steps. The first step is the
construction of the map $\phi$, that associates a cohomology class to a given 
derived Azumaya algebra. For this, we start with a local study of derived Azumaya 
algebras over a base commutative ring, and we prove that any derived Azumaya algebra $B$ is locally 
Morita equivalent, for the \'etale topology, to the trivial derived Azumaya algebra $\mathcal{O}$
(see proposition \ref{p4}). This local triviality property 
is deduced from a deformation theory argument: we first show that any 
derived Azumaya algebras is trivial over an algebraically closed field, and then 
extends this result over an arbitrary stricly Henselian ring by proving that the
stack of trivialization of a given derived Azumaya algebra is algebraic and smooth. This
last step uses the existence of a reasonable moduli stack of dg-modules over dg-algebras and
is based on results from \cite{tova}. It is the only non-purely algebraic part of the
local theory. 
Once the local triviality is known, the construction of the class $\phi$ is rather formal.
We compute the self-equivalence group of the trivial derived Azumaya algebra, and find
that it is the group stack $\mathbb{Z} \times K(\mathbb{G}_{m},1)$. As a consequence the
stack associated to the prestack of derived Azumaya algebras is $K(\mathbb{Z},1)\times K(\mathbb{G}_{m},2)$, 
and it follows formally that any such dg-algebra $B$ over some scheme $X$ provides a class $\phi(B)$ in 
$H^{1}_{et}(X,\mathbb{Z})\times H^{2}_{et}(X,\mathbb{G}_{m})$. The injectivity of the
construction $\phi(B)$ is formal in view of the definition of Morita equivalences between 
derived Azumaya algebras. 

The second step in the proof of theorem \ref{ti1} consists in proving that the map $\phi$ is surjective, which 
is, by very far, the hardest part of the result. The general strategy is 
to deduce the surjectivity from a more general statement stating the existence of  
compact generators in \emph{twisted derived categories}. Our approach is then very similar to 
the approach presented in \cite{dj}, re-interpreting the existence of Azumaya algebras as
the existence of a certain object in a category of twisted sheaves. Our approach has to be done
at the level of derived categories as well as with a much more general notion of \emph{twisted sheaves}. 
We introduce in our second section the notion of \emph{locally presentable dg-categories over schemes}, and
study the stack of those. These are algebraic families of dg-categories, parametrized by some schemes, 
satisfying some set theory conditions
(being locally presentable, or equivalently \emph{well generated} in the sense
of triangulated category, see \cite{ne}), as well as some descent condition with respect to the
flat topology. One key statement is that these locally presentable dg-categories form 
a stack in $\s$-categories $\mathbb{D}g^{lp,desc}$ (i.e. locally presentable dg-categories can  
be glued together with respect to the flat topology, see theorem \ref{t1}). To each scheme $X$, and to any locally presentable dg-category $\alpha$ over
$X$, we can construct its global section $L_{\alpha}(X)$, which is a dg-category obtained by 
integration of $\alpha$ over all affine open sub-schemes in $X$. By definition, 
$L_{\alpha}(X)$ is called the \emph{twisted derived category of $X$ with coefficients in $\alpha$}. 
On the other hand, we construct the stack (again of $\s$-categories) of (quasi-coherent) 
dg-algebras $\mathbb{D}g\mathbb{A}lg$, 
parametrizing quasi-coherent dg-algebras over schemes. There exists a morphism of stacks
$$\phi : \mathbb{D}g\mathbb{A}lg \longrightarrow \mathbb{D}g^{lp,desc}$$
simply obtained by sending a given dg-algebra to the corresponding dg-category of dg-modules. The morphism
$\phi$ provides, for any dg-algebra $B$ over some scheme $X$, a locally presentable dg-category 
$\phi(B)$ over $X$, and thus a corresponding twisted derived category $L_{\phi(B)}(X)$, which 
is nothing else that the dg-category of sheaves of quasi-coherent $B$-dg-modules over $X$. 
The main observation is that we can recognize the locally presentable dg-categories $\alpha$ 
over $X$ of the form $\phi(B)$ for some dg-algebra $B$, by a simple criterion concerning 
the existence of an object $E\in L_{\alpha}(X)$ whose restriction to each 
open affine $U=Spec\, A \subset X$ is a compact generator in $L_{\alpha}(U)$. The major 
result of this work is then the following existence theorem, which is a far reaching 
generalization of the main result of \cite{bv}. 

\begin{thm}\label{ti2}{\emph{(See theorem \ref{t2})}}
Let $X$ be a quasi-compact and quasi-separated scheme, and $\alpha$ a locally 
presentable dg-category over $X$. Assume that there exists 
a fppf covering $X' \longrightarrow X$ such that $L_{\alpha}(X')$ admits 
a compact generator. Then, there exists a compact generator
$E \in L_{\alpha}(X)$ which also a compact local generator.
\end{thm}

The theorem \ref{ti1} is now a direct consequence of the previous result. Indeed, 
the group of self equivalences of the trivial locally presentable dg-category $\mathbf{1}$ on $X$
(which is the family of dg-categories of quasi-coherent complexes over open sub-schemes of $X$), 
is equivalent to $\mathbb{Z}\times K(\mathbb{G}_{m},1)$
(shifts and tensoring with a line bundle are the only $A$-linear self equivalences of the 
dg-category of $A$-dg-modules, for $A$ a commutative ring). Therefore, 
a class $\alpha \in H^{1}_{et}(X,\mathbb{Z})\times H^{2}_{et}(X,\mathbb{G}_{m})$
provides a twisted form of $\mathbf{1}$, that is a locally presentable 
dg-category $\alpha$ over $X$ locally equivalent, for the $fppf$ or \'etale topology, 
to $\mathbf{1}$. By the theorem \ref{ti2}, $L_{\alpha}(X)$ possesses a compact local generator, 
and thus $\alpha$ is of the form $\phi(B)$ for a certain dg-algebra $B$ over $X$, which 
is necessarly a derived Azumaya algebras (because it is so locally). We would like to mention, however, 
that theorem \ref{ti2} seems to us somehow more important than its consequence theorem \ref{ti1}. Indeed, 
there exists many interesting consequences of theorem \ref{ti2}, outside of the 
general question of existence of Azumaya algebras. Some will be given in our last section, and
includes for instance a localization sequence for twisted K-theory, as well as 
construction of smooth and proper dg-categories by direct images of smooth and proper locally 
presentable dg-categories over smooth and proper schemes. \\

\textbf{Related works:}
To finish this introduction, other approaches already exist in order to 
describe the group $H^{2}_{et}(X,\mathbb{G}_{m})$ and particularly its non-torsion part. 
In \cite{js}, the notion of central separable algebras, as well as the \emph{big Brauer group}
is studied. The authors show that the whole group $H^{2}_{et}(X,\mathbb{G}_{m})$
is in bijection with certain equivalence classes of central separable algebras, for 
any quasi-compact Artin stack with an affine diagonal. The present work is in some sense
orthogonal to the approach \cite{js}, as our notion of derived Azumaya algebras
stay as close as possible to the usual notion of Azumaya algebras, whereas the notion
of central separable algebras involves non-unital as well as infinite dimensional 
algebras. Also, the local triviality for the
\'etale topology is taken as part of the definition in the big Brauer group (see \cite[Def. 2.1]{js}), and 
therefore
the results of \cite{js} do not truly make a correspondence between 
$H^{2}_{et}(X,\mathbb{G}_{m})$ and a purely quasi-coherent structure on $X$ which would not mention 
the \'etale topology, contrary to our correspondence \ref{ti1} for which $dBr(X)$ is
defined purely in terms of the theory of perfect complexes on $X$. Another difference, 
it is unclear to me how modules over central separable algebras behave, and how far they form 
a category equivalent to the category of twisted sheaves. One important feature of derived 
Azumaya algebras is that the corresponding derived category of dg-modules gives back the
twisted derived category. There are, however, some similarities between our approach and the
work \cite{js}. The key statement we had to prove is the existence of compact object
with certain nice local properties in the twisted derived category. In the same way, 
the authors \cite{js} deduces their result from the existence of a twisted coherent
sheaf with nice local property. The 
existence of this coherent sheaf is in turn a formal consequence of the extension property 
of coherent sheaves on nice Artin stacks, whereas the existence of a compact local generator
requires extensions properties of compact objects along open immersions (propsition \ref{p8}), as well 
as certain tricky descent statements of compact generators along \'etale and $fppf$ coverings
(propositions \ref{p9}, \ref{p10}).

We should also mention the preprint \cite{brs}, in which the authors define and study 
a topological analog of our notion of derived Azumaya algebras. The definitions of \cite{brs}
has been compared with the one of the present paper in \cite{jo}. A slight adaptation of
our proof of the local triviality for the \'etale topology (proposition \ref{p4}) 
also proves the local triviality for the \'etale topology over connective ring spectra
(this adaptation requires to use \emph{brave new algebraic geometry}, as this is done in \cite[\S 2.4]{hagII}
for instance). Using these techniques, it is for instance possible to prove that the derived Brauer group 
of the sphere spectrum is trivial.

Finally, a recent result concerning the existence of compact generators over algebraic spaces has been 
proved in \cite[Thm. 1.5.10]{lu3}. It can be used in order to extend the theorem \ref{ti2} from 
schemes to algebraic spaces, at least if the \'etale topology is used instead of
the fppf topology.

\bigskip

\textbf{Acknowledgements:} We thank B. Richter for pointing out to us the works \cite{brs} and \cite{jo}. 
The notion of derived Azumaya algebras presented in this work as appeared
during an email conversation with J. Lurie in 2004 (as far as I remember). 
During this correspondence, I have been 
informed of the work of M. Lieblich (see \cite{li}), who proposed a notion of \emph{generalized
Azumaya algebras}, which is rather close, but different, to our notion of derived Azumaya algebra. 
The observation that the existence of generalized Azumaya algebras 
is equivalent to the existence of perfect complexes with nice local properties on $\mathbb{G}_{m}$-gerbes
is already contained in his work. Also, the notion of locally presentable dg-category with descent
of definition \ref{d4} has emmerged during the same email correspondence, and I thank J. Lurie for 
his insights concerning this notion. 

\bigskip

\section{The local theory}

For a simplicial commutative ring $A$ we note $N(A)$ its normalized complex.
The shuffle maps of \cite{ss2} endow the complex $N(A)$ with the structure
of a commutative dg-algebra (over $\mathbb{Z}$). The category of unbounded 
dg-modules over $N(A)$ will be denoted by $A-dg-mod$, and its objects will be called
\emph{$A$-dg-modules}. The category $A-dg-mod$ is a symmetric mono\"\i dal 
model category, for the mono\"\i dal structure induced by the tensor product of dg-modules over $N(A)$, and
for which equivalences are quasi-isomorphisms and fibrations are epimorphisms. The mono\"\i dal 
model category $N(A)-dg-mod$ satisfies moreover the condition \cite{ss1}, and therefore there 
induced model  structure on
mono\"\i ds in $N(A)-dg-mod$ exists (equivalences and fibrations are defined on the underlying objects in $A-dg-
mod$). 
mono\"\i ds in $N(A)-dg-mod$ will be called \emph{$A$-dg-algebras}, and their category is noted $A-dg-alg$.  

The homotopy category of $A$-dg-modules is denoted by $D(A)$, endowed with its natural triangulated structure.
Recall that compact objects in $D(A)$ are precisely the perfect (or dualizable) $A$-dg-modules (see  
\cite{tova}), or equivalently the retracts of finite cell $A$-dg-modules. The category $D(A)$ is 
moreover 
a closed symmetric mono\"\i dal category for the tensor product of $A$-dg-modules. We will note 
$\otimes^{\mathbb{L}}_{A}$
its mono\"\i dal structure and $\mathbb{R}\underline{Hom}_{A}$ the corresponding internal Hom objects. The 
complexes underlying $\mathbb{R}\underline{Hom}_{A}$, by forgetting the $A$-dg-module structure, will be 
denoted $\mathbb{R}Hom_{A}$.  
More generally, if $B$ is a $A$-dg-algebra (or even a $A$-dg-category), the derived category 
$D(B)$ of $B$-dg-modules is naturally enriched over $D(A)$. We denote by 
$\mathbb{R}\underline{Hom}_{B}$ the corresponding Hom's with values in $D(A)$ (note that the
base commutative simplicial ring $A$ is then ambiguous in this notation, but in each situation the base
should be clear from the context). When $B$ is itself commutative, $\mathbb{R}\underline{Hom}_{B}$ 
is then itself a $B$-dg-module, and the two notations agree. 

\subsection{Derived Azumaya algebras over simplicial rings}

We start by the definition of derived Azumaya algebras over simplicial rings. This notion will be later 
generalized to a sheaf-like setting in order to consider derived Azumaya over derived schemes and more 
generally derived stacks. In this section we will concentrate on the basic definition as well as its formal 
properties. 

\begin{df}\label{d1}
\emph{Let $A$ be a commutative simplicial ring. A} derived Azumaya algebra over $A$
(a deraz $A$-algebra \emph{for short) is 
a $A$-dg-algebra $B$ satisfying the following two conditions.}
\begin{enumerate}
\item[$\mathbf{(Az-1)}$] \emph{The underlying $A$-dg-module is a compact generator of the triangulated category $D(A)$}.
\item[$\mathbf{(Az-2)}$] \emph{The natural morphism in $D(A)$}
$$B\otimes_{A}^{\mathbb{L}}B^{op}\longrightarrow \mathbb{R}\underline{Hom}_{A}(B,B),$$
\emph{induced by} 
$$(b,b') \mapsto (c\mapsto bcb')$$
\emph{is an isomorphism}. 
\end{enumerate}
\end{df}

\begin{rmk}\label{r0}
\begin{enumerate}
\item
\emph{When $A$ is a nonsimplicial ring, considered as a constant simplicial ring, then 
Azumaya $A$-algebras are special cases of deraz $A$-algebras: they correspond, up to quasi-isomorphism, 
to deraz $A$-algebras $B$ whose underlying $A$-dg-module is flat and concentrated in degree $0$. Indeed, 
any such deraz $A$-algebra is quasi-isomorphic to a non-dg $A$-algebra $B$, which is projective and of finite type as an $A$-module (condition $\mathbf{(Az-1)}$), and moreover satisfies $B\otimes_{A}B^{op} \simeq \underline{Hom}_{A}(B,B)$ (condition$\mathbf{(Az-2)}$). Our notion of derived Azumaya algebras over $A$ is therefore a generalization of the usual notion of Azumaya algebras. We will see later that there are examples of deraz $A$-algebras non-quasi-isomorphic (and even not Morita equivalent) to underived ones.}
\item 
\emph{The condition $\mathbf{(Az-1)}$ can be checked using the following criterion: a compact
object $E\in D(A)$ is a compact generator if and only if for any prime ideal $p$ of $\pi_{0}(A)$, 
with residue field $A \rightarrow k(p)$, we have $E\otimes_{A}^{\mathbb{L}}k(p)\neq 0$ in $D(k(p))$. 
Indeed, the condition is equivalent to state that the support of $E$ is the whole derived affine
scheme $\mathbb{R}\underline{Spec}\, A$. By \cite[Lem. 3.14]{tho} this implies that $A$ belongs to 
the thick triangulated sub-category generated by $E$, and thus that $E$ is a compact generator
of $D(A)$.}
\end{enumerate}
\end{rmk}

For any morphism of simplicial commutative rings $A \rightarrow A'$, there are base change functors
$$A'\otimes_{A}^{\mathbb{L}} - : D(A) \longrightarrow D(A'),$$
$$A'\otimes_{A}^{\mathbb{L}} - : Ho(A-dg-alg) \longrightarrow Ho(A'-dg-alg).$$
These base change functors will also be denoted by 
$$B \mapsto B_{A'}:=A'\otimes_{A}^{\mathbb{L}}B.$$
For the next proposition, remind that $A \rightarrow A'$ is faithfully flat if 
$\pi_{0}(A) \rightarrow \pi_{0}(A')$ is a faithfully flat morphism of rings, and if 
$$\pi_{*}(A)\otimes_{\pi_{0}(A)}\pi_{0}(A') \simeq \pi_{*}(A')$$
(see \cite[\S 2.2.2]{hagII} for other characterizations).

\begin{prop}\label{p1}
Let $A \rightarrow A'$ be a morphism of simplicial commutative rings. 
\begin{enumerate}
\item If $B$ is deraz $A$-algebra, then 
$B_{A'}$ is a deraz $A'$-algebra. 
\item If $A'$ is faithfully flat over $A$, and if $B \in Ho(A-dg-alg)$ is such that 
$B_{A'}$ is a deraz $A'$-algebra, then 
$B$ is a deraz $A$-algebra. 
\end{enumerate}
\end{prop}

\textit{Proof:} $(1)$ Let $B$ be a deraz $A$-algebra.
Let $E \in D(A')$, and assume that 
$$\mathbb{R}Hom_{A'}(B_{A'},E)=0.$$
By adjunction, we have
$$\mathbb{R}Hom_{A}(B,E)=0,$$
where $E$ is considered as a $A$-dg-module throught the morphism $A \rightarrow A'$.
As $B$ is a compact generator of $D(A)$ we have $E=0$ in $D(A)$, which implies that $E=0$ also in $D(A')$. 
Therefore $B_{A'}$ is compact generator of $D(A')$ and satisfies condition $\mathbf{(Az-1)}$ of \ref{d1}. 

There exists a natural isomorphism in $D(A')$
$$B_{A'}\otimes_{A'}^{\mathbb{L}}(B_{A'})^{op} \simeq A'\otimes_{A}^{\mathbb{L}}
(B\otimes_{A}^{\mathbb{L}}B^{op}).$$
Also, as $B$ is a compact $A$-dg-module, there exists a natural isomorphism in $D(A')$
$$A'\otimes_{A}^{\mathbb{L}}\mathbb{R}\underline{Hom}_{A}(B,B) \simeq
\mathbb{R}\underline{Hom}_{A'-dg-mod}(B_{A'},B_{A'}).$$
Under these isomorphisms, the morphism of $A'$-dg-modules
$$B_{A'}\otimes_{A'}^{\mathbb{L}}(B_{A'})^{op}\longrightarrow \mathbb{R}\underline{Hom}_{A'}(B',B'),$$
is the image by the functor $A'\otimes_{A}^{\mathbb{L}}-$ of the morphism 
$$B\otimes_{A}^{\mathbb{L}}B^{op}\longrightarrow \mathbb{R}\underline{Hom}_{A}(B,B).$$
We thus see that $B_{A'}$ satisfies condition $\mathbf{(Az-2)}$ of definition \ref{d1}. \\

$(2)$ Let $B$ be a $A$-dg-algebra such that $B_{A'}$ is a deraz $A'$-algebra. 
By \cite[Prop. 1.3.7.4]{hagII} 
we know that an object $E\in D(A)$ is perfect if and only if $A'\otimes_{A}^{\mathbb{L}}E$ is so. 
As perfect $A$-dg-modules are exactly the compact objects in $D(A)$, we see that $B$ is a compact object in 
$D(A)$. Assume now that $E\in D(A)$ is such that $\mathbb{R}\underline{Hom}_{A}(B,E)=0$. As $B$ is compact, we 
have 
$$A'\otimes_{A}^{\mathbb{L}}\mathbb{R}\underline{Hom}_{A}(B,E)\simeq 
\mathbb{R}\underline{Hom}_{A'}(B_{A'},E_{A'}).$$ 
We thus have $\mathbb{R}\underline{Hom}_{A'}(B_{A'},E_{A'})=0$, and as 
$B_{A'}$ is a compact generator of $D(A')$ we must have $E_{A'}=0$. As $A'$ is faithfully flat over $A$, we also 
have $E=0$. This shows that $B$ is a compact generator of $D(A)$, and thus that $B$ satisfies condition 
$\mathbf{(Az-1)}$ of definition \ref{d1}. 

Now, as we have seen during the proof of the first point, the morphism
$$B_{A'}\otimes_{A'}^{\mathbb{L}}(B_{A'})^{op}\longrightarrow \mathbb{R}\underline{Hom}_{A'-dg-mod}(B',B')$$
is obtained from 
$$B\otimes_{A}^{\mathbb{L}}B^{op}\longrightarrow \mathbb{R}\underline{Hom}_{A}(B,B)$$
by base change from $A$ to $A'$. As $A'$ is a fully 
faithfully flat $A$-dg-algebra, and as by hypothesis the first 
of these morphisms is an isomorphism, we see that $B$ must satisfy condition $\mathbf{(Az-2)}$ of definition 
\ref{d1}, and therefore is a deraz $A$-algebra. 
\hfill $\Box$ \\

We are now going to see that the property of being a derived Azumaya algebra is stable by Morita
equivalences. For this we will show that the conditions $\mathbf{(Az-1)}$ and $\mathbf{(Az-2)}$ 
can be expressed by 
stating that $B$ is an invertible object in a certain closed symmetric mono\"\i dal category. A formal consequence 
of this fact will be the Morita invariance of the notion of deraz algebras, as well as its stability by 
(derived) tensor products. 

Let $A$ be a simplicial commutative ring.
We let $Hom_{A}$ be the homotopy category of $A$-dg-categories up to derived
Morita equivalences. The precise definition of $Hom_{A}$ is as follows. Its objects are categories enriched
into the symmetric mono\"\i dal category of $A$-dg-modules (which by definition are dg-modules over the
commutative dg-algebra $N(A)$), which we simply call \emph{$A$-dg-categories}. For a $A$-dg-category $T$ we
can form its derived category of $T$-dg-modules, $D(T)$, which is, by definition, the homotopy category of 
left $T$-dg-modules.
A morphism of $A$-dg-categories is simply an enriched functor, and such a functor $f : T \longrightarrow T'$
is a \emph{Morita equivalence} if the induced functor on the level of derived categories of dg-modules
$$f^{*} : D(T') \longrightarrow D(T)$$
is an equivalence of categories. The category $Hom_{A}$ is then obtained from the category of 
$A$-dg-categories and morphisms of $A$-dg-categories by inverting the Morita equivalences. The sets of 
morphisms in $Hom_{A}$ are usually denoted by $[-,-]$, or $[-,-]_{Hom_{A}}$.
According to 
 \cite[Cor. 7.6]{to1} these sets of morphisms  can be described as follows: for $T$ and $T'$ we form
$T\otimes_{A}^{\mathbb{L}}(T')^{op}$, and consider $D_{pspe}(T\otimes_{A}^{\mathbb{L}}(T')^{op})$, the full
sub-category of $D(T\otimes_{A}^{\mathbb{L}}(T')^{op})$ consisting of all $E$ such that for any object 
$x\in T$, $E(x,-)$ is a compact object in $D((T')^{op})$. Then, we have a natural identification
$$[T,T']_{Hom_{A}}\simeq D_{pspe}(T\otimes_{A}^{\mathbb{L}}(T')^{op})/isom,$$
between the set of morphisms from $T$ to $T'$ in the category $Hom_{A}$, and the set of isomorphism 
classes of objects in $D_{pspe}(T\otimes_{A}^{\mathbb{L}}(T')^{op})$. 

The category $Hom_{A}$ is made into a symmetric mono\"\i dal category by the derived tensor product of 
$A$-dg-categories
$$\otimes_{A}^{\mathbb{L}} : Hom_{A} \times Hom_{A} \longrightarrow Hom_{A}.$$
By the results of \cite[Cor. 7.6]{to1} this mono\"\i dal structure is known to be closed. 
As in any symmetric mono\"\i dal category, we can therefore make sense of dualizable, as well as invertible, 
objects in $Hom_{A}$, whose definitions are now breifly recalled. An object $T\in Hom_{A}$ 
is \emph{dualizable} (or \emph{with duals}) if there is $T' \in Hom_{A}$ and a morphism 
$$ev : T \otimes_{A}^{\mathbb{L}}T' \longrightarrow \mathbf{1}=A,$$
such that for any two other objects $U,U'\in Hom_{A}$ the induced map 
$$\xymatrix{
[U,T'\otimes_{A}^{\mathbb{L}}U'] \ar[r]^-{T\otimes_{A}^{\mathbb{L}}} & 
[T\otimes_{A}^{\mathbb{L}}U,T\otimes_{A}^{\mathbb{L}}T'\otimes_{A}^{\mathbb{L}}U'] \ar[r]^-{ev} & 
[T\otimes_{A}^{\mathbb{L}}U,U']}$$
is bijective. This condition is equivalent to ask for the existence of a morphism
$$i : A \longrightarrow T \otimes_{A}^{\mathbb{L}}T',$$
such that the two composite
$$\xymatrix{
T' \ar[r]^-{id_{T'}\otimes_{A}^{\mathbb{L}}i} & T'\otimes_{A}^{\mathbb{L}}T \otimes_{A}^{\mathbb{L}}T'
\ar[r]^-{ev\otimes_{A}^{\mathbb{L}}id_{T'}} & T'} \qquad \xymatrix{
T \ar[r]^-{i\otimes_{A}^{\mathbb{L}}id_{T}} & T\otimes_{A}^{\mathbb{L}}T' \otimes_{A}^{\mathbb{L}}T
\ar[r]^-{id_{T}\otimes_{A}^{\mathbb{L}}ev} & T}$$
are identities. Such a triple $(T,i,ev)$ will be called a \emph{duality datum for $T$}. An important fact about dualizable objects is the uniqueness a duality datum.
If $T$ is dualizable,
with $(T',i,ev)$ and $(T'',i',ev')$ two duality data, then there is a unique isomorphism
$\phi : T' \simeq T''$, intertwining $(i,ev)$ and $(i',ev')$. To finish,
an object $T\in Hom_{A}$ is invertible if there is $T'\in Hom_{A}$ such that $T\otimes_{A}^{\mathbb{L}}T'$ is
isomorphic to the unit $\mathbf{1}=A$. Invertible objects are dualizable, as the choice of an isomorphism
$ev : T\otimes_{A}^{\mathbb{L}}T' \simeq A$ satisfies the condition for $(T',ev^{-1},ev)$ to be a duality 
datum. By uniqueness of these duality datum we see that a dualizable object $T$, with 
duality datum $(T',ev,i)$, is invertible if and only if $ev$ or $i$ are isomorphisms. \\

An $A$-dg-algebra $B$ defines canonically an object $B\in Hom_{A}$, by considering $B$ as 
a $A$-dg-category with a unique object $*$, and $B$ as being the endomorphisms of $*$. This constructions 
defines a functor
$$Ho(A-dg-alg) \longrightarrow Hom_{A},$$
which is compatible with base change and forgetful functors, as well as with the derived tensor
product $\otimes_{A}^{\mathbb{L}}$. 

\begin{df}\label{d2}
\emph{Two deraz $A$-algebras are} Morita equivalent \emph{if their images in $Hom_{A}$ are isomorphic}. 
\end{df}

We will use often the following easy criterion (whose proof is
left to the reader): two $A$-dg-algebras $B$ and $B'$ are isomorphic
in $Hom_{A}$ if and only if there exists a compact generator $E\in D(B')$ whose 
derived endomorphism $A$-dg-algebra $\mathbb{R}\underline{End}_{B'}(E)$ is quasi-isomorphic
to $B$. \\

We now have the following proposition, characterizing deraz $A$-algebras as invertible objects in 
$Hom_{A}$. As a side result we also show that dualizable objects in $Hom_{A}$ are precisely the 
smooth and proper $A$-dg-algebras (as stated without proofs in 
\cite[Prop. 5.4.2]{swisk}).

\begin{prop}\label{p2}
Let $A$ be a simplicial commutative ring. 
\begin{enumerate}
\item An object $T$ in $Hom_{A}$ is dualizable if and only if it is isomorphic to 
(the image of) an $A$-dg-algebra $B$ satisfying the following two conditions.
\begin{enumerate}
\item $B$ is a compact object in $D(A)$ (i.e. $B$ is proper over $A$).
\item $B$ is a compact object in $D(B\otimes_{A}^{\mathbb{L}}B^{op})$ (i.e. $B$ is smooth over $A$).
\end{enumerate}
\item 
An object $T$ in $Hom_{A}$ is invertible if and only if is isomorphic to 
(the image of) a deraz $A$-algebra. 
\end{enumerate}
\end{prop}

\textit{Proof:} $(1)$ We introduce a bigger symetric mono\"\i dal category $\widetilde{Hom}_{A}$ containing 
$Hom_{A}$ as a mono\"\i dal sub-category. The objects of $\widetilde{Hom}_{A}$ are 
$A$-dg-categories. The set of morphisms from $T$ to $T'$ in $\widetilde{Hom}_{A}$ is by definition
the set of isomorphism classes of objects in $D(T\otimes_{A}^{\mathbb{L}}(T')^{op})$. The composition
of morphisms is then induced by the (derived) convolution tensor product construction
$$-\otimes_{T'}^{\mathbb{L}}- : D(T\otimes_{A}^{\mathbb{L}}(T')^{op}) \times
D(T'\otimes_{A}^{\mathbb{L}}(T'')^{op}) \longrightarrow D(T\otimes_{A}^{\mathbb{L}}(T'')^{op}).$$
The reader might worry about the fact that the Hom sets of the category 
$\widetilde{Hom}_{A}$ are not small sets. One obvious solution to this problem is to work 
with universes: if our dg-categories belongs to a universe $\mathbb{U}$, then 
$\widetilde{Hom}_{A}$ will be a $\mathbb{V}$-small category, for $\mathbb{V}$ a universe with
$\mathbb{U}\in \mathbb{V}$. Another solution, avoiding to be forced to believe that universes exist, 
is to check that we are going to use a very small part of the category $\widetilde{Hom}_{A}$, which consists of
all the tensor powers of a given dg-category and its opposite, as well as certain diagonal dg-modules. All of these
obviously live in a small sub-category of $\widetilde{Hom}_{A}$ in which the argument below can be
done. 
 
The category $\widetilde{Hom}_{A}$ is endowed with a symmetric mono\"\i dal structure induced by the
derived tensor product of dg-categories. The nice property of the symmetric mono\"\i dal category 
$\widetilde{Hom}_{A}$ is that all of its objects are dualizable. Indeed, for any $A$-dg-category $T$, 
the $T\otimes_{A}^{\mathbb{L}}T^{op}$ dg-module $(a,b) \mapsto T(b,a)$, provides
an object $T\in D(T\otimes_{A}^{\mathbb{L}}T^{op})$, and therefore two morphisms in $\widetilde{Hom}_{A}$
$$i : A \longrightarrow T\otimes_{A}^{\mathbb{L}}T^{op} \qquad
ev : T\otimes_{A}^{\mathbb{L}}T^{op} \longrightarrow A,$$
that can be checked to be a duality datum in $\widetilde{Hom}_{A}$. Now, $Hom_{A}$ is a symmetric mono\"\i dal 
sub-category $\widetilde{Hom}_{A}$, containing all isomorphisms. Therefore, the unicity of a duality datum 
implies that $T\in Hom_{A}$ is dualizable if and only if the two morphisms $i$ and $ev$ as defined above
in $\widetilde{Hom}_{A}$ lies in the sub-category $Hom_{A}$. That $i$ belongs to $Hom_{A}$ is equivalent to 
state that $T$ is a compact object in $D(T\otimes_{A}^{\mathbb{L}}T^{op})$. In the same way, that 
$ev$ belongs to $Hom_{A}$ implies that $T(b,a)$ is compact in $D(A)$ for any $a,b\in T$. Therefore, in order 
to prove our point $(1)$ it only remains to show that $T$ is moreover isomorphic, in $Hom_{A}$, to 
a $A$-dg-algebra. But this follows from the following lemma.

\begin{lem}\label{l1}
If $T \in Hom_{A}$ is such that 
$T$ is a compact object in $D(T\otimes_{A}^{\mathbb{L}}T^{op})$, then $T$ is isomorphic to 
an $A$-dg-algebra.
\end{lem}

\textit{Proof of the lemma:} We have a natural triangulated functor
$$D(T\otimes_{A}^{\mathbb{L}}T^{op}) \times D(T) \longrightarrow D(T),$$
induced by the composition in $\widetilde{Hom}_{A}$. The object $T$ acts by the identity on 
$D(T)$, and by assumption $T$ lies in the thick triangulated sub-category generated by 
representable objects on $D(T\otimes_{A}^{\mathbb{L}}T^{op})$. These representable objects are given 
by two objects $a,b \in T$, and they act on $D(T)$ by the functor
$$E \mapsto E(a)\otimes_{A}^{\mathbb{L}}T(b,-).$$
This shows that $E=id(E)$ lies in the thick triangulated sub-category of $D(T)$ generated by 
the representable objects $T(b,-)$. As only a finite number of $b$'s are necessary to write 
the object $T$ by retracts, sums, cones and shifts in $D(T\otimes_{A}^{\mathbb{L}}T^{op})$, we see
that $D(T)$ is generated, as a triangulated category, by a finite number of representable objects. 
The direct sums of these representable modules provides a compact generator in $E \in D(T)$, and thus $T$ 
becomes isomorphic in $Hom_{A}$ to $\mathbb{R}\underline{Hom}_{T}(E,E)$, 
the derived endomorphism $A$-dg-algebra of the object $E$.
\hfill $\Box$ \\

$(2)$ An invertible object is a dualizable object. From $(1)$ we see that invertible objects are 
all isomorphic to smooth and proper $A$-dg-algebras. Moreover, the uniqueness of duality data implies that 
such a smooth and proper $A$-dg-algebra $B$ is invertible in $Hom_{A}$ if and only if the morphism
$$ev : B\otimes_{A}^{\mathbb{L}}B^{op} \longrightarrow A,$$
is an isomorphism.
This morphism is itself determined by the diagonal bi-dg-module $B\in D(B\otimes_{A}^{\mathbb{L}}B^{op})$, 
and it is an isomorphism if and only the induced functor
$$\phi : D(B\otimes_{A}^{\mathbb{L}}B^{op}) \longrightarrow D(A)$$
is an equivalence of categories. By construction, this functor sends
$B\otimes_{A}^{\mathbb{L}}B^{op}$ to $B\in D(A)$, and always preserves 
compact objects. Therefore, $\phi$
is fullyfaithful if and only if the morphism
$$B\otimes_{A}^{\mathbb{L}}B^{op}\simeq \mathbb{R}\underline{End}_{B\otimes_{A}^{\mathbb{L}}B^{op}}(B\otimes_{A}^{\mathbb{L}}B^{op})\longrightarrow 
\mathbb{R}\underline{End}_{A}(B,B)$$
is an isomorphism in $D(A)$. In the same way, the essential image of $\phi$ generates
$D(A)$ in the sense of triangulated categories, if and only if $B$ is a compact generator
of $D(A)$. We thus see that $\phi$ is an equivalence of categories if and only 
if $B$ is a deraz $A$-algebra.
\hfill $\Box$ \\

\begin{rmk}\label{r1}
\emph{The proposition \ref{p2} also has an interpretation purely in terms of the category 
$\widetilde{Hom}_{A}$ we have introduced during its proof, at least if this category 
is enhanced into a symmetric mono\"\i dal $2$-category in the usual fashion (the category of morphisms
between $T$ and $T'$ is the category $D(T\otimes_{A}^{\mathbb{L}}(T')^{op})$). Then, it can be proved that 
the an object $T$ in this $2$-category is equivalent to a smooth and proper $A$-dg-algebra if and only if 
it is a fully dualizable object in the sense of \cite{lu2}. This explains the importance of smooth and proper dg-categories in the context of two dimensional field theories. Invertible objects are the same if one considers
$\widetilde{Hom}_{A}$ as a mono\"\i dal category or as a mono\"\i dal $2$-category.}
\end{rmk}

An immediate corollary of the last proposition is the following result (point $(5)$ is 
more a corollary of the proof). 

\begin{cor}\label{c1}
Let $A$ be a simplicial commutative ring.
\begin{enumerate}
\item An $A$-dg-algebra Morita equivalent to a deraz $A$-algebra is itself a deraz $A$-algebra. 
\item If $B$ and $B'$ are two deraz $A$-algebras, then so is $B\otimes_{A}^{\mathbb{L}}B'$.
\item Any deraz $A$-algebra is a smooth and proper $A$-dg-algebra.
\item An $A$-dg-algebra $B$ is a deraz $A$-algebra if and only if there is 
another $A$-dg-algebra $B'$ such that $B\otimes_{A}^{\mathbb{L}}B'$ is Morita equivalent to 
$A$.
\item An $A$-dg-algebra $B$ is a deraz $A$-algebra if and only if the 
$A$-dg-algebra $B\otimes_{A}^{\mathbb{L}}B^{op}$ is Morita equivalent to $A$.
\end{enumerate}
\end{cor}

\subsection{Derived Azumaya algebras over a field}

In this paragraph we fix a base fiel $k$, and will show that all deraz $k$-algebras are 
Morita equivalent to underived Azumaya $k$-algebras. We start by the algebraically closed case.

\begin{prop}\label{p2'}
If $k$ is algebraically closed, then any deraz $k$-algebra is Morita equivalent to 
the trivial $k$-dg-algebra $k$.
\end{prop}

\textit{Proof:} Let $B$ be a deraz $k$-algebra. We need to prove that $B$ is isomorphic, in $Hom_{k}$, to 
the trivial $k$-algebra $k$. 
For this, we let $H^{*}(B)$ be the graded cohomology $k$-algebra of $B$. As we are working over a field, 
we have natural identifications of graded $k$-vector spaces. 
$$H^{*}(B\otimes_{k}B^{op})\simeq H^{*}(B)\otimes_{k}H^{*}(B)^{op}$$
$$H^{*}(\mathbb{R}\underline{Hom}_{k}(B,B))\simeq Hom^{*}_{k}(H^{*}(B),H^{*}(B)).$$
These identifications, and conditions $\mathbf{(Az-1)}$ and $\mathbf{(Az-2)}$ imply that 
$H^{*}(B)$ is itself an Azumaya $k$-algebra, and more precisely that $H^{*}(B)$, considered as 
a graded bi-$H^{*}(B)$-module, induces an equivalence between the categories
of graded modules
$$H^{*}(B)\otimes_{k}H^{*}(B)^{op}-Mod^{gr} \simeq  k-Mod^{gr}.$$
This implies that in particular that $H^{*}(B)$ is a graded projective 
$H^{*}(B)\otimes_{k}H^{*}(B)^{op}$-module, 
and therefore that $H^{*}(B)$ is a semi-simple graded $k$-algebra (i.e. that 
$H^{*}(B)-Mod^{gr}$ is semi-simple $k$-linear abelian category), and thus that any 
graded $H^{*}(B)$-module is graded projective. 

We have a graded functor 
$$\phi : D(B) \longrightarrow H^{*}(B)-Mod^{gr}$$
sending a $B$-dg-module $E$ to its cohomology $H^{*}(E)$. 

\begin{lem}\label{l2}
The functor above $\phi$ is an equivalence of (graded) categories. 
\end{lem}

\textit{Proof:} Let call an object $E\in D(B)$ \emph{graded free} if it is a (possibly infinite) 
sum of shifts of $B$. An object of $D(B)$ will then be called \emph{graded projective} if it is a
retract of a graded free object. It is clear by definition that the functor $\phi$ is fully faithful when 
restricted to graded free objects, and thus also when restricted to graded projective objects. In order to 
show that $\phi$ is fully faithful it is then enough to show that any object in $D(B)$ is graded projective. 

Let $E \in D(B)$, and let $H^{*}(E)$ be the corresponding graded $H^{*}(B)$-module, which is graded projective 
as we saw. Let $p$ be a projector on some graded free graded $H^{*}(B)$-module $P$ with image $H^{*}(E)$. 
If we write 
$$P=\bigoplus_{d}(H^{*}(B)[p]^{(I_{p})}),$$
then $p$ defines a projector $q$ in $D(B)$ on the graded free object
$$Q:=\bigoplus_{d}(B[p]^{(I_{p})}) \in D(B).$$
As $D(B)$ is Karoubian complete (see for instance \cite[Sous-lemme 3]{to4}), we now that this projector splits
$$q : \xymatrix{
Q \ar[r]^-{k} & E' \ar[r]^-{j} & Q,}$$
with $E' \in D(B)$. We claim that $E$ and $E'$ are isomorphic objects in $D(B)$. Indeed, by the choice of 
$p$, there exists a factorisation in $H^{*}(B)-Mod^{gr}$
$$p : \xymatrix{
P \ar[r]^-{r} & H^{*}(E) \ar[r]^-{i} & P.}$$
The images of $k$ and $j$ by $\phi$ provides another such factorisation
$$\phi(q)=p : \xymatrix{
P \ar[r]^-{\phi(k)} & H^{*}(E') \ar[r]^-{\phi(l)} & P.}$$
By the uniqueness of such factorisation, we have that 
the composite morphism 
$$\xymatrix{
H^{*}(E') \ar[r]^-{\phi(l)} & P \ar[r]^-{r} & H^{*}(E)}$$
is an isomorphism.
Now, from the definition of $\phi$, the morphism $r$ lifts uniquily to a morphism 
$r' : Q \longrightarrow E$
with $\phi(r')=r$. We thus get a morphism in $D(B)$
$$\xymatrix{E' \ar[r]^-{j} & Q \ar[r]^-{r'} & E,}$$
whose image by $\phi$ is $r\circ \phi(l)$, which we have seen to be an isomorphism. 
As the functor $\phi$ is conservative we have $E \simeq E'$ in $D(B)$. 

We thus have that all objects in $D(B)$ are graded projective, and therefore that 
$\phi$ is fully faithful. All the objects $H^{*}(B)[p]$ belong to the essential image of $\phi$, and 
therefore so do all graded free $H^{*}(B)$-module. As $\phi$ is fully faithful and as 
$D(B)$ is  Karoubian closed, we deduce that any 
graded projective $H^{*}(B)$-module also belongs to the essential image of $\phi$. But, as all
graded $H^{*}(B)$-module are graded projective we have that $\phi$ is essentially surjective and thus an 
equivalence of categories. 
\hfill $\Box$ \\

We use the previous lemma as follows. Let $E\in H^{*}(B)-Mod^{gr}$ be a simple object, which is
automatically finite dimensional over $k$, as  $H^{*}(B)$ is so.  
Because $k$ is algebraically closed we have 
$$Hom^{0}_{H^{*}(B)-Mod^{gr}}(E,E) \simeq k.$$
Moreover, for any $n\neq 0$, we have
$$Hom^{n}_{H^{*}(B)-Mod^{gr}}(E,E)=Hom^{0}_{H^{*}(B)-Mod^{gr}}(E,E[n])=0,$$
unless $E$ would be isomorphic to $E[n]$ which would a condradiction with the 
finite dimensionality of $E$. Let now $E' \in D(B)$ be an object such that 
$\phi(E')\simeq E$. As the graded functor $\phi$ is an equivalence, we have
$$[E',E']\simeq k \qquad [E',E'[n]]=0 \; if \, n\neq 0.$$
This implies that, the $k$-dg-algebra
$\mathbb{R}\underline{Hom}_{k}(E',E')$ is quasi-isomorphic to $k$, and thus that the pair of adjoint 
functors
$$E \otimes_{k} - : D(k) \longrightarrow D(B) \qquad 
\mathbb{R}\underline{Hom}_{k}(E,-) : D(B) \longrightarrow D(k)$$
defines a full embedding $D(k) \hookrightarrow D(B)$. Therefore, when considered
as a morphism $k \rightarrow B$ in $Hom_{k}$, the object $E$ exhibits $k$ as a retract of $B$. 
By dualizing $B$ we get that $B^{op}$ is a retract of $k$ in $Hom_{k}$. The object 
$B^{op}$ is then the image of a projector $p$ on $k$ in $Hom_{k}$, which corresponds to 
a compact object $P \in D(k)$ with the property that $P\otimes_{k}P \simeq P$. The only possibility is
$P\simeq k$, giving that the projector $p$ is in fact the identity. This implies that $B^{op}$ is isomorphic
to $k$ in $Hom_{k}$, and by duality that $B$ is isomorphic to $k$.  \hfill $\Box$ \\

A traduction of the last proposition is the following more explicit statement.

\begin{cor}\label{c2}
Any derived Azumaya algebra over an algebraically closed field $k$ is quasi-isomorphic 
to a finite dimensional graded matrix algebra $\underline{End}^{*}_{k}(V)$ (for $V$ a graded finite 
dimensional $k$-vector space). In particular, a deraz $k$-algebra is always formal (i.e. quasi-isomorphic to 
its cohomology algebra). 
\end{cor}

We have seen that any deraz $k$-algebra is trivial when $k$ is an algebraically closed field, we will now
show that for any field $k$ the deraz $k$-algebra, up to Morita equivalence, are in one-to-one correspondence
with equivalence classes of nonderived Azumaya algebras. 

For this we introduce underived analog of our symmetric mono\"\i dal category $Hom_{k}$. We let 
$Cm_{k}$ be the category whose objects are $k$-linear categories, and for which the set of morphisms
from $C$ to $C'$ in $Cm_{k}$ is the set of isomorphism classes of $(C\otimes_{k}(C')^{op})$-bimodules $P$, 
such that $P(c,-)$ is a projective $(C')^{op}$-module for any $c\in C$. The category $Cm_{k}$ is
made into a symmetric mono\"\i dal category by means of the tensor product of linear categories over $k$. 
Any $k$-algebra $B$ defines an object $B \in Cm_{k}$, and it is easy to check (by following the same argument
as for our proposition \ref{p1}), that such an algebra is Azumaya over $k$ (in the sense of \cite{gr}), 
if and only if it defines an invertible object in $Cm_{k}$. Moreover, equivalence classes of 
Azumaya algebras over $k$ are in one-to-one correspondence with isomorphism classes of 
invertible objects in $Cm_{k}$.

There is a natural functor
$$Cm_{k} \longrightarrow Hom_{k}$$
which consists of considering a $k$-linear category as a $k$-dg-category is an obvious way (i.e. concentrated in degree $0$). This functor is a symmetric mono\"\i dal functor (because $k$ is a field here), and thus
induces a map on isomorphism classes of invertible objects. 

\begin{prop}\label{p3}
The above functor 
$$Cm_{k} \longrightarrow Hom_{k}$$
induces a bijection between isomorphism classes of invertible objects. In other words, the natural inclusion 
from $k$-algebras to $k$-dg-algebras induces a one-to-one correspondence between the set of equivalent classes 
of Azumaya $k$-algebras and the set of deraz $k$-algebras up to Morita equivalence.
\end{prop}

\textit{Proof:} We need to prove two statements:
\begin{enumerate}
\item If $B$ and $B'$ are two Azumaya $k$-algebras which are isomorphic in $Hom_{k}$, then there also are
isomorphic in $Cm_{k}$.
\item Any deraz $k$-algebra is isomorphic in $Hom_{k}$ to an Azumaya $k$-algebra. 
\end{enumerate}

We start by proving $(1)$. Let $B$ and $B'$ two Azumaya $k$-algebras which are isomorphic
in $Hom_{k}$. We set $C:=B\otimes_{k}(B')^{op}$. We define two stacks $Eq_{0}(B,B')$ and
$Eq(B,B')$ over $Spec\, k$ (we will work with the fppf topology, see however 
\cite{to3}). The stack $Eq_{0}(B,B')$ sends a commutative $k$-algebra 
$A$ to the groupo\"\i d of $C_{A}$-modules $M$ whose corresponding functor
$$\begin{array}{ccc}
B^{op}_{A}-Mod & \longrightarrow & (B_{A}')^{op}-Mod \\
  N & \mapsto & N\otimes_{B}M 
\end{array}$$
is an equivalence of categories. In the same way, the stack 
$Eq(B,B')$ sends a commutative $k$-algebra $A$ to the groupo\"\i d of objects
$M \in D(C_{A})$ whose corresponding functor
$$\begin{array}{ccc}
D(B^{op}_{A}) & \longrightarrow & D((B_{A}')^{op}) \\
  N & \mapsto & N\otimes^{\mathbb{L}}_{B}M 
\end{array}$$
is an equivalence of categories (see \cite[Cor. 3.21]{tova} for a justification that this is
indeed a stack).
The natural inclusion 
$C_{A}-Mod \longrightarrow D(C_{A})$ defines a morphism of stacks
$$f : Eq_{0}(B,B') \longrightarrow Eq(B,B').$$
These stacks have natural actions by the group stacks
$$G_{0}:=Eq_{0}(B,B) \qquad G:=Eq(B,B),$$
and the fact that everything becomes trivial over $\overline{k}$ implies that 
$Eq_{0}(B,B')$ and $Eq(B,B')$ are torsors, for the fppf topology, 
over $Eq_{0}(B,B)$ and $Eq(B,B)$ respectively. 
Now, as $B$ is an invertible object in $Cm_{k}$, we have natural equivalences of group stacks
$$G_{0}\simeq Eq_{0}(k,k) \qquad G\simeq Eq(k,k).$$
These group stacks can be computed explicitely, and we have
$$G_{0}\simeq K(\mathbb{G}_{m},1) \qquad G\simeq \mathbb{Z}\times K(\mathbb{G}_{m},1).$$
The two torsors $Eq_{0}(B,B')$ and $Eq(B,B')$ are thus represented by 
cohomology classes 
$$x_{0} \in H^{1}_{fppf}(Spec\, k,G_{0})\simeq H^{2}_{fppf}(Spec\, k,\mathbb{G}_{m}) 
\qquad x\in H^{1}_{fppf}(Spec\, k,G)\simeq H^{1}_{fppf}(Spec\, k,\mathbb{Z})\times H^{2}_{et}(Spec\, k,\mathbb{G}_{m}).$$
Moreover, the natural inclusion $G_{0} \longrightarrow G$ sends the class $x_{0}$ to the class $x$.
Under the above identification the morphism induced by the inclusion becomes the 
inclusion of the second factor
$$H^{2}_{fppf}(Spec\, k,\mathbb{G}_{m}) \hookrightarrow H^{1}_{fppf}(Spec\, k,\mathbb{Z})\times H^{2}_{fppf}(Spec\, k,\mathbb{G}_{m})$$
and is thus injective.
By assumption $Eq(B,B')$ has a $k$-rational point, and thus $x=0$. We therefore have 
$x_{0}=0$, or equivalentely $Eq_{0}(B,B')(k)\neq \emptyset$. This implies that 
$B$ and $B'$ are by definition isomorphic in $Cm_{k}$, and thus are equivalent as 
Azumaya $k$-algebras. \\

We now prove the second point $(2)$. It is enough to show that a given deraz $k$-algebra 
$B$ is isomorphic in $Hom_{k}$ to the image of a non-dg $k$-algebra, as this 
$k$-algebra will then automatically satisfies assumptions $\mathbf{(Az-1)}$
and $\mathbf{(Az-2)}$, which will imply that it is an Azumaya $k$-algebra (we use here that 
$k$ is a field). 

Let $B$ be a fixed deraz $k$-algebra. If $\overline{k}$ is a given algebraic closure of $k$, we now 
from proposition \ref{p2'} that  
$B_{\overline{k}}$ is Morita equivalent to $\overline{k}$. This is equivalent to 
say that $D(B_{\overline{k}})$ possesses a compact generator whose endomorphism 
dg-algebra is quasi-isomorphic to $\overline{k}$. From our corollary \ref{c1} we know that 
$B$ is a smooth and proper $k$-dg-algebra, and therefore,  
as it is shown in \cite[Thm. 3.6]{tova}, the $(\infty)$-stack 
of perfect $B$-dg-modules is locally of finite type over $k$. A consequence of this is the existence of
a finite field extension $k'/k$, and of a compact generator $E' \in D(B_{k'})$, whose endomorphism 
dg-algebra is quasi-isomorphic to $k'$. In other words, $B_{k'}$ becomes isomorphic in 
$Hom_{k'}$ to $k'$. We consider the natural adjunction morphism of $k$-dg-algebras
$B \longrightarrow B_{k'}$, and the corresponding forgetful functor
$$D(B_{k'}) \longrightarrow D(B).$$
We let $E$ be the image of $E'$ by this functor. The underlying complex of $E$ is perfect over $k$, and 
as $B$ is smooth this implies that $E$ is a compact object in $D(B)$. Moreover $B$ is a direct factor
of $B_{k'}\simeq B^{d}$ in $D(B)$ (here $d$ is the degree of the extension 
$k'/k$) and thus belongs to the thick triangulate sub-category generated by $E$. 
The object $E$ is thus a compact generator of $D(B)$. We now show that the $k$-dg-algebra 
$\mathbb{R}\underline{End}_{B}(E)$ is quasi-isomorphic to a $k$-algebra. 
Indeed, for this it is enough to show that the natural morphism
$$k' \longrightarrow \mathbb{R}\underline{End}_{B}(E)\otimes_{k}k'$$
is a quasi-isomorphism, or equivalently that 
$$k' \simeq \mathbb{R}\underline{End}_{B_{k'}}(E\otimes_{k}k').$$
But we have isomorphisms in $D(B_{k'})$
$$E\otimes_{k}k' \simeq E'\otimes_{k'}k'\otimes_{k}k' \simeq (E')^{d},$$
where $d=deg(k'/k)$.
By assumption on $E'$ we have $\mathbb{R}\underline{End}_{B_{k'}}(E')\simeq k'$, showing that 
$\mathbb{R}\underline{End}_{B_{k'}}(E\otimes_{k}k')$
is quasi-isomorphic to the matrix $k'$-algebra $M_{d}(k')$. 
This implies that the derived endomorphism $k$-dg-algebra of the $B$-dg-module $E$ is 
quasi-isomorphic to a non-dg algebra $B_{0}$ over $k$. As $E$ is a compact generator we have
that $B_{0}$ and $B$ are isomorphic in $Hom_{k}$.
\hfill $\Box$ \\

\begin{rmk}\label{r2}
\emph{The bijection on the set of invertible objects between the categories
$Cm_{k}$ and $Hom_{k}$ does not extend to a stronger form of equivalence. For instance, 
the unit object in $Hom_{k}$ has more automorphisms than in $Cm_{k}$ (the one corresponding to the  
objects $k[n]\in D(k)$, for $n \in  \mathbb{Z}$). The proposition is also wrong over more general rings, 
as we will see that there exist local rings over which there are deraz algebras non Morita equivalent to 
non-derived Azumaya algebras. We will however prove that the natural map from equivalent classes of 
Azumaya algebras to deraz algebras up to Morita equivalences stays injective over any scheme (or more 
generally any stack)}.   
\end{rmk}

\subsection{Local triviality for the \'etale topology}

To finish this first section we prove the local triviality of deraz $A$-algebras for the \'etale
topology. The proof starts with the case of a base field, and proceeds with a 
deformation theory argument. For the statement, and the proof, of the next proposition 
we use the notion of derived stacks, and we refer to \cite[\S 2.2]{hagII} or
\cite{to2} for references. \\

\begin{prop}\label{p4}
Let $A$ be a commutative simplicial ring and $B$ be a deraz $A$-algebra. Then, there exists an 
\'etale covering $A \longrightarrow A'$ such that $B_{A'}$ is Morita equivalent to $A'$. 
\end{prop}

\textit{Proof:}
Let $A$ and $B$ as in the statement of the proposition. We work in the category of 
derived stacks over $A$, for the fppf topology, which we simply denote by $dSt(A)$
(see \cite{to3} for considerations on the fppf topology). 
We consider the derived stack $Eq(A,B)$, of Morita equivalences from 
$A$ to $B$. By definition, its values over a simplicial $A$-algebra $A'$ is the nerve of 
the category of quasi-isomorphisms between all compact $B^{op}_{A'}$-dg-modules $E$ whose corresponding 
morphism $A' \rightarrow B_{A'}$ is an isomorphism in $Hom_{A'}$.
This is an open sub-stack of $\mathcal{M}_{B}$, of all compact $B^{op}$-dg-modules, constructed and
studied in \cite{tova} (we use here Cor. \ref{c1} to insure that $B$ is smooth and proper in order to
apply theorem 3.6 of \cite{tova}). 
In particular, $Eq(A,B)$ is a geometric derived stack which is locally of finite presentation over $A$. 
Let $A'$ be any simplicial $A$-algebra, and $E : 
\mathbb{R}\underline{Spec}\, A' \longrightarrow Eq(A,B)$ a given $A'$-point. 
Then, according to \cite[Cor. 3.17]{tova}, the tangent complex of $Eq(A,B)$ at the point $E$ is given by 
$$\mathbb{T}_{Eq(A,B),E}\simeq \mathbb{R}\underline{Hom}_{B_{A'}^{op}}(E,E)[1],$$
where the morphism $E$ is considered as an object in $D(B_{A'}^{op})$. But, as $E$ induces an isomorphism
$A' \simeq B_{A'}$ in $Hom_{A'}$, the functor
$$\phi_{E} : D(A') \longrightarrow D(B_{A'}^{op}),$$
sending $M$ to $M\otimes_{A'}^{\mathbb{L}}E$ is an equivalence of triangulated categories. 
Therefore we have
$$\mathbb{R}\underline{Hom}_{B_{A'}^{op}}(E,E)\simeq A'.$$
This implies that $\mathbb{T}_{Eq(B,A),E}\simeq A'[1]$ is of negative Tor amplitude, and thus that the derived stack $Eq(B,A)$ is smooth over $A$ (see \cite[\S 2.2.5]{hagII}). Moreover, the proposition 
\ref{p2'} shows that $Eq(B,A) \longrightarrow \mathbb{R}\underline{Spec}\, A$ is surjective over all 
geometric points, and it is therefore a smooth covering. 
This map possesses then local sections for the 
\'etale topology (see \cite[Lem. 2.2.3.1]{hagII}), or equivalentely, $B$ is locally, for the \'etale topology on $A$, Morita 
\'equivalent to $A'$. \hfill $\Box$ \\

The previous proposition implies the following refinement of Prop. \ref{p2'}.

\begin{cor}\label{c3}
Let $k$ be a separably closed field. Then, any deraz $k$-algebra is quasi-isomorphic to 
a graded finite dimensional matrix algebra $\underline{End}^{*}_{k}(V)$. 
\end{cor}

\section{Derived Azumaya algebras over derived stacks}

This section uses the notion of $\infty$-categories (truly only $(1,\infty)$-categories are considered), 
and we will use the approach based on Segal categories (see \cite{sim})
and simplicially enriched categories recalled in \cite[\S 1]{tove2}. 
Therefore, for us $\infty$-\emph{category} is synonimus to \emph{Segal category}, and 
\emph{strict $\infty$-category} to \emph{$\mathbb{S}$-category}. Needless to say that any other, reasonable 
enough, theory of $\infty$-categories could be used equivalently (for instance any of the one presented 
in \cite{ber}). We leave the reader to chose is favorite and to make its own translations of the constructions 
and  results presented in the sequel.

\subsection{The derived stack of locally presentable dg-categories}

We have, up to now, ignored set theoritical issues. In this paragraph we will make some constructions 
involving non-small objects, that will be called \emph{big}, but also bigger object that will be called 
\emph{very big}, in a way that big objects are small with respect to very big one. It is unclear to me how to 
avoid using universes in order to make these constructions meaningful without going into messy cardinality
bound issues. We will therefore assume here that we have added the universe axioms to 
the axiom of sets theory.
I do not think that the results we prove in the sequel really depend on axioms of universes, but rather 
that it is 
an easy way to write things down. Once again, at the very end, when these constructions will be applied to 
prove our main theorems about existence of derived Azumaya algebras, only a very small part of the objects 
that we will describe will be truly used, a part which can be check to exits in ZFC without any 
universe axioms. \\

In this paragraph we use the terminology and notations from the homotoy theory of dg-categories. We
refer the reader to \cite{tab,to1} for more details. In particular, for a dg-category $T$ we note
$[T]$ the associated homotopy category, obtained by replacing all complexes of morphisms by their
$H^{0}$. When $T$ is a triangulated dg-category (for instance when $T$ has small limits and 
colimits), then $[T]$ is naturally endowed with a triangulated structure. \\

We let $\mathbb{U} \in \mathbb{V} \in \mathbb{W}$ be three universes. Sets elements of $\mathbb{U}$ will be called
\emph{small} (which will be implicit if no other adjective is specified), elements of $\mathbb{V}$ \emph{big}, and the ones
of $\mathbb{W}$ \emph{very big}. We let $s\z-CAlg$ be the category of (small) commutative rings, endowed its usual 
model structure. For any $A\in s\z-CAlg$ we $Dg(A)$ be the category of big $A$-dg-categories. We note that 
$Dg(A)$ possesses all limits and colimits along big sets, and is a $\mathbb{V}$-combinatorial model catgeory in the
sense of \cite[Appendix]{hagI} (the underlying category of $Dg(A)$ is itself is very big). For 
$u : A \rightarrow B$ a morphims of simplicial commutative rings, we have a direct image functor between very big categories
$$u_{*} : Dg(B) \longrightarrow Dg(A)$$
which consists of viewing a $B$-dg-category as a $A$-dg-category throught the morphism $A \rightarrow B$. This defines
a functor
$$Dg : s\z-CAlg^{op} \longrightarrow Cat_{\mathbb{W}},$$
from $s\z-CAlg$ to the category of very big categories. For $A\in s\z-CAlg$, the category 
$Dg(A)$ is equipped with a model structure and in particular with a notion of equivalences (here these are
the quasi-equivalences of big $A$-dg-categories). Localizing along these equivalences we get 
a new functor
$$LDg : s\z-CAlg^{op} \longrightarrow \s-Cat_{\mathbb{W}},$$
defined by the obvious formula $LDg(A):=L(Dg(A))$. 
We perform the Grothendieck's construction of \cite[\S 1.3]{tove2} to turn 
the functor $LDg$ into a fibered very big $\s$-category
$$\pi : \int LDg \longrightarrow s\z-CAlg.$$
Objects of the total space $\int LDg$ are pairs $(A,T)$, where 
$A$ is a simplicial commutative ring and $T$ is a big $A$-dg-category. The simplicial set of morphisms in 
$\int LDg$, between $(A,T)$ and $(B,U)$ are given by (see \cite[\S 1.3]{tove2} for more details)
$$Map_{\int LDg}((A,T),(B,U))=\coprod_{u : A \rightarrow B} L(Dg(A))(T,u_{*}(U)).$$

Before going further we need to remind here the notion of homotopy colimits in dg-categories. 
For $T$ a big dg-category (over some commutative simplicial ring $A$), we consider
the model category $T-Mod$, of big $T$-dg-modules. The $T$-dg-modules quasi-isomorphic to dg-modules of the form 
$T(a,-)$, for some object $a\in T$, are called \emph{co-representable}.
We then say that \emph{$T$ possesses small (homotopy) colimits}, 
if the co-representable $T$-dg-modules are stable by small homotopy limits in the
model category $T-Mod$, as well as by the operations $E\otimes_{A}^{\mathbb{L}}$, for
some small $A$-dg-module $E$. 
Dually, we can make sense of the expression  \emph{$T$ possesses small (homotopy) limits}.
As we will never consider limits and colimits in dg-categories in the strict sense, we will allow ourselves
not to mention the word "homotopy".

We will also use  the notion of adjunction between 
dg-categories. As for limits and colimits we will never mean a strict adjunction 
between dg-categories, but only a weak form of it, defined as follows. Let $f : T \longrightarrow T'$ 
be a dg-functor between dg-categories (over some base commutative simplicial ring $A$). 
We say that \emph{$f$ possesses a left adjoint} if there is 
a morphism $g : T' \longrightarrow T$, in the homotopy category of dg-categories, 
such that the two $T\otimes_{A}^{\mathbb{L}}(T')^{op}$-dg-modules
$$T(g(-),-) \qquad T(-,f(-))$$
are isomorphic in the derived category $D(T\otimes_{A}^{\mathbb{L}}(T')^{op})$. This is
equivalent to say that there is a morphism 
$$u : 1 \Rightarrow f\circ g,$$
inside the category $[\mathbb{R}\underline{Hom}(T',T')]$, such that for any $x\in T$ and
$y\in T'$, the induced morphisms (well defined in the derived category $D(A)$ of $A$-dg-modules)
$$T(g(x),y) \rightarrow T'(fg(x),f(y)) \longrightarrow T'(x,f(y))$$
is an isomorphism (here, $\mathbb{R}\underline{Hom}(T',T')$ is the derived internal
Hom of dg-categories of \cite[Thm. 6.1]{to1}). \\

Now, recall that a big $A$-dg-category $T$ is 
\emph{locally presentable} if it possesses a small set of $\lambda$-small objects, for 
$\lambda$ a regular cardinal elements of $\mathbb{U}$, and if it is closed by 
small colimits (in the sense above) (this is the dg-version 
of locally presentable $\s$-categories of \cite{lu1}). Equivalently, there exists 
a small $A$-dg-category $T_{0}$, and a fully faithful
$A$-dg-functor (recall that $\widehat{T_{0}}$ 
is by definition the $A$-dg-category of all cofibrant $T^{op}_{0}$-dg-modules, 
see \cite[\S 7]{to1})
$$T \hookrightarrow \widehat{T_{0}},$$
having a right adjoint (in the sense of the homotopy theory of 
$A$-dg-categories) which commutes furthermore with 
$\lambda$-filtered colimits, for $\lambda$ a (small) regular cardinal. Another equivalent 
definition is to say that there exists a small dg-category $T_{0}$ over $A$, 
and a small set of morphisms $S$ in $\widehat{T_{0}}$, such that 
$T$ is a localization of $\widehat{T_{0}}$ along $S$, \emph{inside the homotopy theory 
of dg-categories with colimits and continuous dg-functors} (sometimes
called \emph{Bousfield localizations}, to avoid possible confusions with the
usual localizations inside the homotopy theory of dg-categories and all dg-functors). 
This, in turn, is equivalent to say that 
$T$ is equivalent to the full sub-dg-category of $\widehat{T_{0}}$ consisting of
all $S$-local objects, that is objects $E$ such that for any 
$x\rightarrow y$ in $S$, the induced morphism on complexes of morphisms
$$\widehat{T_{0}}(y,E) \longrightarrow \widehat{T_{0}}(x,E)$$
is a quasi-isomorphism. Note that the inclusion dg-functor of 
$S$-local objects in $\widehat{T_{0}}$ always has a left adjoint, constructed
by the so-called small objects argument. 

By definition, 
locally presentable dg-categories are big dg-categories, even if we do not mention the adjective
\emph{big} explicitely.
A $A$-dg-functor $T \longrightarrow T'$ between 
locally presentable $A$-dg-categories is \emph{continuous} if it commutes 
with small direct sums, and thus will all small colimits (recall that 
a $A$-dg-functor always commutes with finite limits and colimits when they exist). 

We let $Dg^{lp}(A)$ be the (non full !) sub-category of $Dg(A)$ consisting of locally presentable 
$A$-dg-categories and continuous morphisms between them. For $A \rightarrow B$, the direct image
functor $Dg(B) \longrightarrow Dg(A)$ obviously preserves locally presentable objects as well
as continuous dg-functors. We consider the sub-category $\int LDg^{lp} \subset \int LDg$ consisting of all pairs
$(A,T)$, with $T$ a locally presentable $A$-dg-category, and continuous morphisms between 
those pairs. The natural projection
$$\pi : \int LDg^{lp} \longrightarrow s\mathbb{Z}-CAlg$$
is still a fibered very big $\s$-category.

\begin{lem}\label{l3}
The $\s$-functor $\int LDg^{lp} \longrightarrow s\mathbb{Z}-CAlg$ 
is also a cofibered $\s$-category, so is a bifibered $\s$-category. 
\end{lem}

\textit{Proof of the lemma:} By definition, the results 
states that for $u : A \rightarrow B$ in $s\z-CAlg$, 
the $\s$-functor
$$u_{*} : LDg^{lp}(B) \longrightarrow LDg^{lp}(A)$$
possesses a left adjoint $u^{*}$ (in the sense of $\s$-categories). This is equivalent to the
existence, for an arbitrary $T \in LDg^{lp}(A)$, of an object $T\otimes_{A}^{ct}B \in LDg^{lp}(B)$, 
together with a morphism $T \longrightarrow u_{*}(T\otimes_{A}^{ct}B)$,
such that for any $U\in LDg^{lp}(B)$, the composite morphism
$$Map_{LDg^{lp}(B)}(T\otimes_{A}^{ct}B,U) \longrightarrow
Map_{LDg^{lp}(A)}(u_{*}(T\otimes_{A}^{ct}B),u_{*}(U)) \longrightarrow
Map_{LDg^{lp}(A)}(T,u_{*}(U))$$
is an equivalence (see \cite[\S 1.2]{tove2} for the notion
of adjunction between $\s$-categories). The object $T\otimes_{A}^{ct}B$ can be explicitly 
constructed as follows.
We write $T$ as the localization, in the sense of locally presentable $B$-dg-categories,  
of $\widehat{T_{0}}$, for $T_{0}$ a small $A$-dg-category, along 
a small set of morphisms $S$ in $\widehat{T_{0}}$. That is
we assume that $T$ is the full sub-dg-category of $\widehat{T_{0}}$ consisting of all 
$S$-local objects. 
We consider the small $B$-dg-category 
$T_{0}\otimes_{A}^{\mathbb{L}}B$, which comes equiped with a small set of morphisms
$S\otimes_{A}^{\mathbb{L}}B$ in $\widehat{T_{0}\otimes_{A}^{\mathbb{L}}B}$, consisting of 
all morphisms of $T_{0}^{op}\otimes_{A}^{\mathbb{L}}B$-dg-modules obtained by 
tensoring a morphism in $S$ by $B$ over $A$, or in other words the image of $S$ by the base change
dg-functor
$$-\otimes_{A}^{\mathbb{L}}B : \widehat{T_{0}} \longrightarrow \widehat{T_{0}\otimes_{A}^{\mathbb{L}}B}.$$
 We define $T\otimes_{A}^{ct}B$ to be the localization 
of the locally presentable $B$-dg-category $\widehat{T_{0}\otimes_{A}^{\mathbb{L}}B}$ along the set 
$S\otimes_{A}^{\mathbb{L}}B$, that is the full sub-dg-category of
all $S\otimes_{A}^{\mathbb{L}}B$-local objects in $\widehat{T_{0}\otimes_{A}^{\mathbb{L}}B}$. 
Note that these are precisely the $T_{0}^{op}\otimes_{A}^{\mathbb{L}}B$-dg-modules
$E$ which are $S$-local when simply considered as $T_{0}^{op}$-dg-modules. Therefore, 
the natural $A$-dg-functor $\widehat{T_{0}} \longrightarrow
u_{*}(\widehat{T_{0}\otimes_{A}^{\mathbb{L}}B})$ induces, by localization, a natural $A$-dg-functor
$$T \longrightarrow u_{*}(T\otimes_{A}^{ct}B),$$
and \cite[Thm. 7.1, Cor. 7.6]{to1} imply that the required universal property is satisfied. \hfill $\Box$ \\

We consider the bifibered $\s$-category of locally presentable dg-categories
$$\int LDg^{lp} \longrightarrow s\z-CAlg.$$
Let $A \longrightarrow B_{*}$ be a fpqc hyper-coverings in $s\z-CAlg$. It defines 
a co-augmented cosimplicial object
$$H : \Delta_{+} \longrightarrow s\z-CAlg.$$
Because of the lemma \ref{l3} above, 
given an object $T \in LDg^{lp}(A)$, there is a unique, up to equivalence, $\s$-functor over $s\z-CAlg$
$$F : \Delta_{+} \longrightarrow \int LDg^{lp},$$
such that $F([-1])\simeq T$, and such that $F$ sends every morphism in $\Delta_{+}$ to a 
cocartesian morphism in $\int LDg^{lp}$. As $\mathbb{Z}$ is the
initial object in $s\z-CAlg$, we have a natural $\s$-functor over $s\z-CAlg$
$$\int LDg^{lp} \longrightarrow LDg^{lp}(\mathbb{Z})\times s\z-CAlg,$$
which simply sends every fiber $LDg^{lp}(A')$ to $LDg^{lp}(\mathbb{Z})$ using the direct image $\s$-functor
along $\mathbb{Z} \longrightarrow A'$. Precomposing with $F$ we get a new $\s$-functor
$$\Delta^{+} \longrightarrow LDg^{lp}(\mathbb{Z}).$$
We have thus seen that the consequence of lemma \ref{l3} is the existence, for any locally presentable
$A$-dg-category $T$, and any fpqc hyper-coverings $A\rightarrow B_{*}$, of a co-augmented cosimplicial locally presentable dg-category $T^{*}$, with $T^{-1}=T$. This produces  a natural morphism
$$T=T^{-1} \longrightarrow Holim_{n \in \Delta}T^{n}.$$
Intuitively, this
cosimplicial object is given by 
$$T^{n}=T\otimes_{A}^{ct}B_{n},$$
(which is only an intuitive picture as many homotopy coherences must also be given to describe the object $T^{*}$), 
and the previous natural morphism will therefore be denoted symbolically by
$$T \longrightarrow Holim_{n \in \Delta}T\otimes_{A}^{ct}B_{n}.$$
By construction this morphism actually exists not only in dg-categories over $\mathbb{Z}$ but in 
dg-categories over $A$ (simply replace $s\z-CAlg$ by $sA-CAlg$ in the construction above). 

\begin{df}\label{d4}
\emph{Let $A$ be a commutative simplicial ring and $T$ be a locally presentable $A$-dg-category.
We say that} $T$ has (fpqc) descent \emph{if for any fpqc hyper-coverings $B \rightarrow B_{*}$ in $sA-CAlg$, the natural morphism}
$$T\otimes_{A}^{ct}B \longrightarrow Holim_{n\in \Delta}T\otimes_{A}^{ct}B_{n}$$
\emph{is an equivalence of dg-categories}.
\end{df}

The locally presentable dg-categories with descent are clearly stable by the base change $-\otimes_{A}^{ct}B$, and thus
form a full sub $\s$-category $\int LDg^{lp,desc} \subset \int LDg^{lp}$, which is cofibered over
$s\z-CAlg$. Using \cite[Prop. 1.4]{tove2} we turn this last cofibered $\s$-category into a functor
$$\mathbb{D}g^{lp,desc} : s\z-CAlg \longrightarrow \s-Cat_{\mathbb{W}},$$
that is a \emph{derived pre-stack}. \\

We will need the following behaviour of dg-categories with descent with respect to homotopy limits. 

\begin{lem}\label{l4}
\begin{enumerate}
\item For any commutative simplicial ring $A$ the very big $\s$-category $\mathbb{D}g^{lp,desc}(A)$
has all small limits. 
\item For any morphism $u : A \rightarrow B$ in $s\z-CAlg$, the two $\s$-functors
$$-\otimes_{A}^{ct}B : \mathbb{D}g^{lp,desc}(A) \longrightarrow \mathbb{D}g^{lp,desc}(B)$$
$$u_{*} : \mathbb{D}g^{lp,desc}(B) \longrightarrow \mathbb{D}g^{lp,desc}(A)$$
commute with small limits. 
\end{enumerate}
\end{lem}

\textit{Proof:} We start to observe that, for any $A\in s\z-CAlg$, the $\s$-category $LDg^{lp}(A)$
has small limits. This follows from the dg-analog of the fact that 
locally presentable $\s$-categories are stable by limits (see \cite{lu1}), and can actually be deduced
from it by considering underlying $\s$-categories of dg-categories as follows. 

For any dg-category $T$ over $A$, we can consider its complexes of morphisms and apply the
Dold-Kan correspondence to get simplicial sets. As the Dold-Kan correspondence is 
a weak mono\"\i dal functor, we can use the composition morphisms in $T$ in order to define
the composition morphisms on the corresponding simplicial sets (see \cite[Thm. 4.16]{tab2} for 
more details). We obtain this way
a simplicial category $sT$, functorially associated to $T$, that we consider
as a strict $\s$-category. This defines an $\s$-functor
$$s : LDg(A) \longrightarrow \s-Cat,$$
from the $\s$-category of dg-categories over $A$ to the $\s$-category of $\s$-categories
$\s-Cat$. 
When $T$ has small colimits the $\s$-category $sT$ has small colimits. Moreover, 
if $T$ has small colimits, then $T$ is a locally presentable dg-category over $A$ if and only if
$sT$ is a locally presentable $\s$-category. As the $\s$-functor $T \mapsto sT$ 
commutes with limits (as it is obtained by localizing a right Quillen functor, 
see \cite{tab2}), we deduce from \cite[Prop. 5.5.3.12]{lu1} that $LDg^{lp}(A)$ is stable by small limits
in $LDg(A)$, and thus itself possesses all small limits. 

The next step is to see that, for any morphism of commutative simplicial rings $A \longrightarrow B$, 
the base change $\s$-functor
$$-\otimes_{A}^{ct}B : LDg^{lp}(A) \longrightarrow LDg^{lp}(B)$$
commutes with small limits. Using \cite[Thm. 7.2, Cor. 7.6]{to1}, as well as the explicit 
construction of $T\otimes_{A}^{ct}B$ in terms of localizations (see the proof
of lemma \ref{l3}), we see that there is a natural equivalence
$$T\otimes_{A}^{ct}B \simeq \rh(B,T),$$
where the right hand side denotes the derived internal Hom object in $LDg(A)$. As it is
clear that $T \mapsto \rh(B,T)$ commutes with small limits (because it is a right adjoint), we
deduce that $-\otimes_{A}^{ct}B$ commutes with small limits. 

We are now ready to prove $(1)$ and $(2)$ of lemma \ref{l4}. Let
$$T_{*} : I \longrightarrow \mathbb{D}g^{lp,desc}$$
be a small diagram of locally presentable dg-category with descent over $A$, 
$T=Holim\, T_{i}$, and $A \rightarrow B_{*}$ be a fpqc hyper-covering of commutative simplicial rings. 
Now that we know that base change of locally presentable dg-categories commutes with limits, we have
$$T=Holim\, T_{i} \simeq Holim_{i}\left( 
Holim_{n}T_{i}\otimes_{A}^{ct}B_{n} \right) \simeq 
Holim_{n}\left( T\otimes_{A}^{ct}B_{n}\right),$$
showing that $T$ is a locally presentable dg-category with descent over $A$. This shows $(1)$, and
$(2)$ follows from what we have already seen, that $-\otimes_{A}^{ct}B$ commutes with small limits. 
\hfill $\Box$ \\

The main properties of locally presentable dg-categories with descent is the following gluing property. 

\begin{thm}\label{t1}
The derived prestack $\mathbb{D}g^{lp,desc}$ defined above is a stack for the fpqc topology
on the model category $(s\z-CAlg)^{op}$.
\end{thm}

\textit{Proof:} This follows easily from our descent criterion \ref{ca2}, proved in the appendix. 
The criterion is given in the language of $\s$-sites, but we can simply localize
our model site $(s\z-CAlg)^{op}$ to consider $\mathbb{D}g^{lp,desc}$ as a prestack over 
the $\s$-site $L(s\z-CAlg)^{op}$. The descent condition over the model site $(s\z-CAlg)^{op}$
and over $L(s\z-CAlg)^{op}$ being equivalent, we can simply apply the corollary 
\ref{ca2} over the $\s$-site $L(s\z-CAlg)^{op}$.

The condition $(1)$ is our lemma \ref{l4}, and the right adjoints of $(2)$ are simply the forgetful 
$\s$-functors $u_{*}$, which are clearly conservative. Condition $(3)$ follows from 
the compatibility of the completed tensor product $-\otimes_{A}^{ct}B$ with 
the homotopy base change of rings
$$-\otimes_{A}^{ct}(B\otimes_{A}^{\mathbb{L}}B') \simeq
-\otimes_{B}^{ct}B',$$
which itself follows from the explicit construction of $-\otimes_{A}^{ct}B$, and the usual 
associative property of the derived tensor product of simplicial rings. Finally, condition 
$(4)$ follows easily from the fact that we restrict to dg-categories with descent: we have
for any fpqc hyper-coverings $A \rightarrow B_{*}$, any $T\in \mathbb{D}g^{lp,desc}$, and any 
$K\in \mathbb{D}g^{lp,desc}(\mathbb{Z})$
$$Map_{A}(K\otimes^{ct}A,T) \simeq
Map_{A}(K\otimes^{ct}A,Holim_{n}(T\otimes_{A}^{ct}B_{n})) \simeq
Holim_{n}Map_{A}(K\otimes^{ct}A,(T\otimes_{A}^{ct}B_{n})).$$
We conclude that $\mathbb{D}g^{lp,desc}$ is a stack for the fpqc topology. 
\hfill $\Box$ \\

\begin{rmk}
\emph{We could have considered a weaker descent property, simply requiring 
descent for fppf hypercoverings, and even more restrictively fppf \v{C}ech descent
(i.e. descent with respect of nerves of fppf coverings). It can be proven that 
every locally presentable dg-category has the fppf \v{C}ech descent (Jacob Lurie, private 
communication). The argument follows the same lines as in the proof of our 
theorem \ref{t2}, and proceeds by breaking the fppf descent into two
independant steps: descent for finite flat maps and Nisnevich descent. As a consequence, 
the derived prestack $A \mapsto \mathbb{D}g^{lp}$ is a \v{C}ech stack for the fppf
topology. }

\emph{As we will see in the sequel there many examples of dg-categories with descent, 
and the two classes of all locally presentable dg-categories and
of dg-categories with descent are rather closed.}
\end{rmk}

\begin{df}\label{d4'}
\emph{Let $F$ be a derived stack. A} locally presentable dg-category over $F$ \emph{is a morphism
of derived stacks}
$$F \longrightarrow \mathbb{D}g^{lp,desc}.$$
\end{df}

In intuitive terms, a locally dg-category over $F$ consists of the following data.
\begin{enumerate}
\item For any commutative simplicial ring $A$ and any morphism $u : \mathbb{R}\underline{Spec}\, A \longrightarrow F$,
a locally presentable dg-category $T_{u}$ with descent over $A$.
\item For any commutative simplicial rings $A$ and $B$, and any commutative diagram
$$\xymatrix{
\mathbb{R}\underline{Spec}\, A \ar[r]^-{f} \ar[rd]_-{u} & \mathbb{R}\underline{Spec}\, B \ar[d]^-{v} \\
 & F,}$$
a dg-functor $f^{*} : T_{v} \longrightarrow T_{u}$, of dg-categories over $A$, such that the induced
morphism $T_{v} \otimes_{B}^{ct}A \longrightarrow T_{u}$ is an equivalence.  
\item Homotopy coherences conditions on $f^{*}$ with respect to compositions. 
\end{enumerate}

The first two points above give a rather clear picture of what dg-categories over a derived stack $F$ are, and
the technical part of this notion is really in point $(3)$. There are many ways of writing these
homotopy coherences, which are encoded in the fact that the morphism $F \longrightarrow \mathbb{D}g^{lp,desc}$
is only a morphism in the $\s$-category of derived stacks. \\

There are two trivial dg-categories over any given derived stack, the null dg-category and the unit dg-category.
They both are global sections
$$*=Spec\, \mathbb{Z} \longrightarrow \mathbb{D}g^{lp,desc},$$
corresponding to the $0$ dg-category, and to the unit dg-category $\widehat{\mathbb{Z}}$. Pulling back these
two morphisms along the morphism $F \longrightarrow *$ provides two dg-categories over $F$. Using the
intuitive description of dg-categories over $F$ given above, the null dg-category is the one
given by $T_{u}=0$ for all $u$, and the unit dg-category is given by $T_{u}=\widehat{A}$ for all $u$ (with the
obvious transitions dg-fucntors $f^{*}$). The unit dg-category over $F$ 
will be denoted by $\mathbf{1}_{F}$, or $\mathbf{1}$ if $F$ is clear from the context, 
and should be 
understood as some kind of \emph{categorical structure sheaf of $F$}. \\

The previous theorem shows that dg-categories with descent possess nice local properties. The following
proposition implies that all compactly generated dg-categories are with descent.

\begin{prop}\label{p5}
Let $A$ be a commutative simplicial ring and $T_{0}$ be a small dg-category over $A$. Then, 
the big dg-category $\widehat{T_{0}}$, of small $T_{0}$-dg-modules, is a locally presentable
dg-category with descent.
\end{prop}

\textit{Proof:} The dg-categories of the form $\widehat{T_{0}}$ are clearly locally presentable. Moreover, 
for any simplicial commutative $A$-algebra $B$, we have
$$\widehat{T_{0}}\otimes_{A}^{ct}B\simeq \widehat{T_{0}\otimes_{A}^{\mathbb{L}}B}.$$
Therefore, to show that $\widehat{T_{0}}$ is with descent it is enough to show that 
for any fpqc hyper-coverings $A \rightarrow B_{*}$ we have
$$\widehat{T_{0}} \simeq Holim_{n\in\Delta}\widehat{T_{0}\otimes_{A}^{\mathbb{L}}B_{n}}.$$
But, from the results of \cite{to1} we have, for an arbitrary small $A$-dg-category $T_{0}'$, and
any $A'\in sA-CAlg$
$$\widehat{T_{0}'\otimes_{A}^{\mathbb{L}}A'}\simeq \mathbb{R}\underline{Hom}((T'_{0})^{op},\widehat{A'}),$$
where $\mathbb{R}\underline{Hom}$ is the derived internal Hom between $A$-dg-categories. 
Therefore, as $\mathbb{R}\underline{Hom}$ commutes with homotopy limits (because it is 
a right adjoint in the sense of $\s$-category), 
it is enough to treat the case $T_{0}=A$. We thus have to show that 
$$\widehat{A} \simeq Holim_{n\in\Delta}\widehat{B_{n}},$$
which is nothing else than fpqc descent for quasi-coherent dg-modules, as proved in \cite[Lem. 3.1]{tova}, but 
directly at the level of dg-categories rather than only at the level of nerves of quasi-isomorphisms.
\hfill $\Box$ \\

An interesting corollary of the previous proposition together with lemma \ref{l4} is that
any homotopy limit of compactly generated $A$-dg-categories is with descent (thought it might be
non-compactly generated).

\begin{cor}\label{c4}
Any dg-category over a commutative simplicial ring $A$, which is equivalent to a homotopy limit
of dg-categories of the form $\widehat{T_{0}}$ (for $T_{0}$ a small $A$-dg-category) is a locally 
presentable $A$-dg-category with descent.
\end{cor}

\begin{rmk}\label{r3}
\begin{enumerate}
\item 
\emph{It should be noted that the proposition \ref{p5} together with theorem \ref{t1} provides an 
explicit description of the derived stack $\widetilde{\mathbb{D}g^{c}}$, defined in
\cite[\S 4.3, \S 4.4]{tove2} 
by brutal stackyfication of the prestack of compactly generated dg-categories. Indeed, 
$\widetilde{\mathbb{D}g^{c}}$ can be identified with the full sub-stack of 
$\mathbb{D}g^{lp,desc}$ consisting of all $A$-dg-categories $T$ such that there is a 
fppf covering $A \rightarrow B$ with $T\otimes_{A}^{ct}B$ equivalent to 
$\widehat{T_{0}}$ for some small $b$-dg-category $T_{0}$. In other words, 
$\widetilde{\mathbb{D}g^{c}}$ is equivalent to the derived stack of locally presentable
dg-categories with descent and with local (for the fppf topology) compact generators. Objects
in $\widetilde{\mathbb{D}g^{c}}$ are then} twisted forms of compactly generated 
dg-categories, \emph{and do form a rigid symmetric mono\"\i dal $\s$-category.}
A consequence of our theorem \ref{t2} (more precisely of 
techniques of proof of proposition \ref{p8}, \ref{p9} and \ref{p10}) is that $\mathbb{D}g^{c}$ is in fact a stack, and that 
any twisted form of a compactly generated dg-category (for the fppf topology) is itself
compactly generated.
\item \emph{There are examples of locally presentable dg-categories which do not have descent. By corollary
\ref{c4} these are not compactly generated. For instance, for a general topological space $X$, 
the dg-category $L(X)$, of complexes of sheaves on $X$, is a locally presentable dg-category which does not
satisfies descent in general.}
\end{enumerate}
\end{rmk}

We finish this paragraph with the description of the diagonal 
of the derived stack $\mathbb{D}g^{lp,desc}$. For this we let 
$A$ a commutative simplicial ring and 
$$T_{1},T_{2} : \mathbb{R}\underline{Spec}\, A \longrightarrow \mathbb{D}g^{lp,desc}$$
be two locally presentable dg-categories with descent over $A$. As $\mathbb{D}g^{lp,desc}$
is a stack in $\infty$-categories, we have a derived stack of simplicial sets 
$Map(T_{1},T_{2})$ of morphisms between the two objects $T_{1}$ and $T_{2}$. We will 
assume for simplicity that $T_{1}$ is compactly generated, that is of the form
$\widehat{T_{0}}$ for some small dg-category $T_{0}$ over $A$. Then, 
using \cite[Cor. 7.6]{to1} it is possible to prove that this derived stack can be described as follows.
We let 
$$T:=\widehat{T_{0}^{op}}\otimes_{A}^{ct}T_{2},$$
which is another dg-category with descent over $A$. Then, for any $A$-algebra $A'$, 
we have a natural equivalence
$$Map(T_{1},T_{2})(A')\simeq Map_{A-dg-cat}(A,T\otimes_{A}^{ct}A').$$
In other words, the simplicial set $Map(T_{1},T_{2})(A')$ is a classifying space 
for objects in the dg-category $T\otimes_{A}^{ct}A'$. As a consequence, we have the following 
description for the homotopy groups of $Map(T_{1},T_{2})(A')$ (see \cite[Cor. 7.6]{to1})
$$\pi_{0}(Map(T_{1},T_{2})(A'))\simeq [T\otimes_{A}^{ct}A']/iso \qquad
\pi_{1}(Map(T_{1},T_{2})(A'),x)\simeq aut(x)$$
$$\pi_{i}(Map(T_{1},T_{2})(A'),x)\simeq Ext^{1-i}(x,x),$$
where the automorphisms and ext groups are computed in the triangulated category 
$[T\otimes_{A}^{ct}A']$, the homotopy category associated to $T\otimes_{A}^{ct}A'$.

\subsection{The derived prestack of derived Azumaya algebras}

In the previous paragraph we have defined the derived stack 
$\mathbb{D}g^{lp,desc}$, of locally presentable dg-categories with descent. We will now 
define the derived stack $\mathbb{D}g\mathbb{A}lg$, of dg-algebras, and relate these two derived
stacks by means of a morphism
$$\phi : \mathbb{D}g\mathbb{A}lg \longrightarrow \mathbb{D}g^{lp,desc},$$
sending a dg-algebra to its dg-category of dg-modules. The image of $\phi$, restricted to 
deraz algebras, is a derived subprestack of $\mathbb{D}g^{lp,desc}$, which by 
definition, will be the derived prestack of deraz algebras. \\

For $A\in s\z-CAlg$, we set $A-dg-alg$, the category of (small) $A$-dg-algebras. Localizing along 
the quasi-isomorphisms we obtain an $\s$-category 
$$LDgAlg(A):=L(A-dg-alg).$$
For $u : A \rightarrow A'$ a morphism of simplicial commutative rings, we have the direct image 
functor $A'-dg-alg \longrightarrow A-dg-alg$, which after localization induces 
an $\s$-functor
$$u_{*} : LDgAlg(A') \longrightarrow LDgAlg(A).$$
This defines a functor $DgAlg : s\z-CAlg^{op} \longrightarrow \s-Cat$, of which we take the
Grothendieck's construction to obtain a fibered $\s$-category
$$\int LDgAlg \longrightarrow s\z-CAlg.$$
The existence of left adjoints  to the $\s$-functors $u_{*}$,
$$-\otimes_{A}^{\mathbb{L}}A' : LDgAlg(A) \longrightarrow LDgAlg(A'),$$
implies that $\int LDgAlg$ is a bifibered $\s$-category over $s\z-CAlg$. By 
\cite[\S 1.3]{tove2} it thus corresponds to a functor
$$\mathbb{D}g\mathbb{A}lg : s\z-CAlg \longrightarrow \s-Cat,$$
which is a derived prestack of $\s$-categories. 

\begin{prop}\label{p6}
The derived prestack $\mathbb{D}g\mathbb{A}lg$ is a stack
for the fpqc topology.
\end{prop}

\textit{Proof:} The proposition is proven using the same argument as
for the proof that quasi-coherent complexes form a derived stack (e.g. as in 
\cite[Lem. 3.1]{tova}), but replacing the model categories of dg-modules by model categories
of dg-algebras. We leave the details to the reader.  \hfill $\Box$ \\

We now define a morphism of derived stacks in $\s$-categories
$$\phi : \mathbb{D}g\mathbb{A}lg \longrightarrow \mathbb{D}g^{lp,desc}.$$
For this, it is enough to construct an $\s$-functor
$$\widetilde{\phi} : \int LDgAlg \longrightarrow \int LDg^{lp,desc},$$
covering the natural projections to $s\z-CAlg$, and with the property that 
$\widetilde{\phi}$ takes cocartesian morphisms to cocartesian morphisms. 
As $\int LDgAlg$ and $\int LDg^{lp,desc}$ are obtained by the Grothendieck
construction applied to natural functors
$$LDgAlg, LDg^{lp,desc} : s\z-CAlg^{op} \longrightarrow \s-Cat,$$
we start by defining a natural transformation
$h : LDgAlg \Rightarrow LDg^{lp,desc}$, and we will check that the induced $\s$-functor
$$\int h : \int LDgAlg \longrightarrow \int LDg^{lp,desc}$$
possesses the required properties. For this, let $A$ be a simplicial commutative ring, 
and denote by $Dg_{*}(A)$ the category of big $A$-dg-category $T$ endowed with a 
distinguished object $x \in T$ and dg-functor preserving the distinguished objects. 
We have
$$Dg_{*}(A)\simeq A/Dg(A),$$
where $A$ is considered as an $A$-dg-category with a unique object. By localization we get 
a functor
$$LDg_{*} : s\mathbb{Z}-CAlg^{op} \longrightarrow \s-Cat,$$
which is the pointed version of $LDg$. We consider $LDg^{lp,desc}_{*}$
the sub-functor of $LDg_{*}$ consisting of pointed dg-categories $(T,x)$ with
$T$ a locally presentable dg-category with descent and continuous dg-functors. We have a diagram
of functors
$$\xymatrix{
LDgAlg & LDg^{lp,desc}_{*} \ar[r]^-{p} \ar[l]_-{a} & LDg^{lp,desc},}$$
where $p$ forgets the distinguished objects and $a$ sends a pair $(T,x)$ 
to the dg-algebra $T(x,x)$. We now consider $\mathcal{D}$ the sub-functor 
of $LDg^{lp,desc}_{*}$ consisting of pairs equivalent to $(B-dg-mod^{c},B)$, 
for $B$ a dg-algebra (where $B-dg-mod^{c}$ is the dg-category of cofibrant 
$B$-dg-modules and $B$ is considered as a dg-module over itself). We have another diagram
$$\xymatrix{
LDgAlg & \mathcal{D} \ar[r]^-{p} \ar[l]_-{a} & LDg^{lp,desc},}$$
and it is easy to check that $a$ is now an equivalence. Therefore, we obtain 
a well defined morphism of $\s$-functors, well defined up to equivalence
$$h : LDgAlg \longrightarrow LDg^{lp,desc},$$
from which we obtain an $\s$-functor of fibered $\s$-categories
$$\int h : \int LDgAlg \longrightarrow \int LDg^{lp,desc},$$
which by construction preserves the cartesian morphisms. This morphism
sends a dg-algebra $B$ to the dg-category with descent $\widehat{B}$.
We claim that it does also 
preserve the cocartesian morphisms. Coming back to the definitions we see that this
is equivalent to state that for $A \rightarrow A'$ a morphism in $s\z-CAlg$, and
for $B$ a $A$-dg-algebra, the natural adjunction morphism
$$\widehat{B^{op}}\otimes_{A}^{ct}A' \longrightarrow \widehat{B^{op}\otimes_{A}^{\mathbb{L}}A'},$$
is a quasi-equivalence. But, this is true by the definition of the base 
change $\otimes_{A}^{ct}A'$, and by the results of \cite{to1}. 

We have thus defined 
$$\int h : \int LDgAlg \longrightarrow \int LDg^{lp,desc},$$
which preserves cocartesian morphisms, and therefore induces
a morphism of derived stacks of $\s$-categories
$$\phi : \mathbb{D}g\mathbb{A}alg \longrightarrow \mathbb{D}g^{lp,desc}.$$

We now work in the $\s$-category $dSt_{fppf}$, of derived stacks (of simplicial sets)
over the model category $dAff:=s\z-CAlg^{op}$, endowed with the fppf topology (see \cite{to3}). 
The derived stacks of $\s$-categories $\mathbb{D}g\mathbb{A}alg$ and $\mathbb{D}g^{lp,desc}$
give rise to derived stacks of simplicial sets by considering their \emph{underlying spaces} as in 
\cite[\S 1]{tove2}\footnote{In the sequel we will often neglect writing "$\mathcal{I}$", and assume
implicitely that the underlying space functor has been applied if necessary.}
$$\mathcal{I}(\mathbb{D}g\mathbb{A}alg) \qquad \mathcal{I}(\mathbb{D}g^{lp,desc}).$$
These derived stacks can also be considered as stacks in $\s$-groupo\"\i ds (using the equivalence
between simplicial sets and $\s$-groupo\"\i ds), and they then correspond to the maximal 
substacks in groupo\"\i ds of $\mathbb{D}g\mathbb{A}alg$ and $\mathbb{D}g^{lp,desc}$ respectively (in other
words the symbol $\mathcal{I}$ simply means that we restrict to equivalences between objects, and throw away
any non invertible morphism). 
The morphism $\phi$ induces a morphism of derived stacks
$$\phi : \mathcal{I}(\mathbb{D}g\mathbb{A}alg) \longrightarrow \mathcal{I}(\mathbb{D}g^{lp,desc}).$$
We restrict this morphism to the derived substack $\mathbb{D}g\mathbb{A}alg^{Az} \subset \mathcal{I}(\mathbb{D}g\mathbb{A}alg)$, 
consisting of all deraz algebras (this is a substack by \ref{p1} $(2)$), and its stacky image 
$\mathbb{D}g^{Az}\subset \mathcal{I}(\mathbb{D}g^{lp,desc})$. By definition, for any $A\in s\z-CAlg$, 
$\mathbb{D}g\mathbb{A}alg^{Az}(A)$ is equivalent to the nerve of the category of quasi-isomorphisms
between deraz $A$-algebras, whereas $\mathbb{D}g^{Az}$ is equivalent to the nerve of the
category of quasi-equivalences between locally presentable big $A$-dg-categories with descent, locally
equivalent (for the fppf topology) to $\widehat{B}$ for $B$ a deraz algebra. The morphism 
$\phi$ simply sends a deraz algebra $B$ to $\widehat{B}$. \\

We are now ready to define derived Azumaya algebras over an arbitrary derived stack.

\begin{df}\label{d5}
\emph{Let $F\in dSt_{fppf}$ be a derived stack}. 
\begin{enumerate}
\item \emph{The} classifying space of derived Azumaya algebras over $F$ \emph{is the full sub-simplicial set
$\mathbb{D}eraz(F)$ of $Map_{dSt_{fppf}}(F,\mathbb{D}g^{Az})$, consisting of all 
morphisms $F \rightarrow \mathbb{D}g^{Az}$ that factors throught 
$\phi : \mathbb{D}g\mathbb{A}lg^{Az} \longrightarrow \mathbb{D}g^{Az}$. }
\item \emph{The} classifying space of derived Azumaya dg-categories over $F$ \emph{is}
$$\mathbb{D}eraz^{dg}(F):=Map_{dSt_{fppf}}(F,\mathbb{D}g^{Az}).$$
\end{enumerate}
\end{df}

We now finish by the following description of the derived stack $\mathbb{D}g^{Az}$, which is
a simple corollary of the results of our first section. 

\begin{cor}\label{c5}
There exist a natural equivalence of derived stacks
$$\psi : K(\mathbb{Z},1)\times K(\mathbb{G}_{m},2) \simeq \mathbb{D}g^{Az}.$$
In particular, for any $F\in dSt_{fppf}$, there is a natural monomorphism of simplicial sets
$$\psi : \mathbb{D}eraz(F) \hookrightarrow \mathbb{H}_{fppf}^{1}(F,\mathbb{Z})\times \mathbb{H}_{fppf}^{2}(F,
\mathbb{G}_{m}):=
Map_{dSt_{fppf}}(F,K(\mathbb{Z},1))\times Map_{dSt_{fppf}}(F,K(\mathbb{G}_{m},2)),$$
and thus an injective map
$$\pi_{0}(\mathbb{D}eraz(F))\hookrightarrow H^{1}_{fppf}(F,\mathbb{Z}) \times H^{2}_{fppf}(F,\mathbb{G}_{m}).
$$
\end{cor}

\textit{Proof:} By definition, any object in $\mathbb{D}g^{Az}(A)$ is locally, for the \'etale topology, 
equivalent to
some object $\widehat{B}$, for $B$ a deraz $A$-algebra. By proposition 
\ref{p4}, $B$ is itself locally Morita equivalent, for the
\'etale topology, to the trivial deraz algebra $A$. Therefore, any object in $\mathbb{D}g^{Az}(A)$ is
locally equivalent to the trivial object $\widehat{A}$. This implies that there exists a natural 
equivalence of derived stacks
$$K(G,1) \simeq \mathbb{D}g^{Az},$$
where $G$ is the derived group stack of auto-equivalences of the trivial object. 
By \cite[Thm. 8.15]{to1}, the derived group 
stack $G$ is simply the group of invertible perfect dg-modules, which are all of the form
$\mathcal{L}[n]$, for $\mathcal{L}$ a line bundle and $n$ some integer
(the group law being given by the tensor product). We thus have a natural equivalence of 
derived group 
stacks
$$G\simeq \mathbb{Z}\times K(\mathbb{G}_{m},1).$$ 
We get this way the required equivalence
$$\psi : K(G,1)\simeq K(\mathbb{Z},1)\times K(\mathbb{G}_{m},2) \simeq \mathbb{D}g^{Az}.$$
\hfill $\Box$ \\

In \cite[\S 4.3]{tove2}, we have defined $\mathbb{D}g^{c}$, the derived prestack of 
compactly generated dg-categories, 
as a stack in symmetric mono\"\i dal $\s$-categories. The exact same construction provides a 
refined versions of the derived stacks $\mathbb{D}g\mathbb{A}lg$ and $\mathbb{D}g^{lp,desc}$ as
derived stacks of symmetric mono\"\i dal $\s$-categories. The morphism 
$$\phi : \mathbb{D}g\mathbb{A}lg \longrightarrow \mathbb{D}g^{lp,desc}$$
can also be refined as a symmetric mono\"\i dal morphism. For any derived stack $F$, these symmetric mono\"\i dal 
structures
induce symmetric mono\"\i dal group-like structures on the classifying spaces
$\mathbb{D}eraz(F)$ and $\mathbb{D}g^{Az}(F)$, and $\phi$ induces a symmetric mono\"\i dal morphism between 
these two objects. By the equivalence between symmetric mono\"\i dal group-like $\s$-groupo\"\i ds and
connective spectra (see \cite[\S 2.1]{tove2}), another equivalent way to state the existence of these
group structures is simply by stating that
$\mathbb{D}eraz(F)$ and $\mathbb{D}g^{Az}(F)$ are the $0$-th space of natural spectra, 
and that $\phi$ is itself the restriction to the $0$-th spaces of a morphism of spectra.  
As a consequence, we see that the sets
$\pi_{0}(\mathbb{D}eraz(F))$ and $\pi_{0}(\mathbb{D}g^{Az}(F))$ come equiped with natural 
abelian group structures, both induced by the tensor products of locally presentable dg-categories. 

\begin{rmk}\label{r4}
\emph{The equivalence of corollary \ref{c5} can be improved as an equivalence between symmetric mono\"\i dal 
derived stacks, where the mono\"\i dal structures is the standard one on $K(\mathbb{Z},1)\times K(\mathbb{G}_{m},2)$, and
induced by the tensor product of locally presentable dg-categories on the right hand side, as explained above.}
\end{rmk}

We now make the following definitions, generalizing the two well known definitions of Brauer group 
of a scheme.

\begin{df}\label{d6}
\emph{Let $F \in dSt_{fppf}$ be a derived stack.}
\begin{enumerate}
\item The derived algebraic Brauer group of $F$ \emph{is defined by}
$$dBr(F):=\pi_{0}(\mathbb{D}eraz(F)),$$
\emph{where the group structure on the right hand side is induced by the tensor product 
of deraz algebras. }
\item The derived categorical Brauer group of $F$ \emph{is defined by}
$$dBr_{cat}(F):=\pi_{0}(\mathbb{D}g^{Az}(F)),$$
\emph{where the group structure on the right hand side is induced by the tensor product 
of locally presentable dg-categories. }
\item The big derived categorical Brauer group of $F$ \emph{is defined by}
$$dBr_{cat,big}(F):=\pi_{0}(\mathbb{D}g^{lp,desc,inv}(F)),$$
\emph{where the group structure on the right hand side is induced by the tensor product 
of locally presentable dg-categories, and where $\mathbb{D}g^{lp,desc,inv}$ denotes the mono\"\i dal 
sub-stack of $\mathbb{D}g^{lp,desc}$ consisting of all invertible objects and equivalences between 
them.}
\item The big derived cohomological Brauer group of $F$ \emph{is defined by}
$$dBr'_{big}(F):=H^{1}_{fppf}(F,\mathbb{Z})\times H^{2}_{fppf}(F,\mathbb{G}_{m}).$$
\end{enumerate}
\end{df}

We obviously have natural injective morphisms of groups
$$dBr(F) \subset dBr_{cat}(F)\simeq dBr'_{big}(F) \subset dBr_{cat,big}(F),$$
where the isomorphism in the middle follows from corollary \ref{c5}, the first inclusion is induced by 
the morphism $\phi$, and the last inclusion follows from the fact the inclusion 
of symmetric mono\"\i dal substacks (due to proposition \ref{p2} $(2)$)
$$\mathbb{D}g^{Az} \subset \mathbb{D}g^{lp,desc,inv}.$$

We consider now two copies of the morphism $\phi$
$$\mathbb{D}g\mathbb{A}lg^{Az} \times \mathbb{D}g\mathbb{A}lg^{Az} \longrightarrow
\mathbb{D}g^{Az}\times \mathbb{D}g^{Az}.$$
Pulling-back the diagonal along this morphism provides a natural morphism of derived stacks
$$\pi : \mathcal{M} \longrightarrow \mathbb{D}g\mathbb{A}lg^{Az} \times \mathbb{D}g\mathbb{A}lg^{Az}.$$
It is easy to describe the projection $\pi$. Indeed, if $X=\mathbb{R}\underline{Spec}\, A$ is
an affine derived scheme, and $X \longrightarrow \mathbb{D}g\mathbb{A}lg^{Az} \times \mathbb{D}g\mathbb{A}lg^{Az}$, 
is a pair of deraz $A$-algebras $B_{1}$ and $B_{2}$, the fiber product
$$\mathcal{M}_{B_{1},B_{2}}:=\mathcal{M}\times_{\mathbb{D}g\mathbb{A}lg^{Az} \times \mathbb{D}g\mathbb{A}lg^{Az}}^{h}
X,$$
is the derived stack over $Spec\, A$ of Morita equivalences between $B_{1}$ and $B_{2}$, that we have denoted
by $Eq(B_{1},B_{2})$ during our proof of 
proposition \ref{p4}. This description of the correspondence $\pi$
implies the following description of the group $dBr(F)$, for a derived stack $F$. Its elements are 
simply morphisms of derived stacks $F \longrightarrow \mathbb{D}g\mathbb{A}lg^{Az}$, 
up to equivalences, which by definition are called \emph{deraz algebras over $F$}. Two such elements
are declared to be \emph{Morita equivalent} if there exists a commutative diagram of derived stacks
$$\xymatrix{
 & \mathcal{M} \ar[d]^-{\pi} \\
F \ar[r]_-{B_{1},B_{2}} \ar[ru] & \mathbb{D}g\mathbb{A}lg^{Az} \times \mathbb{D}g\mathbb{A}lg^{Az},}$$
a terminology justified by the description of $\pi$ above.
We then have a natural identification between the group 
$dBr(F)$ and the Morita equivalent classes of deraz algebras over $F$. This last remark justifies the 
fact that $\pi_{0}(\mathbb{D}eraz(F))$ is called the derived algebraic Brauer group of $F$. When 
$F$ is moreover a scheme, then $dBr(F)$ can even be described more concretely in termes
of sheaves of quasi-coherent dg-algebras, locally satisfying the two conditions $\mathbf{(Az-1)}$
and $\mathbf{(Az-2)}$ of definition \ref{d1}. 
These sheaves of dg-algebras are then taken up to Morita equivalences, 
a Morita equivalence being here a sheaf of quasi-coherent bi-dg-modules inducing local equivalences on derived
categories of dg-modules, as this is done in the concrete treatment of categorical sheaves
of \cite{tove1}.

\section{Existence of derived Azumaya algebras}

We now arrive at the question of existence of derived Azumaya algebras, which is the problem
of surjectity of the natural embeddings (see definition \ref{d6})
$$dBr(F) \subset dBr_{cat}(F) \subset dBr_{cat,big}(F),$$
for a given derived stack $F$. 
We remind here that we have proved that the derived stack
$\mathbb{D}g^{Az}$ is naturally equivalent to 
$K(\mathbb{Z},1)\times K(\mathbb{G}_{m},2)$, and therefore that we have a natural isomorphism
$$dBr_{cat}(F)\simeq H^{1}_{fppf}(F,\mathbb{Z})\times H^{2}_{fppf}(F,\mathbb{G}_{m}).$$
The question of surjectivity of the first map is therefore the question of representability of
cohomology classes in terms of derived Azumaya algebras, and is a derived generalization
of the original question of Grothendieck of representability of torsion classes in 
$H^{2}(X,\mathbb{G}_{m})$ by Azumaya algebras over $X$ (see \cite{gr}). The new cohomological
term $H^{1}(X,\mathbb{Z})$ is somehow anecdotal, but the fact that non torsion elements
are taken into account is the true new feature of the derived setting. The
surjectivity of the second map seems to us a more exotic question, a bit outside of the scope of
the present paper. In this section we will show our main theorem, which implies that the 
morphism $\phi$ is surjective for quasi-compact and quasi-separated schemes (and more generally 
for a large class of Deligne-Mumford stacks). \\

\subsection{Derived Azumaya algebras and compact generators}

Before starting to study the surjectivity of the inclusion $dBr(F) \subset dBr_{cat}(F)$, we will 
give a simple and purely categorical criterion for a class $\alpha \in dBr_{cat}(F)$ to belongs
to $dBr(F)$. For this we will first study a slightly more general problem. Let $F$ be a fixed derived stack
and 
$$\alpha : F \longrightarrow \mathbb{D}g^{lp,desc}$$
be a dg-category with descent on $F$
(here we really mean a morphism of derived stacks from $F$ to 
$\mathcal{I}(\mathbb{D}g^{lp,desc})$, the underlying derived stacks of simplicial sets 
of the stacks of $\s$-categories $\mathbb{D}g^{lp,desc}$). 
We first study the existence of a dg-algebra $B_{\alpha}$ over $F$, 
with $\phi(B)=\alpha$. In other words, we ask the question of existence of a commutative
diagram of derived stacks
$$\xymatrix{
 & \mathbb{D}g\mathbb{A}lg \ar[d]^-{\phi} \\
 F \ar[r]_-{\alpha} \ar[ru]^-{B_{\alpha}} & \mathbb{D}g^{lp,desc}.}$$
When $B_{\alpha}$ exists we will say that \emph{$B_{\alpha}$ realizes $\alpha$}. 
For this, we let $L_{\alpha}(F)$ be the dg-category of global sections of $\alpha$
$$L_{\alpha}(F):=\Gamma(F,\alpha) \in \mathbb{D}g^{lp,desc}(\mathbb{Z}),$$
defined as follows. The morphism  $\alpha$ corresponds, by \cite[Prop. 1.4]{tove2}, to a cartesian 
$\s$-functor
$$u_{\alpha} : \int F \longrightarrow \int LDg^{lp,desc}.$$
We then use the natural $\s$-functor
$$\int LDg^{lp,desc} \longrightarrow LDg(\mathbb{Z}),$$
obtained by functoriality because $\mathbb{Z}$ is initial in $s\z-CAlg$. By definition we have
$$L_{\alpha}(F):=Holim_{\int F} u_{\alpha} \in LDg(\mathbb{Z}).$$
The homotopy category 
of $L_{\alpha}(F)$ is a triangulated category denoted by $D_{\alpha}(F)$. 

\begin{df}\label{dtwist}
\emph{With the same notations as above, the triangulated category $D_{\alpha}(F)$
is called the} $\alpha$-twisted derived category of $F$. \emph{The dg-category 
$L_{\alpha}(F)$ is it self called the} $\alpha$-twisted derived dg-category of $F$.
\end{df}

Note that by definition $L_{\alpha}(F)$ is a priori only a big dg-category. When 
$F$ is not too big, then $L_{\alpha}(F)$ can be proven
to be a locally presentable dg-category (see the corollary \ref{clglu} below). 

The construction $L_{\alpha}(F)$ is functorial in $F$ as follows. 
If $F' \rightarrow F$ is a morphism of derived stacks, and
$\alpha$ is a locally presentable dg-category over $F$. There is a natural 
$\s$-functor
$$\int F' \longrightarrow \int F,$$
and thus an induced morphism on the corresponding limits
$$L_{\alpha}(F')\longrightarrow L_{\alpha}(F),$$
where we still denote by $\alpha$ the dg-category $\alpha$ pulled back on $F'$. 
With a bit more work we can enhanced this construction to an $\s$-functor
$$L_{\alpha} : (dSt/F)^{op} \longrightarrow \mathbb{D}g(\mathbb{Z}),$$
from the $\s$-category of derived stacks over $F$ to the $\s$-category 
of big $\z$-dg-categories. 

\begin{lem}\label{lglu}
Let $F$ be a derived stack and $F_{i}$ be a small diagram of objects in $dSt/F$. Then, for 
any locally presentable dg-category $\alpha$ over $F$ the natural dg-functor
$$L_{\alpha}(Hocolim_{i}F) \longrightarrow Holim_{i} L_{\alpha}(F_{i})$$
is a quasi-equivalence.
\end{lem} 

\textit{Proof:} The $\s$-category $dSt/F$ is an $\s$-topos having the $\s$-site
$dAff/F$, of derived affine schemes over $F$, as a generating site. Therefore, 
there is an equivalence between the $\s$-category of $\s$-functors
$$(dSt/F)^{op} \longrightarrow \mathbb{D}g(\z)$$
sending colimits to limits, and $\s$-functors $(dAff/F)^{op} \longrightarrow \mathbb{D}g(\z)$
satisfying the descent condition for fppf hypercoverings. The lemma follows from the fact that 
$$L_{\alpha}(-) : (dSt/F)^{op} \longrightarrow \mathbb{D}g(\z)$$
 is the left Kan extension of the 
$\s$-functor 
$$u_{\alpha} : \int F \simeq (dAff/F)^{op} \longrightarrow \mathbb{D}g(\z),$$
and that this last $\s$-functor does have descent for the fppf topology (because it 
has descent for the fpqc topology by definition of dg-categories with descent). 
\hfill $\Box$ \\

\begin{cor}\label{clglu}
Let $F$ be a small derived stacks (i.e. a small colimit of representable
derived stacks), and $\alpha$ a locally presentable dg-category over $F$. Then
the dg-category $L_{\alpha}(F)$ is locally presentable with descent.  
\end{cor}

\textit{Proof:} Follows from the previous lemma and from lemma \ref{l4}. \hfill $\Box$ \\

Note that when $\alpha=\mathbf{1}$ is the trivial dg-category over $F$, 
$D_{\alpha}(F)$ is simply the derived category of quasi-coherent complexes on $F$ (this can be taken as a 
definition of the 
derived category of quasi-coherent complexes on $F$, or can be proven to be 
equivalent to other definitions, see e.g. \cite{pr,to2,lu4}). \\

We come back to a small derived stack $F$ together with a 
locally presentable dg-category $\alpha$ over $F$. 
By adjunction, any
object $E\in D_{\alpha}(F)$ gives rise to a unique dg-functor $E : \widehat{\mathbb{Z}} \longrightarrow 
L_{\alpha}(F)$ pointing $E$, and thus to a morphism
$$\mathbf{1}_{F}  \longrightarrow \alpha,$$
of locally presentable dg-categories over $F$ (here 
$\mathbf{1}_{F}$ is the unit dg-category over $F$). Evaluating this morphism at 
a morphism $u : \mathbb{R}\underline{Spec}\, A \longrightarrow F$, we get a dg-functor
$A \longrightarrow u^{*}(\alpha)$, and therefore an object
$E_{u}\in u^{*}(\alpha)$. Keeeping these notations we have the following definition.

\begin{df}\label{dloc}
\emph{Let $F$ be a small derived stack together with a locally presentable
dg-category $\alpha$ over $F$. An object $E\in D_{\alpha}(F)$
is a} a compact local generator \emph{if for any commutative simplicial ring $A$,  
and any morphism }
$$u : \mathbb{R}\underline{Spec}\, A \rightarrow F$$ 
\emph{the object $E_{u}$ is compact 
generator for the triangulated dg-categor $[u^{*}(\alpha)]$}. 
\end{df}

Note that, despite the terminology, a compact local generator is not 
necessarily a compact object in $D_{\alpha}(F)$. Note also that 
a important point in the above definition is that $E$ is globally defined, and that 
a compact local generator should not be confused with the notion of
\emph{local compact generators}, that would rather be compact generators
defined locally on $F$. The question to know if local data of compact generators gives rise
to a global compact local generator is the main question concerning the existence of
derived Azumaya algebras, as we will see later.  \\

The following result shows that, over an affine object, a compact generator is precisely
the same thing as a compact local generator. Moreover, for a given global object, being
a compact generator is a local property for the fpqc topology.

\begin{lem}\label{l5}
Let $T$ be a locally presentable dg-category with descent over a commutative simplicial ring $A$. 
\begin{enumerate}
\item An object 
$E \in T$ is a compact local generator if and only if it is a compact generator of the triangulated 
category $[T]$. 
\item Let $A \longrightarrow B$ be a fpqc covering of commutative simplicial rings. Then, a compact 
object $E$ of $T$ is a compact generator for $[T]$ if and only if its base change 
$E\otimes_{A}^{\mathbb{L}}B$ is a compact generator for $[T\otimes_{A}^{ct}B]$.
\end{enumerate}
\end{lem}

\textit{Proof:} $(1)$ Assume first that $E$ is a compact generator of $[T]$. We need to prove that for
any morphism of simplicial commutative rings $A \rightarrow A'$, the object 
$E\otimes_{A}^{\mathbb{L}}A' \in T\otimes_{A}^{ct}A'$, obtained by base change, is
a compact generator. For this we use the triangulated adjunction
$$-\otimes_{A}^{\mathbb{L}}A' : [T] \longleftrightarrow [T\otimes_{A}^{ct}A'] : f,$$
given on the one side from the adjunction dg-functor $T \longrightarrow T\otimes_{A}^{ct}A'$, and 
on the other side from the forgetful dg-functor 
$T\otimes_{A}^{ct}A' \longrightarrow T\otimes_{A}^{ct}A\simeq T$. The fact that these dg-functors are
adjoints to each others simply follows from the usual adjunction of locally presentable $A$-dg-categories
$$-\otimes_{A}^{\mathbb{L}}A' : \widehat{A} \longleftrightarrow \widehat{A'} : f,$$
tensored over $A$ by $T$ (or it also follows from the explicit construction 
of $\otimes_{A}^{ct}A'$ obtained by writing $T$ as a localisation of some $\widehat{T_{0}}$).
We claim moreover that the right adjoint $f$ is a conservative dg-functor. Indeed, 
write $T$ as a localisation of some dg-categories $\widehat{T_{0}}$, for 
$T_{0}$ a small dg-category over $A$, and along a small set of maps $S$ in $\widehat{T_{0}}$. 
Then the forgetful dg-functor $f$ is simply the forgetful dg-functor
$$\widehat{T_{0}\otimes_{A}^{\mathbb{L}}B} \longrightarrow \widehat{T_{0}}$$
restricted to $S$-local objects, which is clearly conservative. 
It is then a formal consequence
of this adjunction, and of the fact that $f$ is continuous and conservative, 
that $E\otimes_{A}^{\mathbb{L}}A'$ is a compact generator for
$[T\otimes_{A}^{ct}A']$. 

Conversely, is $E$ is a compact local generator then, taking the identity $A \rightarrow A$ gives that 
$E$ is a compact generator of $[T]$. \\

$(2)$ We again use the adjunction
$$-\otimes_{A}^{\mathbb{L}}B : [T] \longleftrightarrow [T\otimes_{A}^{ct}B] : f.$$
We let $T_{E}$ be the smallest sub-dg-category of $T$ containing $E$ and which 
is stable by colimits. By assumption on $E$ the inclusion dg-functor
$T_{E} \subset T$ induces an equivalence after base change
$$T_{E} \otimes_{A}^{ct}B \simeq T\otimes_{A}^{ct}B.$$
But both dg-categories $T_{E}$ and $T$ are locally presentable with descent over $A$,
as this follows from proposition \ref{p5} for the case of $T_{E}$ and by assumption for $T$. 
The theorem \ref{t1} implies that the base change $\s$-functor
$$-\otimes_{A}^{ct}B : \mathbb{D}g^{lp,desc}(A) \longrightarrow \mathbb{D}g^{lp,desc}(B)$$
is conservative. We thus have that the inclusion $T_{E} \subset T$ is an 
equivalence of dg-categories, or equivalently that $E$ is a compact generator of $[T]$. 
\hfill $\Box$ \\

The following proposition is a simple, but useful, observation.

\begin{prop}\label{p6'}
Let $F$ be any derived stack and $\alpha$ be a locally presentable dg-category over $F$.
Then, $\alpha$ is realizable by a dg-algebra $B_{\alpha}$ over $F$, if and only if the triangulated
category $D_{\alpha}(F)$ admits a compact local generator.
\end{prop}

\textit{Proof:} As we have defined the derived stack $\mathbb{D}g^{lp,desc}$, we can define
the derived stack $\mathbb{D}g^{lp,desc}_{*}$, of \emph{pointed locally presentable dg-categories with 
descent}.
This derived stack simply classifies locally presentable dg-categories with descent $T$, together with 
an object $E \in T$, as well as continuous dg-functors that preserves this object. We have 
a sub derived prestack $\mathbb{D}g^{lp,desc,cp} \subset \mathbb{D}g^{lp,desc}_{*}$ of 
pointed dg-categories for which the distinguished object is a compact generator. We observe that 
the morphism $\phi$ factorizes as
$$\mathbb{D}g\mathbb{A}lg \longrightarrow \mathbb{D}g^{lp,desc,cp} \longrightarrow \mathbb{D}g^{lp,desc},$$
because for a dg-algebra $B$, the dg-category $\widehat{B^{op}}$ possesses a canonical 
compact generator given by $B$ itself. Moreover, the induced morphism 
$$\mathbb{D}g\mathbb{A}lg \longrightarrow \mathbb{D}g^{lp,desc,cp}$$
is an equivalence of derived prestacks (so in particular the right hand side is
itself a stack). 

Now, given $\alpha : F \longrightarrow \mathbb{D}g^{lp,desc}$, a factorization
of $\alpha$ throught $\mathbb{D}g^{lp,desc}_{*}$ is precisely equivalent to the data of 
an object $E \in L_{\alpha}(F)$. Moreover, this object $E$ is a compact local generator 
exactly when this factorization lands in the substack $\mathbb{D}g^{lp,desc,cp} \subset 
\mathbb{D}g^{lp,desc}_{*}$. 
\hfill $\Box$ \\

The main existence statement concerning compact local generators is
the following theorem. 

\begin{thm}\label{t2}
Let $X$ be quasi-compact and quasi-separated derived scheme, and
$\alpha$ be a locally presentable dg-category with descent over $X$.
We assume that there is a fppf covering $X' \longrightarrow X$,
such that $D_{\alpha}(X')$ possesses a compact local generator. Then 
the triangulated category $D_{\alpha}(X)$ possesses a compact generator which 
is also a compact local generator. 
\end{thm}

The proof of the previous theorem will take us some time and will be given is the next paragraphes: it will be
a direct consequences of propostions \ref{p8}, \ref{p9} and \ref{p10} below. We however
already state the following main corollary.

\begin{cor}\label{c6}
Let $X$ be a quasi-compact and quasi-separated 
derived scheme, and $\alpha$ be a  locally presentable dg-category with descent over $X$. 
Then, $\alpha$ possesses a compact local generator locally for the fppf topology on $X$ if and
only if it is realized by a dg-algebra $B_{\alpha}$ on $X$. In particular we have natural isomorphisms
$$dBr(X) \simeq dBr_{cat}(X)\simeq H^{1}_{fppf}(X,\mathbb{Z})\times H^{2}_{fppf}(X,\mathbb{G}_{m}),$$
so any class in $H^{2}_{fppf}(X,\mathbb{G}_{m})$ is realizable by a derived 
Azumaya algebra over $X$. 
\end{cor}

\textit{Proof:} Assume first that $\alpha$ possesses a compact local generator locally for the
fppf topology on $X$, that is there exist a fppf covering $u : X' \longrightarrow X$, 
and an object $E\in D_{u^{*}(\alpha)}(X')$, which is a compact generator when 
restricted to any affine derived scheme over $X'$. Then the theorem \ref{t2} implies that 
$D_{\alpha}(X)$ possesses a compact local generator, which by proposition \ref{p6'}
implies that $\alpha$ is realized by a dg-algebra $B_{\alpha}$ over $X$.
Conversely, if $\alpha$ is realized by $B_{\alpha}$ then proposition \ref{p6'}
tells us that $D_{\alpha}(X)$ admits a compact local generator. \\

The second part of the corollary, concerning deraz algebras, is a formal consequence
of the first part. Indeed, we already know that 
the natural morphism
$$\phi : dBr(X) \longrightarrow dBr_{cat}(X),$$
is injective. The surjectivity of this morphism is equivalent to state that 
any locally dg-category $\alpha$ on $X$, which is locally (for the fppf topology) equivalent to the unit
dg-category, is realized by a deraz algebra over $X$. By the corollary, such a dg-category 
$\alpha$ is realized by some dg-algebra $B_{\alpha}$. This dg-algebra, is then locally on $X_{fppf}$, 
Morita equivalent to the unit dg-algebra, and thus is locally a deraz dg-algebra. The proposition 
\ref{p1} $(2)$ therefore implies that $B_{\alpha}$ is a deraz dg-algebra over $X$, which is
an antecedant of $\alpha$ by the morphism $\phi$. 

Finally, the isomorphism $dBr_{cat}(X)\simeq H^{1}_{fppf}(X,\mathbb{Z})\times H^{2}_{fppf}(X,\mathbb{G}_{m})$, 
has been already proved in corollary \ref{c5}. \hfill $\Box$ \\

\subsection{Gluing compact generators for the Zariski topology}

In this paragraph we start proving our theorem \ref{t2} by first proving the
following special case. 

\begin{prop}\label{p8}
Let $X$ be quasi-compact and quasi-separated derived scheme, and
$\alpha$ be a locally presentable dg-category with descent over $X$.
We assume that there is a Zariski covering $\{U_{i}\}$ of $X$,
such that each $D_{\alpha_{U_{i}}}(U_{i})$ possesses a compact local generator. Then
the triangulated category $D_{\alpha}(X)$ possesses a compact generator which is also 
a compact local generator.  
\end{prop}

\textit{Proof:} This is essentially the same proof as the proof of the existence of 
a compact generator for the quasi-coherent derived categories of quasi-compact and quasi-separated
schemes, as this is done in \cite{bv}. \\

We start by some local construction over affine derived schemes. Assume 
that $X:=\mathbb{R}\underline{Spec}\, A$, for
some commutative simplicial ring $A$, and that $\alpha$ is realized by a 
dg-algebra $B$ over $A$. Let $U\subset X$ be a quasi-compact open affine 
derived sub-scheme. The sub-scheme $U$ is given by 
certain elements $f_{1},\dots,f_{n} \in A_{0}$, as the union of the elementary
opens
$$U_{i}:=X_{f_{i}}=\mathbb{R}\underline{Spec}\, A[f_{i}^{-1}] \subset X.$$ 
We let $K=K(A,f)$ be the  simplicial $A$-module obtained by 
freely adding $1$-simplices $h_{i}$ to $A$ with the property
$$d_{0}(h_{i})=f_{i} \qquad d_{1}(h_{i})=0 \qquad \forall \, i.$$
Considered as an $A$-dg-module, $K$ is perfect, and thus a compact object in $D(A)$. It is
moreover a compact generator for the full triangulated sub-category
of $D(A)$ consisting of all $A$-dg-modules set-theorically supported
on the closed sub-scheme $X-U$ (i.e. $A$-dg-modules $E$ such that 
$E\otimes_{A}^{\mathbb{L}}A[f_{i}^{-1}]\simeq 0$ for all $i$), as shown by 
the following lemma. 

\begin{lem}\label{lp8}
The right orthogonal complement of $K$ in $D(A)$ is generated, by homotopy colimits and limits by 
the objects $A[f_{i}^{-1}]$.
\end{lem}

\textit{Proof of the lemma:} We denote by 
$D_{U}(A)$ the full sub-category of $D(A)$ generated by 
the $A[f_{i}^{-1}]$ and homotopy colimits and limits. It is the essential image of the
direct image (fully faithful) functor 
$$D(U) \longrightarrow D(X)\simeq D(A),$$
induced by the inclusion $U\hookrightarrow X$.

The $A$-dg-module $K$ can be written as a
tensor product
$$K\simeq \otimes^{\mathbb{L}}_{i}K_{i},$$
where, for each $i$, $K_{i}$ is the cone of the morphism
$\times f_{i} : A \longrightarrow A$. We write
$$K_{\neq j}:=\otimes^{\mathbb{L}}_{i\neq j}K_{i}.$$
The object $A[f_{i}^{-1}]$ is right orthogonal to $K_{i}$, and thus
to $K$ because of the formula for any $M\in D(A)$
$$\rh_{A}(K,M)\simeq \rh_{A}(K_{\neq i},\rh_{A}(K_{i},M)).$$
As $K$ is compact this implies that $D_{U}(A)$ is contained in the
right orthogonal complement of $K$.

Conversely, assume that $M$ is right orthogonal to $K$. There is a morphism
$$u : M\longrightarrow M',$$
with the property that $M'$ belongs to $D_{U}(A)$, and that the cone $C$ of $u$ restricts 
to zero over the open $U$ (take $M'=u_{*}u^{*}(M)$, where $u^{*}$ and $u_{*}$ are
the inverse and direct image functors between $D(U)$ and $D(X)$). 

As $M$ and $M'$ are both right orthogonal to $K$, so is $C$. We thus have
$$\rh_{A}(K,C)\simeq \rh_{A}(K_{n},\rh_{A}(K_{\neq n},C))\simeq 0.$$
This implies that $\rh_{A}(K_{\neq n},C)$ is naturally a dg-module over $A[f^{-1}_{n}]$, and that 
we have
$$\rh_{A}(K_{\neq n},C) \simeq \rh_{A}(K_{\neq n},C)\otimes_{A}^{\mathbb{L}}A[f^{-1}_{n}] 
\simeq \rh_{A}(K_{\neq n},C\otimes_{A}^{\mathbb{L}}A[f^{-1}_{n}])\simeq 0.$$
By a descending induction on $n$ we thus see that 
$$\rh_{A}(K_{1},C)\simeq 0.$$
Therefore, we have
$$C\simeq C\otimes_{A}^{\mathbb{L}}A[f_{1}^{-1}]\simeq 0,$$
and thus $M=M'$, and $M\in D_{U}(A)$. This finishes the proof of the lemma. \hfill $\Box$ \\

We now consider $B_{-K}:=B\otimes_{A}^{\mathbb{L}}K \in D(B)$, the free $B$-dg-module
over $K$. By construction $B_{-K}$ is a compact object in $D(B)$, which is a generator
for the full sub-category $L_{\alpha}(X,K)$, of $B$-dg-modules which are set theorically supported
on $X-U$ (that is restrict to zero over $U$). We denote by $L_{\alpha}(X,U)$ the full sub-category of 
$L_{\alpha}(X)=\widehat{B^{op}}$, consisting of $B$-dg-modules $E \in L_{\alpha}(X)$
with $K\otimes_{A}^{\mathbb{L}}E\simeq 0$. As
$K$ is a perfect $A$-dg-module this is equivalent to say that 
$\rh_{A-dg-mod}(K,E)\simeq 0$, that is that
$E$ is right orthogonal to $B_{-K}$ in $L_{\alpha}(X)$.  The sub-category 
$L_{\alpha}(X,U)$ therefore consists of all 
$B$-dg-modules whose underlying $A$-dg-module belongs to $D_{U}(A) \subset D(A)$. 

We now consider the natural inverse image dg-functor
$$L_{\alpha}(X) \longrightarrow L_{\alpha}(U),$$
precomposed by the natural inclusion $L_{\alpha}(X,U)\subset L_{\alpha}(X)$. This
provides a natural dg-functor 
$$L_{\alpha}(X,U) \longrightarrow L_{\alpha}(U).$$
By definition of $L_{\alpha}(U)$, and by the gluing lemma \ref{lglu} we have
$$L_{\alpha}(U)\simeq Holim L_{\alpha}(X_{f}),$$
where the limit is taken over all elementary open affines
$$X_{f}=\mathbb{R}\underline{Spec}\, A[f^{-1}] \subset U\subset X.$$
We claim that the induced dg-functor
$$L_{\alpha}(X,U) \longrightarrow L_{\alpha}(U) \simeq Holim L_{\alpha}(X_{f})$$
is an equivalence of dg-categories. This simply is the descent for quasi-coherent 
$B$-dg-modules for the Zariski topology, and can be either proven directly or 
deduced from the descent of quasi-coherent complexes as follows. We have natural 
equivalence of dg-categories
$$L_{\alpha}(X_{f}) \simeq \rh(B,L(X_{f})) 
\qquad L_{\alpha}(X,U) \simeq \rh(B,L(U)),$$
where the internal Hom's are relative to the $\s$-category of locally presentable
dg-categories over $A$, and $L(-)$ stands for $L_{\mathbf{1}}(-)$, the dg-category 
of quasi-coherent complexes. We therefore have a commutative square of
$A$-dg-categories
$$\xymatrix{
L_{\alpha}(X,U) \ar[r] \ar[d]_-{\sim} & Holim L_{\alpha}(X_{f}) \ar[d]^-{\sim} \\
\rh(B,L(U)) \ar[r] & Holim \rh(B,L(X_{f})).}$$
The vertical morphisms are equivalences. The bottom horizontal morphism is also an equivalence
as we have
$$L(U) \simeq Holim L(X_{f}),$$
because quasi-coherent complexes form a stack for the fpqc topology, and thus
for the Zariski topology (see e.g. \cite[Lem. 3.1]{tova}). 

The consequence of this first local discussion is the existence of 
a semi-orthogonal decomposition of compactly generated locally presentable dg-categories 
$$<B_{-K}> \subset L_{\alpha}(X) \longrightarrow L_{\alpha}(U),$$
where $<B_{-K}> \subset L_{\alpha}(X)$ is the full sub-dg-category generated
by $B_{-K}$, colimits and shifts. 
In particular, the localisation theorem of Thomason-Neeman \cite[Thm. 2.1]{ne}, and
\cite[Cor. 3.2.3]{bv}, applies. 
We therefore have that for any compact object $E \in L_{\alpha}(U)$, there
exists a compact object $E'\in L_{\alpha}(X)$, whose image by 
$$L_{\alpha}(X) \longrightarrow L_{\alpha}(U)$$
is $E\oplus E[1]$. \\

We now prove the proposition \ref{p8}, by induction on the number of open affines needed
to cover the derived scheme $X$. First of all if $X=\mathbb{R}\underline{Spec}\, A$ is affine, 
we have $L_{\alpha}(X)\simeq \widehat{B^{op}}$, for some $A$-dg-algebra $B$ realizing $\alpha$.
Then $B$, as a dg-module over itself is both, a compact generator and a compact local generator
thanks to lemma \ref{l5}.  

Assume now that the proposition ($(1)+(2)$) holds for all quasi-compact and quasi-separated derived schemes which are
union of $m$ derived affine schemes, for $m>0$.
Let $X$ be a quasi-compact derived scheme, 
covering by two open sub-schemes $U$ and $V$, with $V$ affine and $U$ 
covered by $m$ affine open sub-schemes. We let $W=U\cap V$. Lemma \ref{lglu} tells us that
we have a homotopy cartesian square 
of dg-categories
$$\xymatrix{
L_{\alpha}(X) \ar[r] \ar[d] & L_{\alpha}(U) \ar[d] \\
L_{\alpha}(V) \ar[r] & L_{\alpha}(W).
}$$
By the induction hypothesis, let $E_{U}$ a compact generator of $L_{\alpha}(U)$ which is
also a compact local generator.  By what we have seen in the affine case, 
the image of $E_{U}\oplus E_{U}[1]$ in $L_{\alpha}(W)$ lifts to a
compact object in $L_{\alpha}(V)$. Replacing $E_{U}$ by $E_{U}\oplus E_{U}[1]$ 
we simply assume that there is a compact object in $L_{\alpha}(V)$ whose image
in $L_{\alpha}(W)$ is equivalent to the image of $E_{U}$. By the homotopy cartesian square
above there is an object $E_{X}\in L_{\alpha}(X)$ whose image 
is equivalent to $E_{U}$ in $L_{\alpha}(U)$. Moreover, because of the same
homotopy cartesian of dg-categories the object $E_{X}$ has compact 
components in $L_{\alpha}(V)$, $L_{\alpha}(U)$ and $L_{\alpha}(W)$ and thus
is easily seen to be itself compact. 

The derived scheme $V$ is affine,  
and $W$ is a quasi-compact open in $V$. By what we have seen in the beginning of the proof
the dg-functor
$$L_{\alpha}(V) \longrightarrow L_{\alpha}(W)$$
is a localization obtained by killing a compact object $E_{0}$
of $L_{\alpha}(V)$ (denoted by $B_{-K}$ in our discussion). By construction, this
compact object goes to zero in $L_{\alpha}(W)$, and is a compact generator for
the full sub-dg-category of $L_{\alpha}(V)$ of objects sent to zero in $L_{\alpha}(W)$. The 
homotopy cartesian square above implies that $E_{0}$ naturally lifts to a compact
object $F_{X} \in L_{\alpha}(X)$, whose restriction to $U$ is zero. We let 
$$E:=E_{X}\oplus F_{X}.$$
We claim that $E$ is a compact generator of $L_{\alpha}(X)$, as well as 
a compact local generator. Indeed,  
an object $E'$ is right orthogonal to $E$ in $L_{\alpha}(X)$, if and only if
it is right orthogonal to $E_{X}$ and $F_{X}$. But, by construction, we have
$$L_{\alpha}(X)(F_{X},E')\simeq 0\times_{0}^{h}L_{\alpha}(V)(E_{0},E'_{|V})\simeq
L_{\alpha}(V)(E_{0},E'_{|V}).$$
Therefore, $E'$ is right orthogonal to $F_{X}$ if and only if its restriction $E'_{|V}$
to $V$ is right orthogonal to $E_{0}$. By definition of $E_{0}$ this
is in turn equivalent to the fact that $E'_{|V}$ is supported on the open $W\subset V$, or 
equivalently that for any object $F \in L_{\alpha}(V)$, the restriction morphism
$$L_{\alpha}(V)(F,E'_{|V}) \longrightarrow L_{\alpha}(W)(F_{|W},E'_{|W})$$
is a quasi-isomorphism. In particular, we must have
$$L_{\alpha}(X)(E_{X},E') \simeq L_{\alpha}(U)(E_{U},E'_{|U}) \times_{L_{\alpha}(W)(E_{|W},E'_{|W})}^{h}L_{\alpha}(V)(E_{|V},E'_{|V})\simeq L_{\alpha}(U)(E_{U},E'_{|U}).$$
This implies that if $E'$ is orthogonal to both, $E_{X}$ and $F_{X}$, 
then we have $E'_{|U}=0$, thus $E'_{|W}=0$ and therefore $E'_{|V}=0$ because
$E'_{|V}$ is supported on $W$. This shows that $E'$ must be zero, and that 
$E$ is indeed a compact generator for $L_{\alpha}(X)$. 

By construction, $E$ is a compact generator such that $E_{|U}\simeq E_{U}$
is a compact local generator of $L_{\alpha}(U)$. Moreover, 
$E_{|V}$ is a compact object in $L_{\alpha}(V)$ with the following
two properties:
\begin{enumerate}
\item its restriction to each open affine $V_{f} \subset W$ is 
a compact generator of $L_{\alpha}(V_{f})$
\item $E_{|V}$ contains $E_{0}$ as a direct factor. 
\end{enumerate}

We claim that these two properties, together with the fact that $V$ is affine
implies that $E_{|V}$ is a compact generator of $L_{\alpha}(V)$. Indeed, let 
us write $V\simeq \mathbb{R}\underline{Spec}\, A$, and $W=\cup V_{f_{i}}$, for
a finite number of elements $f_{i}\in A_{0}$. Then, if an object 
$F$ is right orthogonal to $E_{|V}$ in $L_{\alpha}(V)$, by the property
$(2)$ above the object $F$ must be supported on $W$. Moreover, we have
$$L_{\alpha}(V)(E_{|V},F) \otimes_{A}^{\mathbb{L}}A[f^{-1}] \simeq
L_{\alpha}(V_{f})(E_{|V_{f}},F_{|V_{f}}).$$
Therefore, if $F$ is right orthogonal to $E_{|V}$ we must have 
$F_{|V_{f}}=0$, and as $F$ is supported on $W$ this implies that $F\simeq 0$. 
\hfill $\Box$ \\

\subsection{Gluing compact generators for the fppf topology}

We will proceed in two steps, using our previous result about gluing compact generators for
the Zariski topology \ref{p8}. We start by proving that existence of compact generators locally for
the \'etale topology implies the existence of compact generators locally for the 
Zariski topology. Then, we prove that existence of compact generators locally for the fppf
topology implies the existence of compact generators locally for the \'etale topology. In both cases we
will use the follwing lemma. 

\begin{lem}\label{l6}
Let $A \longrightarrow B$ be a finite and faithfully flat morphism between commutative simplicial rings, and
$T$ be a locally presentable dg-category with descent over $A$. Then $[T]$ possesses a compact generator
if and only if $[T\otimes_{A}^{ct}B]$ does so.
\end{lem}

\textit{Proof:} We use the adjunction 
$$-\otimes_{A}^{\mathbb{L}}B : [T] \longrightarrow [T\otimes_{A}^{ct}B] : f.$$
The right adjoint $f$ is easily seen to be conservative, as we already have seen 
this during the proof of lemma \ref{l5} $(1)$. During the same proof we also have seen that 
$E\otimes_{A}^{\mathbb{L}}B$ is a compact generator of $[T\otimes_{A}^{ct}B]$ when
$E$ is a compact generator of $[T]$. 

Conversely, assume that $E$ is a compact generator of $[T\otimes_{A}^{ct}B]$.
We claim that $f(E)$ is a compact generator of $[T]$. For this we use that 
$f$ does itself have a right adjoint 
$$f^{!} : [T] \longrightarrow [T\otimes_{A}^{ct}B].$$
This right adjoint is a priori not a continuous functor, but $B$ being finite
flat over $A$ implies that it is continuous. Indeed, the direct image dg-functor
$$p : \widehat{B} \longrightarrow \widehat{A},$$
possesses a right adjoint 
$$p^{!} : \mathbb{R}\underline{Hom}_{A}(B,-) \otimes^{\mathbb{L}} - : 
\widehat{A} \longrightarrow \widehat{B}.$$
This right adjoint commutes with direct sums because $B$ is finite flat over $A$, thus
projective and of finite type (see \cite[Lem. 2.2.2.2]{hagII}), and thus compact as a $A$-dg-module.
In particular, if $B^{\vee}:=\mathbb{R}\underline{Hom}_{A}(B,A)$ denotes the $A$-dg-module dual of $B$
then we have
$$p(p^{!}(M)) \simeq B^{\vee}\otimes^{\mathbb{L}} M,$$
for all $M \in D(A)$. 
Tensoring the continuous dg-functor $p^{!}$ with $T$ gives a continuous dg-functor
$$f^{!} : [T] \longrightarrow [T\otimes_{A}^{ct}B],$$
which is right adjoint to the direct image dg-functor $f$. It follows formally from the
existence of $f^{!}$ that $f$ preserves compact objects, and thus that 
$f(E)$ is compact in $[T]$. Moreover, $B^{\vee}$ is
a faithfully flat $A$-dg-module, and thus $f\circ f^{!}\simeq B^{\vee} \otimes^{\mathbb{L}} -$ is
a conservative functor. As
$f$ is conservative, it follows that so is $f^{!}$. This, again formally, implies that 
$f(E)$ is a generator of $[T]$.  
\hfill $\Box$ \\

We now show that local compact generators for the \'etale topology can be glued to 
get compact local generators for the Zariski topology. The proof of the proposition we give below is an 
adaptation of Gabber's proof of the existence of Azumaya algebras over affine schemes (see \cite{gab}). \\

\begin{prop}\label{p9}
Let $A \longrightarrow B$ be a \'etale covering of commutative simplicial rings, and
$T$ be a locally presentable dg-category over $A$. Then $[T]$ possesses a compact generator
if and only if $[T\otimes_{A}^{ct}B]$ does so.
\end{prop}

\textit{Proof:} By lemma \ref{l5} the base change of a compact generator for $[T]$ is a compact
generator for $[T\otimes_{A}^{ct}B]$. So let $E$ be a compact generator
for $[T\otimes_{A}^{ct}B]$. We are allowed to refine our \'etale covering $A \rightarrow B$. Moreover, 
using proposition \ref{p8} we are 
also allowed to work locally for the Zariski topology on $A$. We can therefore
assume that our \'etale morphism $A \longrightarrow B$ is a standard \'etale morphism: there exists
a commutative diagram of commutative simplicial rings 
$$\xymatrix{
A \ar[r]^-{p} \ar[dr] & C \ar[d]^-{j} \\
 & B,}$$
where $j$ is a Zariski open immersion (that can be taken to be principal, that is of the
form $C[c^{-1}]$ for some $c\in \pi_{0}(C)$), and where $C$ is equivalent to
$A[X]/p$, for $p$ a monic polynomial in $\pi_{0}(A)$. Before going further let us 
give more details about $C=A[X]/p$. Here $A[X]$ is simply the free commutative simplicial $A$-algebra, and
$p$ is a monic polynomial over the ring $\pi_{0}(A)$
$$p=x^{d}+\sum_{i}x^{i}a_{i} \in \pi_{0}(A)[X].$$
As homotopy classes of morphisms $A[X] \rightarrow A'$ are in one-to-one correspondence
with $\pi_{0}(A')$, we can represent $p$ by a morphism of commutative simplicial rings
$$p : A[X] \longrightarrow A[X]$$
well defined in the homotopy category of commutative simplicial rings. By definition
$A[X]/p$ is the commutative simplicial ring fitting in the following homotopy cocartesian 
square
$$\xymatrix{
A[X] \ar[r]^-{p} \ar[d] & A[X] \ar[d] \\
A \ar[r] & A[X]/p.
}$$
We now consider the following morphism of commutative rings
$$\mathbb{Z}[X_{1},\dots,X_{d}]^{\Sigma_{d}} \hookrightarrow \mathbb{Z}[X_{1},\dots,X_{d}].$$
The left hand side is isomorphic to a polynomial ring 
$\mathbb{Z}[T_{1},\dots,T_{d}]$, where the identification is by sending 
$T_{i}$ to the $i$-th elementary function on the set of variables $X_{j}$
$$T_{i} \mapsto \phi_{i}(X)=(-1)^{d-i}\sum_{n_{1}+\dots+n_{d}=i} X_{1}^{n_{1}}.\dots X_{d}^{n_{d}}.$$
The polynomial $p$ is determined by a morphism of commutative simplicial rings
(again, well defined in the homotopy category $Ho(s\z-CAlg)$)
$$\mathbb{Z}[T_{1},\dots,T_{d}] \longrightarrow A,$$
sending $T_{i}$ to the element $a_{i-1} \in \pi_{0}(A)$. We form the homotopy push-out
$$\xymatrix{
\mathbb{Z}[T_{1},\dots,T_{d}] \ar[r] \ar[d] & A \ar[d]\\
\mathbb{Z}[X_{1},\dots,X_{d}] \ar[r] & C'.}$$
As $\mathbb{Z}[X_{1},\dots,X_{d}]$ is finite and free as a $\mathbb{Z}[T_{1},\dots,T_{d}]$-module, 
$A \rightarrow C'$ is a finite and flat morphism of commutative simplicial rings, coming equiped 
with a natural action of $\Sigma_{d}$.  
There exists a natural commutative diagram, in $Ho(s\z-CAlg)$
$$\xymatrix{
A \ar[r] \ar[dr] & C \ar[d] \\
 & C'.}$$
To see this, let $q\in A_{0}[X]$ be a lift of 
$p\in \pi_{0}(A)[X]$.
Then, the datum of 
a morphism in $Ho(s\z-CAlg)$
$$C \longrightarrow C',$$
is equivalent to the data of a pair $(c,h)$, where
$c\in C'_{0}$, $h\in C'_{1}$ with
$$d_{0}(h)=q(c) \qquad d_{1}(h)=0,$$ 
two such data $(x,h)$ and $(y,k)$ considered as being equivalent
if they are homotopic (in a rather obvious meaning). The morphism 
of commutative simplicial $A$-algebras $C \longrightarrow C'$
is then obtained by chosing, say, the image $x_{1}$, of $X_{1}$ in $C'_{0}$, 
which satisfies by construction that $q(x_{1})$ is naturall homotopic to zero.
In fact, by the construction of $C'$, the image of the polynomial $q$ in $C'_{0}$
comes equiped with a natural homotopy to the polynomial 
$\prod_{i}(X-x_{i})$. 

We now base change the Zariski open $j : C \longrightarrow B$ to get 
another Zariski open
$$C' \longrightarrow B':=C'\otimes_{C}^{\mathbb{L}}B.$$
By hypothesis $[T\otimes_{A}^{ct}B]$ has a compact generator, which by 
lemma \ref{l5} implies that $[T\otimes_{A}^{ct}B']$ also possesses a compact generator. 
By using the action of $\Sigma_{d}$ on $C'$, by autormophisms over $A$, 
we obtain a finite family of Zariski open 
$$C' \longrightarrow B'_{\sigma} \qquad \sigma\in \Sigma_{d}.$$
This family of Zariski open has the property that each 
triangulated category $[T\otimes_{A}^{ct}B'_{\sigma}]$ possesses a compact generator. 

We now let $T':=T\otimes_{A}^{ct}C'$, which is a locally presentable dg-category 
with descent over $C'$. It is such that each $T'\otimes_{C'}^{ct}B'_{\sigma}$
has a compact generator. Moreover, the Zariski open $B'_{\sigma}$ form 
a Zariski open covering of $C'$. This can be checked on the level 
of geometric points, and thus becomes a non simplicial statement. Indeed, 
let $u : C' \longrightarrow L$ be morphism with $L$ an algebraically closed field. 
By definition of $C'$ 
the morphism $u$ is determined by the data of elements $\alpha_{i} \in L$ 
such that $p(X)=\prod_{i}(X-\alpha_{i})$, where $p$ is considered as
a polynomial with coefficients in $L$ using the morphism $A \rightarrow C' \rightarrow L$. 
The point $A \rightarrow L$ lifts to a point $B \rightarrow L$ (because $A \rightarrow B$
is an \'etale covering), and thus to a point $v : C \rightarrow L$. The morphism 
$v$ is determined by $\alpha \in L$ which is a root of $p(X)$, which must be equal
to $\alpha_{i}$ for some $i$. If we chose $\sigma \in \Sigma_{d}$ with
$\sigma^{-1}(i)=1$, then we have a commutative diagram
$$\xymatrix{
C \ar[r] \ar[d] & B \ar[r] & L \\
C' \ar[d]_-{\sigma} & & \\
C', \ar[rruu]_-{u} & & }$$
where the morphism at the top is the point $v$. This implies that $\sigma^{-1}$
sends the Zariski open $B'$ of $C'$ to a new Zariski open containing the point $u$. 
This finishes to show that $B'_{\sigma}$ form a Zariski covering of $C'$ such that 
$[T'\otimes_{A}^{ct}B'_{\sigma}]$ possesses a compact generator. 
We then conclude using proposition \ref{p8} that $T'$ has a compact generator, which by 
lemma \ref{l6} implies that $T$ also has a compact generator.
\hfill $\Box$ \\

We now finish the last step of the proof of our theorem \ref{t2} by proving that the local compact
generators for the fppf topology can be glued to a local 
compact generator for the Zariski topology, and thus to a global 
compact generator by our previous proposition \ref{p8}.

\begin{prop}\label{p10}
Let $A \longrightarrow B$ be a fppf covering of commutative simplicial rings, and
$T$ be a locally presentable dg-category over $A$. Then $[T]$ possesses a compact generator
if and only if $[T\otimes_{A}^{ct}B]$ does. 
\end{prop}

\textit{Proof:} As for the proofs of \ref{p8} and \ref{p9} one implication is already known. So let 
$E$ be a compact generator for $[T\otimes_{A}^{ct}B]$. From \cite[Lem. 2.14]{to3} 
(only the easiest affine case is 
needed here) we know the existence of 
a smooth covering $p : X \rightarrow \mathbb{R}\underline{Spec}\, A$, 
classifying quasi-smooth quasi-sections of $A \rightarrow B$. 
Taking local sections of $p$ for the \'etale topology, we do get the existence of 
a commutative diagram of commutative simplicial rings
$$\xymatrix{
A \ar[r] \ar[d] & B \ar[d] \\
A' \ar[r] & B',}$$
where $A \rightarrow A'$ is an \'etale covering, and $A' \rightarrow B'$
is a finite flat morphism of a certain rank $m$ (and we have moreover that $B \rightarrow B'$ is
a quasi-smooth morphism, but we will not use this additional property). Therefore, 
using lemma \ref{l5} and the hypothesis, $[T\otimes_{A}^{ct}B']$ possesses a compact generator. By the lemma
\ref{l6} applied to the finite flat map $A' \rightarrow B'$, 
we have that $[T\otimes_{A}^{ct}A']$ also have a compact generator, which by the proposition 
\ref{p9} implies that $T$ has a compact generator. This finishes the proof of the proposition.
\hfill $\Box$ \\

\section{Applications and complements}

\subsection{Existence of derived Azumaya algebras}

In this subsection we refine a bit the existence statement for derived azumaya algebras corollary 
\ref{c6}, 
by also considering certain algebraic spaces and stacks. 

\begin{cor}\label{c7}
Let $X$ be a derived stack. Then, the natural map
$$dBr(X) \longrightarrow H^{1}_{et}(X,\mathbb{Z}) \times H^{2}_{et}(X,\mathbb{G}_{m})$$
is bijective in each of the following cases.
\begin{enumerate}
\item $X$ is a quasi-compact and quasi-separated derived scheme.
\item $X$ is a smooth and separated algebraic space of finite type over a field of characteristic 
zero $k$.
\item $X$ is a quasi-compact, separated and derived Deligne-Mumford stack 
whose coarse moduli space is a derived scheme.
\end{enumerate}
\end{cor}

\textit{Proof:} $(1)$ This is our corollary \ref{c6}. For $(2)$, let $p : X' \longrightarrow X$
be a proper birational morphism with $X'$ smooth. As $X$ is smooth
we have $H^{1}_{et}(X,\mathbb{Z})=0$. Let $\alpha\in H^{2}_{et}(X,\mathbb{G}_{m})$, that 
we represent by a locally presentable dg-category over $X$. 
We have a pull-back 
functor
$$p^{*} : L_{\alpha}(X) \longrightarrow L_{\alpha}(X').$$
This functor possesses a right adjoint 
$$p_{*} : L_{\alpha}(X') \longrightarrow L_{\alpha}(X),$$
which is fully faithful, as this easily follows from the fact that $X$ has rational singularities
(and thus $\mathbb{R}p_{*}(\mathcal{O}_{X'})\simeq \mathcal{O}_{X}$). Therefore, 
$L_{\alpha}(X)$ becomes a direct summand (in the sense of orthogonal decomposition) 
of $L_{\alpha}(X')$. It is easy to deduce from this that the image by $p_{*}$ of 
a compact generator of $L_{\alpha}(X')$ is a compact generator of $L_{\alpha}(X)$. In the same
way, the image by $p_{*}$ of a compact local generator of $L_{\alpha}(X')$ is a compact
local generator of $L_{\alpha}(X)$, by using a base change to a affine scheme $U$ \'etale over $X$. 
As we know that such compact local generator
exists for $L_{\alpha}(X')$ we see that $L_{\alpha}(X)$ has a compact local generator, and thus
that $\alpha$ is represented by a derived Azumaya algebra over $X$. 

We now consider the case $(3)$. We let $\pi : X \longrightarrow M$ be the projection of 
$X$ to its coarse moduli space $M$, which by assumption is a quasi-compact and separated 
derived scheme. We let $\alpha \in H^{1}_{et}(X,\mathbb{Z})\times H^{2}_{et}(X,\mathbb{G}_{m})$, 
that we represent by a locally presentable dg-category $\alpha$ over $X$. We then consider
$\pi_{*}(\alpha)$, which is a locally presentable dg-category over $M$ defined by 
$$(U\mapsto M) \mapsto L_{\alpha}(\pi^{-1}(U)).$$
We have 
$$L_{\alpha}(X)\simeq L_{\pi_{*}(\alpha)}(M).$$
Now, locally for the \'etale topology on $M$, the derived stack $X$ is
equivalent to a quotient stack $[U/G]$, for $G$ a finite group acting on 
$U=\mathbb{R}\underline{Spec}\, A$ a derived affine scheme. More precisely, 
there is an \'etale covering $\{V_{i} \longrightarrow M\}$, and pull-back diagramms
$$\xymatrix{
[U_{i}/G_{i}] \ar[r] \ar[d] & X \ar[d] \\
V_{i} \ar[r] & M,}$$
with $[U_{i}/G_{i}]$ as above ($U_{i}\simeq \mathbb{R}\underline{Spec}\, A_{i}$ and $G_{i}$
a finite group acting on $U_{i}$). 
If we let $B_{i}:=A_{i}[G_{i}]$, the twisted group simplicial algebra of $A_{i}$ by $G_{i}$, then 
we have
$$L_{\alpha}([U_{i}/G_{i}])\simeq \widehat{B_{i}^{op}}.$$
This shows that the locally presentable dg-category $\pi_{*}(\alpha)$ on $M$ 
possesses compact generators locally for the \'etale topology on $M$. The theorem \ref{t2}
tells us that $L_{\pi_{*}(\alpha)}(M)\simeq L_{\alpha}(X)$ possesses a compact local generator, 
and thus that $\alpha$ is represented by a deraz algebra over $X$. \hfill $\Box$ \\

The previous corollary implies the existence of derived Azumaya algebras 
in the two following concrete examples. 

The first situation is the famous
example of Mumford, of a normal local $\mathbb{C}$-algebra $A$ of dimension $2$, with
a non-torsion class $\alpha \in H^{2}_{fppf}(Spec\, A,\mathbb{G}_{m})$ (see \cite{gr}). This class
is then realized by an $A$-dg-algebra $B$, which, as $\alpha$ is non-torsion, is not 
Morita equivalent to an Azumaya non-derived algebra. It seems an interesting question to 
write down explicit derived Azumaya algebras, say by generators and relations, representing
the class $\alpha$. 

For our second example, let $X$ be any smooth and proper scheme
over some ring of characteristic zero $k$. We assume that there is a class $\alpha \in H^{2}(X,\mathcal{O})$ 
which is non-zero. We see the class $\alpha$ as a first order deformation of the
trivial class in $H^{2}_{fppf}(X,\mathbb{G}_{m})$, and thus consider it as a class
in $H^{2}_{fppf}(X[\epsilon],\mathbb{G}_{m})$ that restricts to zero
in $H^{2}_{fppf}(X,\mathbb{G}_{m})$. The corollary \ref{c7} tells us that the class $\alpha$ is then
realized by a derived Azumaya algebra $B$ over $X[\epsilon]$, whose restriction to $X$ is 
Morita equivalent to $\mathcal{O}$. Again, as the class $\alpha$ is non-zero it is non-torsion, and
the deraz algebra $B$ is thus not Morita equivalent to an Azumaya algebra. \\

An important corollary of the Cor. \ref{c7}, or rather of the proof of point $(2)$, 
is the following existence statement of compact generator for derived categories
of Deligne-Mumford stacks. 

\begin{cor}\label{c7'}
Let $X$ be a quasi-compact and separated Deligne-Mumford stack whose 
coarse moduli space is a scheme. Then the derived category 
$D(X)$ of quasi-coherent sheaves on $X$ admits a compact generator.
\end{cor}

\subsection{Localization theorem for twisted K-theory}

We present here an application of the existence of derived azumaya algebras to 
twisted K-theory, by deducing the existence of a localization exact triangle 
for the K-theory of twisted perfect complexes. For this we need to introduce some
notations. Let $\alpha$ be a locally presentable dg-category over some derived stack $X$. 
The compact object in $L_{\alpha}(X)$ form a small triangulated dg-category $L_{\alpha}(X)_{c}$, 
of which we can construct its non-connective K-theory spectrum (see \cite{schl})
$$K_{\alpha}(X):=K(L_{\alpha}(X)_{c}).$$
For instance, when $\alpha=\mathbf{1}_{X}$ is the unit dg-category over $X$, then 
$K_{\alpha}(X)$ is nothing else than the K-theory spectrum of perfect complexes
over the derived stack $X$. 

We will also use the following notation: if $Y$ is a closed sub-stack of $X$ we note
$L_{\alpha}(X,Y)$  the full sub-dg-category of $L_{\alpha}(X)$ of all objects
that restricts to zero over the open complement $X-Y$. We thus have a sequence of
dg-categories
$$L_{\alpha}(X,Y)_{c} \longrightarrow L_{\alpha}(X)_{c} \longrightarrow L_{\alpha}(X-Y)_{c},$$
which in general is only exact on the left. Exactness on the right (i.e. 
the fact that any object in  $L_{\alpha}(X-Y)_{c}$ is a retract of a the image of
an object in $L_{\alpha}(X)_{c}$)
is precisely what is needed in order to apply Thomason's localization theorem to obtain 
an exact triangle on the level of K-theory spectra. 

\begin{cor}\label{c8}
Let $X$ be a separated and quasi-compact derived Deligne-Mumford stack whose coarse moduli 
space is a derived scheme. Let $\alpha$ be a locally presentable dg-category over $X$, which locally 
for the \'etale topology on $X$ possesses a compact generator. Let $Y$ be a finitely presented
closed sub-stack 
in $X$, then the
sequence of dg-categories 
$$L_{\alpha}(X,Y)_{c} \longrightarrow L_{\alpha}(X)_{c} \longrightarrow L_{\alpha}(X-Y)_{c}$$
induces an exact triangle on the corresponding (non-connective) $K$-theory 
spectra 
$$K_{\alpha}(X,Y) \longrightarrow K_{\alpha}(X) \longrightarrow K_{\alpha}(X-Y).$$
\end{cor}

\textit{Proof:} We have to prove that the sequence of dg-categories
$$L_{\alpha}(X,Y)_{c} \longrightarrow L_{\alpha}(X)_{c} \longrightarrow L_{\alpha}(X-Y)_{c}$$
is exact, that is the right hand side identifies itself naturally with 
the quotient dg-category $L_{\alpha}(X)_{c}/L_{\alpha}(X,Y)_{c}$, where this
quotient has to be understood in the homotopy theory of small dg-categories
up to Morita equivalences (see \cite{tabmor}). 

Let $\pi : X \longrightarrow M$ be the projection of $X$ to its coarse moduli space, which by 
assumption is a quasi-compact and separated derived scheme. By pushing $\alpha$ down to $M$, 
and noticing that $L_{\alpha}(X)\simeq L_{\pi_{*}(\alpha)}(M)$ (as well as
with versions with supports), 
we see that we can in fact assume that $X=M$ and that $X$ is in fact a derived scheme. By the theorem 
\ref{t2} $\alpha$ is realized by a sheaf of dg-algebra $B$ over $X$ with quasi-coherent cohomology. 
The dg-categories
$L_{\alpha}(X)$ and $L_{\alpha}(X-Y)$ are then naturally equivalent to the dg-categories
of $B$-dg-modules over $X$ with quasi-coherent cohomology, and  $B$-dg-modules over $X-Y$
with quasi-coherent cohomology. 
The pull-back functor
$$j^{*} : L_{\alpha}(X) \longrightarrow L_{\alpha}(X-Y)$$
then possesses a fully faithful right adjoint
$$j_{*} : L_{\alpha}(X-Y) \longrightarrow L_{\alpha}(X),$$
simply induced by push-forward on the level of complexes with quasi-coherent cohomology.
This implies that we have a semi-orthogonal decomposition
of locally presentable dg-categories
$$L_{\alpha}(X,Y) \longrightarrow L_{\alpha}(X) \longrightarrow L_{\alpha}(X-Y).$$
We know by our theorem \ref{t2} that $L_{\alpha}(X)$ and 
$L_{\alpha}(X-Y)$ are compactly generated, and it is easy to see that 
the dg-functor $j^{*}$ preserves compact objects. The kernel $L_{\alpha}(X,Y)$ is also 
compactly generated as another application of the theorem \ref{t2} to $\alpha_{Y}$, the kernel 
of the natural morphism 
$$\alpha \longrightarrow j_{*}j^{*}(\alpha),$$
of locally presentable dg-categories over $X$, which we have seen to have compact generators
locally for the Zariski topology on $X$ (it was denoted $B_{-K}$ in the proof of proposition 
\ref{p8}).
It follows formally from this that the sequence induced on compact objects
$$L_{\alpha}(X,Y)_{c} \longrightarrow L_{\alpha}(X)_{c} \longrightarrow L_{\alpha}(X-Y)_{c}$$
is an exact sequence in the Morita theory of small dg-categories. 
\hfill $\Box$ 

\subsection{Direct images of smooth and proper categorical sheaves}

In \cite{tove1} we have introduced a certain class of categorical sheaves, as categorified
versions of vector bundles. By definition these categorical sheaves are 
locally presentable dg-categories over schemes, which locally for the Zariski topology
are smooth and proper in the sense of our proposition \ref{p2}. The result of the present paper have
several important consequences concerning functoriality properties for these
categorical sheaves. The following result states various properties of 
direct images of categorical sheaves along certain kind of morphisms. Probably the most 
important case is that smooth and proper locally presentable dg-categories are stable by 
direct images along smooth and proper morphisms. This gives a general tool in order
to construct smooth and/or proper dg-categories by integration of smooth and/or proper
dg-categories over smooth and/or proper schemes. \\

\begin{cor}\label{c9}
Let $X$ be a quasi-separated and quasi-compact derived scheme, equiped with 
a morphism $X \longrightarrow \mathbb{R}\underline{Spec}\, k$, for some 
commutative simplicial ring $k$.  
Let $\alpha$ be a locally presentable dg-category over $X$, and assume
that there is an \'etale covering $U=\mathbb{R}\underline{Spec}\, A \longrightarrow X$, 
such that $L_{\alpha}(U)$ possesses a compact generator, and therefore is
equivalent to $\widehat{B^{op}}$ for some $A$-dg-algebra $B$. 
\begin{enumerate}
\item The dg-category $L_{\alpha}(X)$ possesses a compact generator which is also 
a compact local generator.
\item If $B$ is smooth as
an $A$-dg-algebra, and if $X$ is smooth over $k$, then 
the compactly generated dg-category $L_{\alpha}(X)$ is smooth over $k$.
\item If $B$ is a proper $A$-dg-algebra (in the sense of proposition \ref{p2}), and if
$X$ is proper and flat over $k$, 
then the compactly generated dg-category $L_{\alpha}(X)$ is proper over $k$.
\item If $B$ is saturated as an $A$-dg-algebra (i.e. smooth and proper), 
and if $X$ is smooth and proper
over $k$, then 
$L_{\alpha}(X)$ is saturated as a $k$-dg-category.
\end{enumerate}
\end{cor}

\textit{Proof:} $(1)$ This is theorem \ref{t1}.  \\

$(2)$ By $(1)$ there is a compact local generator
in $L_{\alpha}(X)$, and thus $L_{\alpha}(X)$ can be written as 
the dg-category of dg-modules with quasi-coherent cohomology over a dg-algebra
$B_{\alpha}$ over $X$. We thus have
$$L_{\alpha}(X)\simeq Holim_{\mathbb{R}\underline{Spec}\, A' \rightarrow X}\widehat{B_{A'}^{op}}$$
where $B_{A'}$ is the $A'$-dg-algebra which is the restriction of $B_{\alpha}$ over
$\mathbb{R}\underline{Spec}\, A'$. As $X$ is quasi-compact and quasi-separated
this limit can be taken to be a finite limit, corresponding to 
a finite hypercoverings of $X$ by affine open sub-schemes. 

As the compact local generator is also 
a compact generator, we have 
$$L_{\alpha}(X)\simeq \widehat{B_{X}^{op}},$$
where $B_{X}:=Holim B_{A'}$ is the dg-algebra of global sections of $B$ over $X$.
We need to prove that $B_{X}$ is a smooth dg-algebra over $k$. 

\begin{lem}\label{lc9}
With the notation as above, the natural dg-functor
$$\widehat{B_{X}} \longrightarrow Holim_{\mathbb{R}\underline{Spec}\, A' \rightarrow X} \widehat{B_{A'}}$$
is a quasi-equivalence.
\end{lem}

\textit{Proof of the lemma:} We start by the quasi-equivalence
$$\widehat{B_{X}^{op}} \longrightarrow Holim_{\mathbb{R}\underline{Spec}\, A' \rightarrow X} \widehat{B_{A'}^{op}}.$$
All the dg-categories envolved in this equivalence
are compactly generated dg-categories. Moreover, all the dg-functors envolved 
in the above diagram are compact dg-functors, that is dg-functors preserving
compact objects. As the homotopty limit on the right hand side is finite, the
compact objects of $Holim_{\mathbb{R}\underline{A'} \rightarrow X} \widehat{B_{A'}^{op}}$
form a dg-category equivalent to $Holim_{\mathbb{R}\underline{Spec}\, A' \rightarrow X} 
(\widehat{B_{A'}^{op}})_{c}$, 
where $(\widehat{B_{A'}^{op}})_{c}$ is the full sub-dg-category of compact objects in 
$\widehat{B_{A'}^{op}}$. We thus have a quasi-equivalence
$$(\widehat{B_{X}^{op}})_{c} \simeq
Holim_{\mathbb{R}\underline{Spec}\, A' \rightarrow X} (\widehat{B_{A'}^{op}})_{c}.$$
Passing to opposite dg-categories we find a new quasi-equivalence of dg-categories
$$(\widehat{B_{X}})_{c}\simeq (\widehat{B_{X}^{op}})_{c}^{op} \simeq
(Holim_{\mathbb{R}\underline{Spec}\, A' \rightarrow X} (\widehat{B_{A'}^{op}})_{c})^{op}\simeq
Holim_{\mathbb{R}\underline{Spec}\, \underline{A'} \rightarrow X} (\widehat{B_{A'}})_{c}.$$
Now, again because the homotopy limit is finite we have 
$$(Holim_{\mathbb{R}\underline{Spec}\, \underline{A'} \rightarrow X} (\widehat{B_{A'}}))_{c} \simeq 
Holim_{\mathbb{R}\underline{Spec}\, \underline{A'} \rightarrow X} (\widehat{B_{A'}})_{c}.$$
The conclusion is then that the natural dg-functor
$$\phi : \widehat{B_{X}} \longrightarrow Holim_{\mathbb{R}\underline{Spec}\, A' \rightarrow X} \widehat{B_{A'}},$$
induces a quasi-equivalence on the full sub-dg-categories of compact objects. By the theorem
\ref{t1}, 
applied to the locally presentable dg-category over $X$ realized by the dg-algebra $B_{\alpha}^{op}$, 
the dg-category $Holim_{\mathbb{R}\underline{Spec}\, A' \rightarrow X} \widehat{B_{A'}}$
is compactly generated. This implies that $\phi$ is thus a quasi-equivalence. \hfill $\Box$ \\

According to the lemma \ref{lc9}, we have
$$\widehat{B_{X}\otimes_{k}^{\mathbb{L}}B_{X}^{op}} \simeq
Holim_{\mathbb{R}\underline{Spec}\, A' \rightarrow X} 
\widehat{B_{A'}\otimes_{k}^{\mathbb{L}}B_{A'}^{op}}\simeq 
L_{\alpha^{\vee}\otimes^{ct}_{k}\alpha}(X\times_{k} X),$$
where $\alpha^{\vee}\otimes^{ct}_{k}\alpha$ is the external product of $\alpha$ by its dual 
$\alpha^{\vee}$ (realized by $B_{\alpha}^{op}$)
on the derived stack $X\times_{k}X$. In other words, 
$\widehat{B_{X}\otimes_{k}^{\mathbb{L}}B_{X}^{op}}$ is naturally equivalent to 
the dg-category of dg-modules over $B_{\alpha}\otimes_{k}^{\mathbb{L}}B^{op}_{\alpha}$, which is
itself a dg-algebra over $X\times_{k}X$.
The diagonal bi-dg-module $B_{X}$ in 
$\widehat{B_{X}\otimes_{k}^{\mathbb{L}}B_{X}^{op}}$ corresponds by these equivalences
to $B_{\alpha}$, considered as a dg-module over $B_{\alpha}\otimes_{k}^{\mathbb{L}}B_{\alpha}^{op}$, 
which is itself a dg-algebra over $X\times_{k}X$. As $X$ is strongly of finite presentation
over $k$, its quasi-coherent cohomology commutes with direct sums. As a consequence, 
to show that $B_{\alpha}$ is a compact object in $L_{\alpha^{\vee}\otimes^{ct}_{k}\alpha}(X\times_{k} X)$, 
it is enough to show that it is locally, for the \'etale topology on $X\times_{k}X$, a compact object. We can therefore 
replace $X$ by the \'etale covering $U=\mathbb{R}\underline{Spec}\, A \rightarrow X$, which is affine and 
smooth over $k$. 

We thus assume that $X=\mathbb{R}\underline{Spec}\, A$ is affine and smooth over $k$, and
that $\alpha$ is realized by a smooth $A$-dg-algebra $B$. In order to finish the proof
of $(2)$ in this special case we need to show that $B$ is moreover smooth over $k$. 
For this, we write the following 
base change formula
$$B\otimes_{A}^{\mathbb{L}}B^{op}\simeq
(B\otimes_{k}^{\mathbb{L}}B^{op})\otimes_{A\otimes_{k}^{\mathbb{L}}A}^{\mathbb{L}}A.$$
As $A$ is smooth over $k$, it is a compact object in $D(A\otimes_{k}^{\mathbb{L}}A)$, and the
previous formula shows then that $B\otimes_{A}^{\mathbb{L}}B^{op}$
is compact in $D(B\otimes_{k}^{\mathbb{L}}B^{op})$, or in other words that 
the direct image functor
$$D(B\otimes_{A}^{\mathbb{L}}B^{op}) \longrightarrow D(B\otimes_{k}^{\mathbb{L}}B^{op})$$
preserves compact objects. In particular it does send $B$ itself to a compact 
object, showing that $B$ is compact in $D(B\otimes_{k}^{\mathbb{L}}B^{op})$. \\

To finish the proof of the corollary we need to prove $(3)$, $(4)$ being 
a direct consequence of $(3)$ and $(2)$. As we saw, we can write
$L_{\alpha}(X)\simeq \widehat{B_{X}^{op}}$, for $B_{X}$ the dg-algebra of 
global section of $B_{\alpha}$. By assumption $B_{\alpha}$ is a perfect complex
on $X$, and thus its global sections form a perfect $k$-dg-module because 
$X$ is proper and flat over $k$. 
\hfill $\Box$ \\

\begin{rmk}\label{rc9}
\emph{With a bit of more work, but using essentially the same arguments, we can 
extend the corollary \ref{c9} to the case where $X$ is more generally 
a seperated quasi-compact tame Deligne-Mumford derived stack, assuming that its
coarse moduli space is a scheme. We get in particular that the dg-category 
$L(X)$, of quasi-coherent complexes on $X$, is smooth (resp. proper) over $k$ when $X$ is so.
In particular, when $X$ is a smooth and proper Deligne-Mumford
derived k-stack, whose coarse moduli space is a derived scheme, 
the dg-category $L(X)$ is saturated, and therefore \cite[Thm. 3.6]{tova} can be applied in order to prove the existence of a locally algebraic derived stacks, 
locally of finite presentation over k, $\mathcal{M}_{X}$, classifiying perfect complexes on $X$. 
For instance, when $X$ is a $\mu_{n}$-gerbe over a smooth
and proper $k$-scheme $M$, 
we deduce the existence of a moduli stack of twisted sheaves and twisted perfect complexes
on $M$.}
\end{rmk}

\begin{appendix}

\section{A descent criterion}

\subsection{Relative limits}

We let $I$ be any big category (see our conventions for universes in our section on locally 
presentable dg-categories, \S 2.1) and we consider 
$\s-Cat/I$, the category of (big) $\s$-categories over 
$I$. Given $\mathcal{C} \rightarrow I$ and
$\mathcal{D} \rightarrow I$ two $\s$-categories over $I$, 
we can form the $\s$-category of relative $\s$-functors
$$\mathbb{R}\underline{Hom}_{I}(\mathcal{C},\mathcal{D}).$$
This $\s$-category is characterized, up to equivalence, by the following
universal property: for any $\s$-category $A$, we have functorial 
bijections
$$[A,\rh_{I}(\mathcal{C},\mathcal{D})]\simeq [A\times \mathcal{C},\mathcal{D}]_{I},$$
where $[-,-]$ denotes the morphisms sets in the homotopy category of 
$\s$-categories, and $[-,-]$ the morphism sets in the homotopy category 
of $\s$-categories over $I$. Using functorial fibrant replacement of 
Segal categories it is possible to define a strictly functorial
model for $\mathbb{R}\underline{Hom}_{I}(\mathcal{C},\mathcal{D})$, by considering
the Segal category of morphisms over $I$,
$$\underline{Hom}_{I}(\mathcal{C},\mathcal{D}'),$$
where $\mathcal{D}' \rightarrow I$ is a functorial fibrant replacement of the
projection $\mathcal{D} \rightarrow I$. 

We now fix $\mathcal{C} \longrightarrow I$, an $\s$-category over $I$.
Let $J$ be a small category, and
denote by $J_{+}$ the cone over $J$, obtained from $J$ by adding a new initial object
$+$ to $J$ (with no morphisms from an object of $J$ to $+$). \\

\begin{df}\label{da1}
\emph{An $\s$-category over $I$, $\mathcal{C} \longrightarrow I$,} has small
relative limits along $J$ \emph{the $\s$-functor}
$$\rh(J_{+},\mathcal{C}) \longrightarrow \rh(J,\mathcal{C})\times_{\rh(J,I)}^{h}\rh(J_{+},I)$$
\emph{possess a fully faithful right adjoint.}
\end{df}

The right adjoint of the previous definition 
will be denoted by 
$$Lim_{J/I} : \rh(J,\mathcal{C})\times_{\rh(J,I)}^{h}\rh(J_{+},I) \longrightarrow 
\rh(J_{+},\mathcal{C}),$$
and is called the \emph{relative limit $\s$-functor}. \\

Unfolding the previous definition gives the following equivalent definition: $\mathcal{C} \rightarrow I$
has relative limits along $J$ if for any commutative diagram in the $\s$-category of 
$\s$-categories
$$\xymatrix{
J \ar[r] \ar[d] & \mathcal{C} \ar[d] \\
J_{+} \ar[r] & I,}$$
there exists a lift $J_{+} \longrightarrow \mathcal{C}$, in the $\s$-category of $\s$-categories, 
which is a final object in the $\s$-category of all possible lifts. 

Suppose that $I=*$, then an $\s$-category $\mathcal{C}$, trivially considered over $*$, 
has relative limits along $J$ if and only if it has limits along $J$. Another
extreme case is when $I$ is arbitrary but $J=*$. Then it is not hard to see that 
$\mathcal{C} \longrightarrow I$ has relative limits along $*$ if and only if
it is a fibered $\s$-category in the sense of \cite[\S 1.3]{tove2}. Finally, 
if $i\in I$, and if we consider the obvious functor $J_{+} \rightarrow \{i\} \subset I$, 
we see that if $\mathcal{C} \rightarrow I$ has relative limits along $J$, then 
its fiber $\mathcal{C}_{i}$ at $i$ is an $\s$-category wich admits limits 
along $J$. These two cases essentially cover anything that can happen, as shown 
in the next lemma. 

\begin{lem}\label{la1}
Let $\mathcal{C} \rightarrow I$ be an $\s$-category over $I$. Then 
$\mathcal{C}$ has relative limits along all small categories $J$ if and only 
if it satisfies the following two conditions.
\begin{enumerate}
\item The $\s$-category $\mathcal{C}$ is fibered over $I$.
\item For any object $i\in I$, the $\s$-category $\mathcal{C}_{i}$, fiber at $i$, 
possesses all small limits. 
\end{enumerate}
\end{lem}

\textit{Proof:} We already have seen that the two conditions are necessary. Assume that they 
hold. Let $J$ be a small category, and let 
$$
\xymatrix{
J \ar[r] \ar[d] & \mathcal{C} \ar[d] \\
J_{+} \ar[r] & I,}$$
be a commutative diagram of $\s$-categories. We must show that a final objects
in the $\s$-category of lifts $J_{+} \longrightarrow \mathcal{C}$ exists. 
By pulling-back everything over $J_{+}$ we can in fact assume
that $J_{+}=I$. As $\mathcal{C}$ is fibered, it can be written 
as $\int F$, for some functor $F : I^{op} \longrightarrow \s-Cat$. The initial object $+$ of $I$ 
induces a natural transformation
$$F \longrightarrow ct(F(i_{0})),$$
from $F$ to the constant presheaf with values $F(i_{0})\simeq \mathcal{C}_{i_{0}}$. By passing to the 
Grothendiek integral we get an $\s$-functor over $I$
$$\phi : \mathcal{C} \longrightarrow \mathcal{C}_{i_{0}}\times I.$$
Precomposing with $J \rightarrow \mathcal{C}$, we get 
an $\s$-functor
$$J \longrightarrow \mathcal{C}_{i_{0}},$$
and we let $x_{0} \in\mathcal{C}_{i_{0}}$ its limit.
By construction, $x_{0}$ comes equiped with a morphism towards the
diagram $J \rightarrow \mathcal{C}_{i_{0}}$, and thus defines a 
lift
$$J_{+} \longrightarrow \mathcal{C},$$
which, as $x_{0}$ is a limit, is final among all possible lifts. \hfill $\Box$ \\

\begin{rmk}
\emph{There is a dual notion of relative colimits, that we leave to the reader. We note also 
that every we said so far about relative limits remain valid when the base category $I$
is replaced by an $\s$-category (see \cite[\S 1.3]{tove2})}. 
\end{rmk}

\subsection{The stack of stacks}

Let $T$ be a (very big) $\s$-topos in the sense of \cite[\S 1.5]{tove2}. We will typically apply the constructions and results of this
part to the case where $T$ is the $\s$-topos of big stacks over the big model 
site $(s\z-CAlg)^{op}$, of small commutative simplicial rings endowed with the fppf topology. 
According to our conventions about universes $T$ will therefore be a very big $\s$-topos. 

We consider $Mor(T):=\rh(\Delta^{1},T)$, the $\s$-category of morphisms in $T$. 
The projection $Mor(T) \longrightarrow T$, associating to every morphism its target, 
makes $Mor(T)$ into a fibered $\s$-category over $T$. It corresponds to
an $\s$-functor
$$\mathbb{S} : T^{op} \longrightarrow \s-\underline{Cat},$$
whose value at an object $x \in T$ is the comma $\s$-category $T/x$. The
$\s$-functor $\mathbb{S}$ is called the \emph{universal family} over $T$. Its main 
property is the following exactness condition.

\begin{lem}\label{la2}
The $\s$-functor $\mathbb{S}$ commutes with limits (i.e. sends
colimits in $T$ to limits in $\s-\underline{Cat}$).
\end{lem}

\textit{Proof:} It is enough to show that the $\mathbb{S}$ commutes with 
(possibly infinite) products, as well as with fibered products. In other words, 
we need to show the following two statements.

\begin{enumerate}
\item Given a (small) family of objects $x_{\alpha}\in T$, $\alpha\in A$, 
with sum $x:=\sum x_{\alpha}$, the base change $\s$-functor
$$\prod_{\alpha}(-\times_{x}x_{\alpha}) : T/x \longrightarrow 
\prod_{\alpha}T/x_{\alpha}$$
is an equivalence of $\s$-categories. 
\item Given a cocartesian diagram of objects in $T$
$$\xymatrix{
t \ar[r] \ar[d] & y \ar[d] \\
x \ar[r] & s}$$
the induced diagram of $\s$-categories
$$\xymatrix{
T/s \ar[r]^-{-\times_{s}x} \ar[d]_-{-\times_{s}y} & T/x \ar[d]^-{-\times_{x}z} \\
T/y \ar [r]_-{-\times_{y}z} & T/z}$$
is (homotopy) cartesian.
\end{enumerate}

Condition $(1)$ above follows from the fact that sums are disjoints in $T$, and
that colimits are stable by base change. Indeed, 
an inverse $\s$-functor is given by sending a family of objects $y_{\alpha} \rightarrow x_{\alpha}$ 
to its sum $\coprod y_{\alpha} \rightarrow x$. 

In order to prove condition $(2)$, we consider the induced $\s$-functor
$$\phi : T/s \longrightarrow T/x\times_{T/z}^{h}T/y.$$
This $\s$-functor possesses a left adjoint, constructed as follows. We consider
the diagram of $\s$-categories
$$\xymatrix{
 T/s \ar[r]^-{-\times_{s}x} \ar[d]_-{-\times_{s}y} & T/x \ar[d]^-{-\times_{x}z} \\
T/y \ar [r]_-{-\times_{y}z} & T/z,}$$
which we consider as a functor 
$$F : I^{op} \longrightarrow \s-Cat,$$
where, obviously $I=\Delta^{1}\times \Delta^{1}$ is the category classifying 
commutative squares. We perform its Grothendieck's construction to get a fibered $\s$-category
$$\mathcal{C}:=\int F \longrightarrow I.$$
We note that $I$ is also the cocone over the category $J=\xymatrix{1 & \ar[r] \ar[l] 0 & 2}$
$$I = \xymatrix{
0 \ar[r] \ar[d] & 2 \ar[d] \\
1 \ar[r] & +.}$$
The $\s$-category $\mathcal{C}$ is also fibered, as the base change $\s$-functors possesses 
left adjoints given by the forgetful $\s$-functors: for $u : x \rightarrow y$ in $T$, the
left adjoint to $-\times_{x}y$ is $T/x \rightarrow T/y$ simply obtained by composition
with $u$. 

There is a natural equivalence of $\s$-categories
$$T/x\times_{T/z}^{h}T/y \simeq \mathbb{R}\underline{Hom}_{I}^{cart}(J,\mathcal{C}),$$
where the right hand side is the $\s$-category of cartesian sections of 
$\mathcal{C} \rightarrow I$ over the sub-category $J\subset I$ (this formally follows
from \cite[Prop. 1.4]{tove2}). Similarly, we have a natural equivalence of $\s$-categories
$$T/s \simeq \mathbb{R}\underline{Hom}_{I}^{cart}(I,\mathcal{C}),$$
and thus a natural inclusion $\s$-functor $T/s \subset \mathbb{R}\underline{Hom}_{I}(I,\mathcal{C})$.
This inclusion $\s$-functor is right adjoint to the evaluation at $+$
$$\mathbb{R}\underline{Hom}_{I}(I,\mathcal{C}) \longrightarrow T/s.$$
Using these identifications, the $\s$-functor $\phi$ becomes the restriction $\s$-functor
$$\mathbb{R}\underline{Hom}_{I}^{cart}(I,\mathcal{C}) \longrightarrow 
\mathbb{R}\underline{Hom}_{I}^{cart}(J,\mathcal{C}).$$
But, as $I=J_{+}$, and as $\mathcal{C}$ possesses relative colimits (see lemma \ref{la1}), this
restriction $\s$-functor possesses a left adjoint 
$$\xymatrix{
\mathbb{R}\underline{Hom}_{I}^{cart}(J,\mathcal{C}) \ar^-{Colim_{J/I}} &
\mathbb{R}\underline{Hom}_{I}(J,\mathcal{C}) \ar[r]^-{ev_{+}} & T/s.}$$
This provides a left adjoint 
$$\psi : T/x\times_{T/z}^{h}T/y \longrightarrow T/s$$
to the $\s$-functor $\phi$. 

The unit and counit of the adjunction $(\psi,\phi)$ are explicitely given as follows. 
Let $u \rightarrow s$ be an object in $T/s$. Then the counit morphism in $T/s$ is
$$(u\times_{s}x)\coprod_{(u\times_{s}z)}(u\times_{s}y) \longrightarrow u.$$
That this morphism is an equivalence in the $\s$-category $T/s$ is simply the 
compatibility of colimits with base change in $T$
$$(u\times_{s}x)\coprod_{(u\times_{s}z)}(u\times_{s}y) \simeq u\times_{s}(x\coprod_{s}y)
\simeq u\times_{s}s\simeq u.$$
The unit of the adjunction is itself given as follows. An object in $T/x\times_{T/z}^{h}T/y$
is represented by a triple 
$$(x' \rightarrow x) \qquad (y' \rightarrow y) \qquad \alpha : x'\times_{x}z \simeq y'\times_{y}z,$$
consisting of two objects $x'\in T/x$, $y'\in T/y$, and an equivalence $\alpha$ between the two 
base changes in $T/z$. The image by $\psi$ of this object is
$$x'\coprod_{x'\times_{x}z}y' \rightarrow s,$$
where $\alpha$ is used in order to define the morphism
$x'\times_{x}z \simeq y'\times_{y}z \rightarrow y'$. 
Therefore, to prove that the unit of the adjunction is an equivalence we need to show that 
both morphisms
$$x'\longrightarrow (x'\coprod_{x'\times_{x}z}y')\times_{s}x  \qquad
y'\longrightarrow (x'\coprod_{x'\times_{x}z}y')\times_{s}y$$
are equivalences in $T$. For this, we write $T$ as an exact localisation of 
an $\s$-category of the form $\widehat{T_{0}}$, of very big prestacks over a big
$\s$-category $T_{0}$. Then, to check that the above two morphisms are always
equivalences in $T$ it is enough to do this replacing $T$ by $\widehat{T_{0}}$. 
Moreover, by evaluation at each object of $T_{0}$ we can further assume that 
$T_{0}=*$. In other words, we can assume that $T$ is simply the $\s$-topos of 
(very big) Kan complexes. But in this case, that the above two morphisms are 
weak equivalences in simply a form of the Van Kampen theorem (see for instance \cite[Prop. 8.1]{rez}). 
\hfill $\Box$ \\

\begin{rmk}\label{ra1}
\emph{The lemma \ref{la2} provides a characterization of $\s$-topos. Indeed, 
a locally presentable $\s$-category $T$ is an $\s$-topos if and only if 
$x \mapsto T/x$ transforms colimits in $T$ into limits in $\s$-cat. This is known as
the} distributive law \emph{(see e.g. \cite[Thm. 6.1.0.6]{lu1})}. 
\end{rmk}

An interesting consequence of the previous lemma is the following construction. 
Suppose that $x_{*} : I \longrightarrow T$ is a (small) diagram of objects in $T$, 
together with a colimit $x_{0}=colim_{i\in I}x_{i}$. Then, the pull-back $\s$-functor
$$T/x_{0} \longrightarrow Holim_{i\in I}T/x_{i}$$
is an equivalence of $\s$-categories. There are two different ways of writing an 
inverse $\s$-functor, one using left adjoints of the base changes functors, as we 
have done during the proof of the lemma, and a second one using the right adjoints
to the base changes. Let us describe now this second construction. For a given 
morphism $p : x \rightarrow y$ in $T$, the base change
$$T/y \longrightarrow T/x$$
possesses a right adjoint 
$$p_{*} : T/x \longrightarrow T/y,$$
which is given by the relative internal Hom object. 

We consider $x \mapsto T/x$ and $p \mapsto p_{*}$ as a 
fibered $\s$-category $\mathcal{C} \longrightarrow T^{op}$. This fibered
$\s$-category possesses relative limits according to lemma \ref{la1}. 
We 
have 
$$Holim_{i\in I}T/x_{i}\simeq \mathbb{R}\underline{Hom}_{T^{op}}^{cart}(I^{op},\mathcal{C}).$$
Therefore, using existence of relative limits along $I$ we produce an $\s$-functor
$$Lim_{I/T} : \mathbb{R}\underline{Hom}_{T}^{cart}(I^{op},\mathcal{C}) 
\longrightarrow \mathbb{R}\underline{Hom}_{T}(I^{op}_{+},\mathcal{C}),$$
where $I_{+} \longrightarrow T$ is fixed using the colimit object $x_{0}$ as the value
at the initial point $+\in I_{+}$. Composing with the evaluation at $+$ we get 
a well defined $\s$-functor
$$Holim_{i\in I}T/x_{i} \longrightarrow T/x_{0}$$
which is a right adjoint to the equivalence $T/x_{0} \simeq Holim T/x_{i}$. This 
right adjoint is therefore also an equivalence of $\s$-categories. It can be described 
intuitively as follows. Consider an object $x'$ in $Holim T/x_{i}$. Such an object is
represented by a system of objects $x_{i}' \rightarrow x_{i}$, together 
with natural equivalences $x_{i}'\times_{x_{i}}x_{j} \simeq x'_{j}$, for 
all $i\rightarrow j$ a morphism in $I$. For each $i$ we have the natural 
morphis $p : x_{i} \rightarrow x_{0}$, and we consider
$p_{*}(x_{i})$. By adjunction, we get that way an $I$-diagram of objects 
$i \mapsto p_{*}(x_{i})$ inside the $\infty$-category $T/x_{0}$. 
The limit of this diagram is the image of $x'$ in $T/x_{0}$. In other words, 
the $\s$-functor $Holim_{i\in I}T/x_{i} \longrightarrow T/x_{0}$ sends $x'$ to
$$lim_{i\in I}p_{*}(x_{i}) \in T/x_{0}.$$
As this $\s$-functor is an equivalence we get that for any $y \rightarrow x_{0}$, 
the natural morphism
$$y \rightarrow lim p_{*}(y\times_{x_{0}}x_{i})$$
is an equivalence in $T$. In the same way, for any $x'$ in $Holim T/x_{i}$, we
the natural morphisms
$$x_{i}' \rightarrow x_{i}'\times_{x_{0}}(lim p_{*}(x'_{i})$$
are all equivalences in $T$. 

\subsection{The descent statement}

We now use our last two paragraph to state and prove our descent statement. It consists of 
a criterion to show that a prestack in $\s$-categories is a stack when 
it comes equiped with good forgetful functors to the stack of stacks $\mathbb{S}$ 
defined our last paragraph. \\

We let $T$ be a very big $\s$-topos, generated by a big $\s$-site $(C,\tau)$. 
We let $j : C \longrightarrow T$ be the natural embedding. The universal family
$\mathbb{S}$, restricted to $C$ is still denoted by $\mathbb{S}$, and is
called the \emph{stack of stacks}. It can be described explicitely as follows. For
each object $U\in C$, $\mathbb{S}(U)$ is the $\infty$-category of stacks over 
$C/U$, the comma $\s$-category endowed with the induced topology. For 
$u : U \rightarrow V$ in $C$, the base change
$$u^{*} : \mathbb{S}(V) \longrightarrow \mathbb{S}(U)$$
is simply obtained by the restriction along $u\circ : C/U \longrightarrow C/V$.
Recall that all the $\s$-functors $u^{*}$ possess right adjoint $u_{*}$.

\begin{prop}\label{pa2}
With the same notations as above, let 
$$F : C^{op} \longrightarrow \s-\underline{Cat}_{\mathbb{W}}$$
be a prestack of $\s$-categories. We assume that the following conditions hold.
\begin{enumerate}
\item For each object $x\in T$ the $\s$-category $F(x)$ has 
small limits. 
\item For each morphism $u : x \rightarrow y$, the $\s$-functor
$$u^{*} : F(y) \longrightarrow F(x)$$
has right adjoint $u_{*}$. 
\item There exists a family of morphism (in the $\s$-category of prestacks of $\s$-categories)
$$\phi_{\alpha} : F \longrightarrow \mathbb{S},$$
such that for each $x\in T$ the corresponding family of $\s$-functors
$$\{\phi_{\alpha}(x) : F(x) \longrightarrow \mathbb{S}(x)\}_{\alpha}$$
is conservative.
\item For each $x\in T$, and all $\alpha$, the $\s$-functor
$$\phi_{\alpha}(x) : F(x) \longrightarrow \mathbb{S}(x)$$
commutes with limits.
\item For each morphism $u : x \rightarrow y$ in $T$, and for all $\alpha$, the natural 
transformation
$$\phi_{\alpha}\circ u_{*} \Rightarrow u_{*}\phi_{\alpha}$$
is an equivalence (in the $\s$-category of $\s$-functors from $F(x)$ to 
$\mathbb{S}(y)$). 
\end{enumerate}
Then $F$ is a stack (i.e. has the descent property of 
all $\tau$-hypercoverings).
\end{prop}

\textit{Proof:} Few words of explanations concerning condition $(5)$ first. 
We have a natural equivalence of $\s$-functors $u^{*}\circ \phi_{\alpha} \simeq \phi_{\alpha}\circ^{*}$, 
coming from the fact that the $\phi_{\alpha}$ are natural transformations. By adjunction these
natural equivalences provide natural transformations
$$\phi_{\alpha}\circ u_{*} \Rightarrow u_{*}\phi_{\alpha},$$
which by condition $(5)$ must be equivalences. \\

Let $H_{*} \rightarrow X$ be a hyper-coverings over $X \in C$. We must show that 
the natural $\s$-functor
$$f : F(X) \longrightarrow Holim_{n\in \Delta} F(H_{n})$$
is an equivalence of $\s$-categories. For this, we use the morphisms $\phi_{\alpha}$, in order to get 
commutative diagrams
$$\xymatrix{
F(X) \ar[r]^-{f} \ar[d]_-{\phi_{\alpha}} & Holim_{n\in \Delta} F(H_{n}) \ar[d]^-{\alpha} \\
\mathbb{S}(X) \ar[r]_-{f'} & Holim \mathbb{S}(H_{n}).}$$ 
Because of lemma \ref{la2} we know that $\mathbb{S}$ is a stack, and therefore that 
the bottom horizontal $\s$-functor is an equivalence. Using relative limits, we can construct, 
for $F$ and $\mathbb{S}$, right adjoints $g$ and $g'$ to the $\s$-functors $f$ and $g$
(as we have done for $\mathbb{S}$ in our last paragraph), and because of conditions $(4)$ and $(5)$
we get adjoints commutative diagrams
$$\xymatrix{
Holim_{n\in \Delta} F(H_{n})  \ar[r]^-{g} \ar[d]_-{\phi_{\alpha}} &  \ar[d]^-{\alpha}F(X) \\
Holim \mathbb{S}(H_{n}) \ar[r]_-{g'} & \mathbb{S}(X).}$$ 
Therefore, to prove that $F$ is stack we simply need to prove that the unit and counit of
the adjunction $(f,g)$ are equivalences. But, as the $\phi_{\alpha}$ form a conservative family this
follows from the corresponding statement for $\mathbb{S}$, which we already know is true
because $f'$ is an equivalence by lemma \ref{la2}.
\hfill $\Box$ \\

A nice consequence of the previous proposition is the following descent criterion. 
We state it using model sites as this is the
version that we will use. It is however true for any $\s$-site. 

\begin{cor}\label{ca2}
Let $(C,\tau)$ be a big $\s$-site and 
$$F : C^{op} \longrightarrow \s-\underline{Cat}_{\mathbb{W}}$$
be a prestack in very big $\s$-categories. 
We assume that the following conditions hold.
\begin{enumerate}
\item For any object $U\in C$, the $\s$-category $F(U)$ has all small limits, and for any 
$u : V \rightarrow U$ the base change $\s$-functor
$$u^{*} : F(U) \longrightarrow F(V)$$
commutes with all small limits. 
\item For any $u : V \rightarrow U$ in $C$, the base change $\s$-functor $u^{*}$ possesses
a right adjoint 
$$u_{*} : F(U) \longrightarrow F(V),$$
which is moreover conservative.
\item For each cartesian square in the $\s$-category $C$
$$\xymatrix{
V' \ar[r]^-{u'} \ar[d]_-{q'} & U' \ar[d]^-{q} \\
V \ar[r]_-{u} & U,}$$
and any object $x\in F(V)$, the natural morphism
$$q^*u_{*}(x) \longrightarrow q'_{*}(u')^{*}(x)$$
is an equivalence in $F(U')$. 
\item For each object $U\in C$, and 
each two objects $x,y \in F(U)$, the simplicial presheaf
$$Map(x,y) : (C/U)^{op} \longrightarrow SSet_{\mathbb{W}},$$
sending $u : V \rightarrow U$ to the simplicial set $Map_{F(V)}(u^{*}(x),u^{*}(y))$,
is a stack for the induced model topology.
\end{enumerate}
Then $F$ is a stack. 
\end{cor}

\textit{Proof:} 
We let $A$ be the set of equivalence classes of objects in the $\s$-category 
$F(*)$. For any $x\in A$, we consider the morphism
$$\phi_{x} : F \longrightarrow \mathbb{S}$$
defines as follows. For each object $U\in C$, we let 
$x_{U}:=p^{*}(x)$, where $p : U \longrightarrow *$ is the natural projection. 
For any $y\in F(U)$ we set 
$$\phi_{x}(y):=Map(x_{U},y),$$
where $Map(x_{U},y)$ is the prestack over $C/U$ defined by sending
$u : V\rightarrow U$ to the simplicial set 
$Map_{F(V)}(x_{V},u^{*}(y))$. By the condition $(4)$, the prestack 
$\phi_{x}(y)$ is a stack over $C/U$. Moreover, we obviously have, 
for any $u : V \rightarrow U$ and $y\in F(U)$
$$\phi_{x}(u^{*}(y))=u^{*}(\phi_{x}(y)),$$
showing that $\phi_{x}$ defines a morphism of prestacks
$$\phi_{x} : F \longrightarrow \mathbb{S}.$$
When $x$ varies in $A$ we get our family of morphisms
$$\{\phi_{x} : F \longrightarrow \mathbb{S}\}_{x\in A}.$$
It is then formal to check that this family satisfies the properties
of the propostion \ref{pa2}. \hfill $\Box$ \\

\end{appendix}

\end{document}